\documentclass[11pt,oneside,a4paper]{article}
\usepackage{blindtext}
\usepackage{amsmath,amsfonts,amssymb,amsthm,eucal}
\usepackage{microtype}
\usepackage[colorlinks, citecolor=blue, linkcolor=blue, urlcolor=Maroon, filecolor=Maroon, linktocpage=true]{hyperref}
\usepackage[usenames,dvipsnames,svgnames,table]{xcolor}
\usepackage[normalem]{ulem}
\usepackage{lmodern}
\usepackage{enumerate}
\usepackage{bbm}
\usepackage{tikz}
\usepackage{tikz-cd}
\usepackage[font=footnotesize,labelfont=bf]{caption}
\usepackage{subcaption}
\usepackage{framed}
\usepackage[nottoc,notlot,notlof]{tocbibind}
\usepackage{tocloft}
\usepackage{morefloats}
\usepackage{float}
\usepackage[nottoc,notlot,notlof]{tocbibind}
\usepackage{tocloft}
\usepackage{hhline}
\DeclareMathOperator{\D}{d\!}
\theoremstyle{definition}
\newtheorem{Th}{Theorem}[section]

\newtheorem{Prop}[Th]{Proposition}

\newtheorem{Cor}[Th]{Corollary}
\newtheorem{Lem}[Th]{Lemma}
\newtheorem{Def}[Th]{Definition}

\theoremstyle{definition}
\newtheorem{rem}[Th]{Remark}
\usepackage[a4paper,
	hcentering,
	text={453pt,686pt},
]{geometry}
\allowdisplaybreaks
\catcode`@=11
\@addtoreset{equation}{section}
\catcode`@=12
\begin{document}

%

%
%
%\rightline
%{\textcolor{red}{\today}}
%\vskip 1.0 true cm  
\begin{center}\noindent
\textbf{\Large Special geometry of quartic curves}
\\[2em]
{\fontseries{m}\fontfamily{cmss}\selectfont \large David\ Lindemann}\\[1em] 
{\small
Department of Mathematics, Aarhus University\\
Ny Munkegade 118, Bldg 1530, DK-8000 Aarhus C, Denmark\\
\texttt{david.lindemann@math.au.dk}
}
\end{center}
\vspace{1em}
\begin{abstract}
\noindent
We classify maximal quartic generalised projective special real curves up to equivalence. A maximal quartic generalised projective special real curve consists of connected components of the intersection of the hyperbolic points of a quartic homogeneous real polynomial $h:\mathbb{R}^2\to\mathbb{R}$ and its level set $\{h=1\}$. Two such curves are called equivalent if they are related by a linear coordinate transformation. As an application of our results we prove that quartic generalised projective special real manifolds that are homogeneous spaces have non-regular boundary behaviour, meaning that the differential of each of these spaces' defining polynomials vanishes identically on a ray in the boundary of the cone spanned by the corresponding manifold. Lastly we describe the asymptotic behaviour of each curve.
\end{abstract}
\textbf{Keywords:} affine differential geometry, centro-affine curves, special real geometry, generalised special real geometry, quartic polynomials\\
\textbf{MSC classification:} 53A15, 51N35 (primary), 51N10, 53C26 (secondary)
\tableofcontents
\section{Introduction}
	In this work we study \textbf{maximal quartic generalised projective special real curves} and determine their equivalence classes up to linear transformations of the ambient space, see Theorems \ref{thm_quartic_CCGPSR_curves_classification} and \ref{thm_incomplete_quartics}. Furthermore we describe their asymptotic behaviour in Proposition \ref{prop_limit_geos}. In order to proceed and make sense of what is to come, we first need to clear up some terminology.
	\begin{Def}\label{def_equivalence_polys}
		A homogeneous polynomial $h:\mathbb{R}^{n+1}\to\mathbb{R}$, $n\geq 1$, of degree $\tau\geq 2$ is called \textbf{hyperbolic} if there exists $p\in\mathbb{R}^{n+1}$, such that $-\partial^2h_p$ has Lorentzian signature. Such a point $p$ is called \textbf{hyperbolic point} of $h$. Two hyperbolic homogeneous polynomials $h,\overline{h}:\mathbb{R}^{n+1}\to\mathbb{R}$ of the same degree are called \textbf{equivalent} if they are related by a linear transformation of $\mathbb{R}^{n+1}$.
	\end{Def}
	
	\begin{Def}
		A \textbf{projective special real} (short: \textbf{PSR}) \textbf{manifold} of dimension $n\geq 1$ is a hypersurface $\mathcal{H}\subset\mathbb{R}^{n+1}$ that is contained in the level set $\{h=1\}$ of a hyperbolic homogeneous cubic polynomial $h:\mathbb{R}^{n+1}\to\mathbb{R}$, such that $\mathcal{H}$ consists only of hyperbolic points of $h$. Similarly, \textbf{generalised projective special real} (short: \textbf{GPSR}) \textbf{manifold} of dimension $n\geq 1$ is a hypersurface $\mathcal{H}\subset\mathbb{R}^{n+1}$ that is contained in the level set $\{h=1\}$ of a hyperbolic homogeneous polynomial $h:\mathbb{R}^{n+1}\to\mathbb{R}$ of degree $\tau\geq 4$ and consists only of hyperbolic points of $h$.
	\end{Def}
	
	Whenever we write ``\textbf{(G)PSR} manifold'' we mean either a PSR manifold, corresponding to $\tau=3$, or a GPSR manifold for arbitrary but fixed $\tau\geq 4$. In the case of $\tau =4$ we will call the corresponding GPSR manifolds \textbf{quartic GPSR manifolds}. Similarly as for their defining polynomials, there is a notion of equivalence for (G)PSR manifolds.
	
	\begin{Def}
		Two (G)PSR manifolds $\mathcal{H}\subset\{h=1\}$, $\overline{\mathcal{H}}\subset\{\overline{h}=1\}$, of the same dimension $n$ with defining polynomials of coinciding degree are called \textbf{equivalent} if there exists a linear transformation of the ambient space $A\in\mathrm{GL}(n+1)$, such that $A(\mathcal{H})=\overline{\mathcal{H}}$.
	\end{Def}
	
	Note that $\mathcal{H}\subset\{h=1\}$ and $\overline{\mathcal{H}}\subset\{\overline{h}=1\}$ being equivalent implies that their defining polynomials $h$ and $\overline{h}$ are also equivalent in the sense of Definition \ref{def_equivalence_polys}, but the converse need not hold. For example, for any open subset $U\subset\mathcal{H}$, such that $U\ne\mathcal{H}$, $U$ and $\mathcal{H}$ are not equivalent but have the same defining polynomial.
	
	\begin{Def}\label{def_maximal_closed}
		A (G)PSR manifold $\mathcal{H}\subset\{h=1\}$ is called
			\begin{enumerate}[(i)]
				\item \textbf{maximal} if it consists of connected components of $\{h=1\}\cap\{\text{hyperbolic points of }h\}$,\label{def_maximal}
				\item \textbf{closed} if it is closed as a subspace of its ambient space.\label{def_closed}
			\end{enumerate}
		If a maximal (G)PSR manifold is connected, we say that it is a \textbf{maximal connected} (G)PSR manifold. Note that $\mathcal{H}$ being closed automatically implies that it is maximal. The other direction need not hold. For concrete examples for the latter see Theorem \ref{thm_incomplete_quartics}.
	\end{Def}
	In the following, we will abbreviate ``closed connected PSR, GPSR, (G)PSR'' with \textbf{CCPSR}, \textbf{CCGPSR}, \textbf{CC(G)PSR}, respectively. Every (G)PSR manifold $\mathcal{H}\subset\{h=1\}$ is a \textbf{centro-affine hypersurface} as the position vector field is automatically transversal along $\mathcal{H}$, which follows from Euler's homogeneous function theorem. The corresponding \textbf{centro-affine fundamental form} $g_{\mathcal{H}}$ is given by
		\begin{equation}
			g_\mathcal{H}=-\frac{1}{\tau}\partial^2 h|_{T\mathcal{H}\times T\mathcal{H}},
		\end{equation}
	which follows from \cite[Prop.\,1.3]{CNS}. The centro-affine fundamental form is determined by the \textbf{centro-affine Gau{\ss} equation}
		\begin{equation*}
			\overline{\nabla}_X Y = \nabla^{\mathrm{ca}}_X Y + g_{\mathcal{H}}(X,Y)\xi
		\end{equation*}
	for all $X,Y\in\mathfrak{X}(\mathcal{H})$, where $\overline{\nabla}$ denotes the flat connection on the ambient space $\mathbb{R}^{n+1}$, $\nabla^{\mathrm{ca}}$ denotes the induced \textbf{centro-affine connection}, and $\xi$ denotes the position vector field in $\mathbb{R}^{n+1}$. Note that in general the centro-affine connection and the Levi-Civita connection of $(\mathcal{H},g_{\mathcal{H}})$ do not coincide. By Euler's homogeneous function theorem and the hyperbolicity of each point in a (G)PSR manifold it follows that the centro-affine fundamental form of a (G)PSR manifold is a Riemannian metric.

	Thus, when we talk about maximal connected quartic GPSR curves, we mean a connected component of a hyperbolic level set of a homogeneous quartic polynomial in two real variables. Before introducing more advanced technical tools in Section \ref{sect_prelim} that we will use in order to prove our main results in Theorems \ref{thm_quartic_CCGPSR_curves_classification} and \ref{thm_incomplete_quartics}, we should ask ourselves why one should be interested in the geometry of such curves and, more generally, quartic GPSR manifolds of higher dimension. There are three main motivations.
	
	Firstly, the field of \textbf{special geometry}. Under the term ``special geometry'' we understand the study of \textbf{projective special K\"ahler (PSK) manifolds} \cite{F} and \textbf{PSR manifolds}. This field is motivated by theoretical physics, more specifically supergravity theories. PSR manifolds are related to PSK manifolds by the supergravity r-map \cite{CHM} and correspond to dimensional reduction from $5$ to $4$ spacetime dimensions \cite{GST,DV}. This allowed a very explicit description of the so-obtained PSK manifolds, since their geometry is completely determined by the underlying PSR manifold. There is currently no known way to generalise the supergravity r-map to GPSR manifolds with defining polynomial of degree $\tau\geq 4$. One ansatz to obtain a possible generalisation is to try to construct such a generalisation in the for the least complicated positive-dimensional GPSR manifolds, which precisely are quartic GPSR curves. Our present work will thus be of good use for future research in that direction.
	
	Secondly, (G)PSR manifolds appear naturally in the \textbf{geometry of K\"ahler cones} \cite{DP,W} as level sets of the intersection forms in the real $(1,1)$-cohomology. In \cite{W} the curvature and examples of so-obtainable PSR manifolds have been studied. Further insight in this topic from an algebraic geometry perspective has been gained in \cite{M1,M2}.
	
	Lastly, the classification of real polynomials in degree $\tau\geq 3$ is a long-standing classical problem. The study of cubics goes back to Newton \cite{N}. Even when restricting to hyperbolic homogeneous polynomials, there is only a general well-known classification of linear and quadratic real homogeneous polynomials, but even for cubics no general classification is known. For cubics there are partial results, see \cite{CHM} for the case of $2$ and \cite{CDL} for the case of $3$ real variables. Further classification results have been obtained in \cite{DV} where the authors classify all hyperbolic homogeneous cubics $h$ in any number of variables that have a homogeneous space as one of the connected components of their level sets $\{h=1\}$. In \cite{CDJL} all reducible hyperbolic homogeneous cubics are classified, and in \cite{L2} it has been shown that the moduli space of hyperbolic homogeneous cubics $h$ that admit a closed connected component of $\{h=1\}$ is compactly generated and connected. For degree $4$ and up there are no known classification results for hyperbolic homogeneous polynomials in at leas $2$ real variables. This work presents the first such results for quartics in $2$ real variables. There are, however, important related results for real quartic curves. In \cite{KW} the isotopy types of affine and projective quartic curves have been classified, and the roots of real quartic polynomials in $1$ variable have been extensively studied in \cite{P}. In \cite[Sect.\,4]{T} an explicit example of a hyperbolic homogeneous quartic polynomial $h:\mathbb{R}^3\to\mathbb{R}$ with regards to its level set $\{h=1\}$ and the curvature induced by its centro-affine fundamental form is studied. This is explicitly related to the results of \cite{W} in the cubic case.
	
	\paragraph*{Acknowledgements} Theorem \ref{thm_quartic_CCGPSR_curves_classification} is based on the authors PhD-thesis \cite{L1}, albeit the proof has been streamlined and all figures have been redone. The other results of this work are new. This work was partly supported by the \textit{German Research Foundation} (DFG) under the RTG 1670 ``Mathematics Inspired by String Theory and Quantum Field Theory'', and partly by a ``Walter Benjamin PostDoc Fellowship'' granted to the author by the DFG. The author would like to thank Vicente Cort\'es and Andrew Swann for helpful discussion during the preparation of this work.

\section{Preliminaries}\label{sect_prelim}
	\begin{Lem}\label{lem_CC_implies_precompact}
		Let $\mathcal{H}$ be a CC(G)PSR manifold and let $p\in\mathcal{H}$ be arbitrary. Then the intersection of the cone spanned by $\mathcal{H}$ and the affinely embedded tangent space of $\mathcal{H}$ at $p$, that is
			\begin{equation*}
				\left(\mathbb{R}_{>0}\cdot \mathcal{H}\right)\cap \left(p + T_p\mathcal{H}\right)\subset\mathbb{R}^{n+1},
			\end{equation*}
		is precompact.
		\begin{proof}
			\cite[Lem.\,1.14]{CNS}
		\end{proof}
	\end{Lem}
	
	\begin{rem}
		A natural question arising from Lemma \ref{lem_CC_implies_precompact} is the following. Suppose we are given a hyperbolic homogeneous polynomial $h:\mathbb{R}^{n+1}\to\mathbb{R}$ of degree $\tau\geq 2$ and a hyperbolic point $p$ of $h$ that is contained in a fixed connected component $\mathcal{H}$ of $\{h=1\}$. Further suppose that the cone spanned by $\mathcal{H}$ intersected with the affinely embedded tangent space of $\mathcal{H}$ at $p$ is precompact in the ambient space. Does this already imply that every point in $\mathcal{H}$ is a hyperbolic point of $h$? For $\tau=2$ this is easily seen to be true as there exists precisely one hyperbolic quadratic polynomial up to equivalence in each dimension and each connected component of its positive level set is isomorphic to the homogeneous space $\mathrm{SO}(n,1)^+/\mathrm{SO}(n)$. For $\tau\geq 4$, which corresponds to GPSR manifolds, this does however in general not hold true. For a counterexample in the quartic case we have Thm. \ref{thm_incomplete_quartics} \hyperref[eqn_incomplete_qGPSRcurves_class_c]{c)}, for which every connected component of $\{h=1\}\cap\{\text{hyperbolic points of }h\}$ has precisely one flat point. Surprisingly, this statement does in fact hold true for $\tau=3$, that is for PSR manifolds. This is the main result of \cite{L2}. In \cite[Thm.\,1.1]{L2} a compact convex generating set $\mathcal{C}_n$ for the moduli space of $n$-dimensional CCPSR manifolds is constructed, given by the affine subset of cubic polynomials in $n+1$ variables
			\begin{equation}\label{eqn_CCPSR_gen_set}
				\mathcal{C}_n=\left\{x^3-x\langle y,y\rangle + P_3(y)\ \left|\ \max\limits_{\|y\|=1} P_3(y)\leq\frac{2}{3\sqrt{3}}\right.\right\},
			\end{equation}
		where $y=(y_1,\ldots,y_n)^T$ denotes linear coordinates on $\mathbb{R}^n$, $\langle\cdot,\cdot\rangle$ is the induced Euclidean scalar product with norm $\|\cdot\|$, and $P_3:\mathbb{R}^n\to\mathbb{R}$ is a cubic homogeneous polynomial. The statement of \cite[Thm.\,1.1]{L2} is that every $n$-dimensional CCPSR manifold $\mathcal{H}$ is, up to equivalence, the connected component of a level set $\{h=1\}$ that contains the point $\left(\begin{smallmatrix}x\\y\end{smallmatrix}\right)=\left(\begin{smallmatrix}1\\0\end{smallmatrix}\right)$ for some $h\in\mathcal{C}_n$. One can now show that the condition for the cone spanned by said connected component intersected with $\left(\begin{smallmatrix}1\\0\end{smallmatrix}\right)+T_{\left(\begin{smallmatrix}1\\0\end{smallmatrix}\right)}\mathcal{H}$ is precompact is precisely the maximality condition for $P_3$ in \eqref{eqn_CCPSR_gen_set}.
	\end{rem}

	\begin{Def}
		In the following we will denote by $G^h$ the \textbf{linear automorphism group} of a given hyperbolic homogeneous polynomial $h$.
	\end{Def}
	
	Note that the linear automorphism group $G^h$ need not coincide with the linear automorphism group of a given (G)PSR manifold $\mathcal{H}\subset\mathbb{R}^{n+1}$ contained in $\{h=1\}$, which is given by all linear transformations $A\in\mathrm{GL}(n+1)$, such that $A(\mathcal{H})=\mathcal{H}$. Similarly as for the generating set of the moduli space of CCPSR manifolds discussed above which allows to assume that the defining polynomial of a CCPSR manifold is of a certain form, there is an analogue statement that holds for all (G)PSR manifolds.
	
	\begin{Prop}\label{prop_standard_form}
		Let $\mathcal{H}$ be a (G)PSR manifold contained in the level set of a hyperbolic homogeneous polynomial $h:\mathbb{R}^{n+1}\to\mathbb{R}$ of degree $\tau\geq 3$, and let further $\left(\begin{smallmatrix}x\\y\end{smallmatrix}\right)=(x,y_1,\ldots,y_n)^T$ denote linear coordinates on $\mathbb{R}^{n+1}$ with induced Euclidean scalar product on $\{0\}\times\mathbb{R}^n$ denotes by $\langle\cdot,\cdot\rangle$. Then for all $p\in \mathcal{H}$ there exist a linear transformation of the ambient space $A(p)\in\mathrm{GL}(n+1)$, such that
			\begin{enumerate}[(i)]
				\item $A(p)^*h=x^\tau-x^{\tau-2}\langle y,y\rangle + \sum\limits_{k=3}^\tau x^{\tau-k}P_k(y)$,\label{prop_standard_form_eqn_i}
				\item $A(p)\cdot\left(\begin{smallmatrix}1\\0\end{smallmatrix}\right)=p$,\label{prop_standard_form_eqn_ii}
			\end{enumerate}
		where for each $3\leq k\leq \tau$, $P_k:\mathbb{R}^n\to\mathbb{R}$ is a homogeneous polynomial of degree $k$. Locally, $A$ can be chosen to be a smooth map. If $\mathcal{H}$ is closed and connected, one can choose $A:\mathcal{H}\to\mathrm{GL}(n+1)$ to be smooth globally.
		\begin{proof}
			\cite[Prop.\,3.1]{L2}
		\end{proof}
	\end{Prop}
	
	In general, the polynomials $P_k$ in Proposition \ref{prop_standard_form} \eqref{prop_standard_form_eqn_i} are not uniquely determined by the manifold $\mathcal{H}$ alone without any additional assumptions, for the cubic case cf. the discussion in the last paragraph in \cite[p.\,10]{L2}. Proposition \ref{prop_standard_form} tells us that at each point, we can assume without loss of generality that the defining polynomial of a given (G)PSR manifold to be of a relatively simple form, which we will now give a name for further referencing.	
	
	\begin{Def}\label{def_standard_form}
		A (G)PSR manifold $\mathcal{H}\subset\{h=1\}$, $h$ of degree $\tau\geq 3$, is in \textbf{standard form} if $\left(\begin{smallmatrix}x\\y\end{smallmatrix}\right)=\left(\begin{smallmatrix}1\\0\end{smallmatrix}\right)\in\mathcal{H}$ and $h$ is of the form
			\begin{equation}\label{eqn_h_standard_form}
				h=x^\tau-x^{\tau-2}\langle y,y\rangle + \sum\limits_{k=3}^\tau x^{\tau-k}P_k(y).
			\end{equation}
		Analogously, if the defining polynomial $h$ is of the form \eqref{eqn_h_standard_form} we will say that $h$ is in \textbf{standard form}.
	\end{Def}
	
		Note that with the above definition a CC(G)PSR manifold $\mathcal{H}\subset\{h=1\}$ being in standard form means that the defining polynomial $h$ is in standard form and that $\mathcal{H}$ is precisely the connected component of $\{h=1\}\cap\{\text{hyperbolic points of }h\}$ that contains the point $\left(\begin{smallmatrix}1\\0\end{smallmatrix}\right)$. Furthermore, observe that the compact convex generating set $\mathcal{C}_n$ of the moduli space of CCPSR manifolds in \eqref{eqn_CCPSR_gen_set} is defined by a condition on the $P_3$-part of the standard forms of the corresponding defining polynomials.
		
		In the proof of Theorems \ref{thm_quartic_CCGPSR_curves_classification} and \ref{thm_incomplete_quartics} we will make heavy use of an infinitesimal version of Proposition \ref{prop_standard_form} which is to be understood as follows. Suppose that we are given a (G)PSR manifold $\mathcal{H}\subset\{h=1\}$ in standard form with corresponding terms $P_3,\ldots,P_\tau$ as in \eqref{eqn_h_standard_form}. By definition, the point $\left(\begin{smallmatrix}1\\0\end{smallmatrix}\right)$ is contained in $\mathcal{H}$, so we can choose a sufficiently small open neighbourhood $U$ of that point in $\mathcal{H}$ and study how the terms $P_3,\ldots,P_\tau$ change when infinitesimally changing the reference point $p\in U$ for the standard form of the defining polynomial. In formulas, that means studying
			\begin{equation}\label{eqn_A_pullback_h_map_differential}
				\left.\D \left(A(\cdot)^*h\right)\right|_{p=\left(\begin{smallmatrix}1\\0\end{smallmatrix}\right)}
			\end{equation}
		where $A:U\to\mathrm{GL}(n+1)$ is chosen to fulfil Proposition \ref{prop_standard_form} \eqref{prop_standard_form_eqn_i} and \eqref{prop_standard_form_eqn_ii} and $A(\cdot)^*h$ is viewed as a smooth map from $U$ to the vector space of homogeneous polynomials of degree $\tau$ which we identify with $\mathrm{Sym}^\tau\left(\mathbb{R}^n\right)^*$, i.e.
			\begin{equation*}\label{eqn_A_pullback_h_map}
				A(\cdot)^*h:U\to \mathrm{Sym}^\tau\left(\mathbb{R}^n\right)^*.
			\end{equation*}
		By construction, $A(\cdot)^*h$ is of the form
			\begin{equation}
				A(p)^*h=x^\tau-x^{\tau-2}\langle y,y\rangle + \sum\limits_{k=3}^\tau x^{\tau-k} P_k(p)(y),
			\end{equation}
		where for all $3\leq k\leq \tau$, $P_k(\cdot)(y):U\to \mathrm{Sym}^k\left(\mathbb{R}^n\right)^*$ is a smooth map. Hence, we obtain that the differential of $A(\cdot)^*h$ at $\left(\begin{smallmatrix}1\\0\end{smallmatrix}\right)$ \eqref{eqn_A_pullback_h_map_differential} is of the form
			\begin{equation}\label{eqn_diff_standard_form_explicit}
				\left.\D \left(A(\cdot)^*h\right)\right|_{p=\left(\begin{smallmatrix}1\\0\end{smallmatrix}\right)}= \sum\limits_{k=3}^\tau x^{\tau-k} \delta P_k(y),\quad \delta P_k(y):= \left.\D P_k(\cdot)(y)\right|_{p=\left(\begin{smallmatrix}1\\0\end{smallmatrix}\right)}\quad \forall\,3\leq k\leq \tau.
			\end{equation}
		In the above equation \eqref{eqn_diff_standard_form_explicit}, the terms $\delta P_k(y)$ for $3\leq k\leq \tau$ are viewed as linear maps
			\begin{equation}
				\delta P_k(y):\mathbb{R}^n\to\mathrm{Sym}^k\left(\mathbb{R}^n\right)^*
			\end{equation}
		where we have identified $T_{\left(\begin{smallmatrix}1\\0\end{smallmatrix}\right)}\mathcal{H}\cong \mathbb{R}^n$ via the choice of linear coordinates $y$ in \eqref{eqn_h_standard_form} on the $\{0\}\times\mathbb{R}^n$-part of the ambient space. Now we have to investigate the possible choices of $A:U\to\mathrm{GL}(n+1)$ as in Proposition \ref{prop_standard_form} which can lead to different explicit formulas for the terms $\delta P_k(y)$. Note that $A$ not being unique follows from the fact that an orthogonal transformation in the $y$-coordinates might change the terms $P_3,\ldots,P_\tau$ but will leave the polynomial $h$ in standard form. To do so we will use a result of \cite{L2} for the explicit form of the map $A$ and its differential at $\left(\begin{smallmatrix}1\\0\end{smallmatrix}\right)$. In the following, let $\partial_x h:=\frac{\partial h}{\partial x}$ and $\partial_y h:=\left(\frac{\partial h}{\partial y_1},\ldots,\frac{\partial h}{\partial y_n}\right)$.
		
		\begin{Prop}\label{prop_A_explicit_form}
			Let $\mathcal{H}$ be a (G)PSR manifold in standard form, $U$ an open neighbourhood of $\left(\begin{smallmatrix}1\\0\end{smallmatrix}\right)$ in $\mathcal{H}$, and $A:U\to\mathrm{GL}(n+1)$ a smooth map fulfilling Proposition \ref{prop_standard_form} \eqref{prop_standard_form_eqn_i} and \eqref{prop_standard_form_eqn_ii}. Additionally, let $A\left(\left(\begin{smallmatrix}1\\0\end{smallmatrix}\right)\right)=\mathbbm{1}$. Then $A$ is of the form
				\begin{equation}\label{eqn_A_explicit}
					A(p)=\left(
						\begin{array}{c|c}
							x(p) & \left.-\frac{\partial_y h}{\partial_x h}\right|_p\cdot \underset{}{E(p)}\\ \hline
							y(p) & \vphantom{\Big(}E(p)
						\end{array}
					\right)
				\end{equation}
			where $E:U\to\mathrm{GL}(n)$ is a smooth map fulfilling
				\begin{equation*}\label{eqn_E_implicit}
					-\frac{1}{2}\partial^2 h_p\left(\left(
						\begin{array}{c}
							\left.-\frac{\partial_y h}{\partial_x h}\right|_p\cdot E(p)y\\
							E(p)y
						\end{array}
					\right),\left(
						\begin{array}{c}
							\left.-\frac{\partial_y h}{\partial_x h}\right|_p\cdot E(p)y\\
							E(p)y
						\end{array}
					\right)\right)=\langle y,y\rangle
				\end{equation*}
			for all $p\in U$ and all $y\in\mathbb{R}^n$, and $E\left(\left(\begin{smallmatrix}1\\0\end{smallmatrix}\right)\right)=\mathbbm{1}$. The differential of $A$ at $\left(\begin{smallmatrix}1\\0\end{smallmatrix}\right)$ fulfils
				\begin{equation}\label{eqn_differential_A_formula}
					\left.\D A(\cdot)\right|_{p=\left(\begin{smallmatrix}1\\0\end{smallmatrix}\right)} = \left(\begin{array}{c|c}
						\underset{}{0} & \frac{2}{\tau}\D y^T\\ \hline
						\D y & \overset{}{\frac{3}{2}}P_3(y,\cdot,\D y)^T + by
					\end{array}\right),
				\end{equation}
			where we have identified $P_3$ with its corresponding trilinear form, $b:\mathbb{R}^n\to \mathfrak{so}(n)$ is a linear map of the form
				\begin{equation}\label{eqn_b_explicit}
					b=\sum\limits_{i=1}^n a_i\otimes \D y_i,
				\end{equation}
			and $by$ is to be understood as $by=\sum\limits_{i=1}^n (a_i\cdot y)\otimes \D y_i$.
			\begin{proof}
				\cite[Prop.\,3.1,\;Prop.\,3.4]{L2} with a slightly different notation.
			\end{proof}
		\end{Prop}
		
		Note that in the definition of the map $b$ we have again identified $\mathbb{R}^n$ with the tangent space of $\mathcal{H}$ at $\left(\begin{smallmatrix}1\\0\end{smallmatrix}\right)$. For explicit formulas in the most general case $\tau\geq 3$ for all $\delta P_k(y)$, $3\leq k\leq \tau$, see \cite[eqn.\,(3.21)]{L2}. For the purpose of this work we will only need a formula for $\delta P_3(y)$ and $\delta P_4(y)$ for the quartic case $\tau=4$. Inserting \eqref{eqn_differential_A_formula} into \eqref{eqn_diff_standard_form_explicit} for $\tau=4$ yields
			\begin{align}
				\delta P_3(y)&= -\langle y,y\rangle\langle y,\D y\rangle + \D P_3|_y\left(by + \frac{1}{4}\partial^2 P_3|_y\cdot \D y^T\right) + \D P_4|_y(\D y),\label{eqn_deltaP3_quartics}\\
				\delta P_4(y)&= \frac{1}{2}P_3(y)\langle y,\D y\rangle + \D P_4|_y\left(by + \frac{1}{4}\partial^2 P_3|_y\cdot \D y^T\right).\label{eqn_deltaP4_quartics}
			\end{align}
		The terms $\delta P_k(y)$ in \eqref{eqn_diff_standard_form_explicit} together with \eqref{eqn_differential_A_formula} give us a convenient way to formulate a condition for a (G)PSR manifold to be a homogeneous space.
	
	\begin{Prop}\label{prop_hom_iff_deltaPs_vanish}
		A maximal connected (G)PSR manifold in standard form of dimension $n$ with defining polynomial of degree $\tau$ is a homogeneous space if and only if there exists a linear map $b:\mathbb{R}^n\to\mathfrak{so}(n)$ as in \eqref{eqn_b_explicit}, such that $\delta P_k(y)=0$ for all $3\leq k\leq \tau$.
		\begin{proof}
			\cite[Prop.\,3.12]{L2}.
		\end{proof}
	\end{Prop}
	
	Since homogeneous (G)PSR manifolds are homogeneous Riemannian spaces, they are automatically geodesically complete. In general it is currently not known whether the condition for a GPSR manifold to be closed in its ambient space is sufficient for geodesic completeness with respect to the centro-affine fundamental form. For closed PSR manifolds this property has first been proven in \cite{CNS} and two alternative proofs can be found in \cite[Prop.\,4.17, Prop.\,5.17]{L1}.
	
	A general classification of homogeneous (G)PSR manifolds in all dimensions for every degree of the corresponding defining polynomials is, at least according to our opinion, currently far out of reach. However, homogeneous PSR manifolds have been completely classified in \cite{DV}. In Theorem \ref{thm_quartic_CCGPSR_curves_classification} we will classify all complete quartic GPSR curves, so we will in particular obtain a list of all homogeneous ones. In higher dimensions, homogeneous quartic GPSR manifolds have not yet been classified but we will show in Proposition \ref{prop_hom_quartics_sing_at_inf} that they all necessarily posses the following property.
	
	\begin{Def}\label{def_sing_at_inf}
		A (G)PSR manifold $\mathcal{H}\subset\{h=1\}$ is called \textbf{singular at infinity} if there exists a point $p\ne 0$ in the boundary of its cone $\mathbb{R}_{>0}\cdot \mathcal{H}$, that is $p\in\partial\left(\mathbb{R}_{>0}\cdot \mathcal{H}\right)\setminus\{0\}$, such that $h(p)=0$ and $\D h_p=0$. 
	\end{Def}
	
	If a (G)PSR manifold $\mathcal{H}\subset\{h=1\}$ is singular at infinity or admits a point $p\in\partial(\mathbb{R}_{>0}\cdot\mathcal{H})\cap\{h=0\}\setminus\{0\}$, such that $-\partial^2h_p$ has at least $2$-dimensional kernel, $\mathcal{H}$ is said to have \textbf{non-regular boundary behaviour}.
	
	Note that for closed (G)PSR manifolds the first condition $h(p)=0$ in Definition \ref{def_sing_at_inf} is always satisfied. In order to check whether or not a given maximal (G)PSR manifold is singular at infinity it actually suffices to find a maximal connected (G)PSR curve contained in said manifold that is singular at infinity.
	
	\begin{Lem}\label{lem_sing_at_inf_iff_restr_to_curve}
		A maximal (G)PSR manifold $\mathcal{H}\subset\mathbb{R}^{n+1}$ is singular at infinity if and only if there exists a maximal connected (G)PSR curved $\widetilde{\mathcal{H}}$ obtained by intersecting a connected component of $\mathcal{H}$ with a plane in $\mathbb{R}^{n+1}$, such that $\widetilde{\mathcal{H}}$ is singular at infinity.
		\begin{proof}
			Follows from \cite[Lem.\,3.5\,\textit{(iii)}]{L2} if $\mathcal{H}$ is closed. If $\mathcal{H}$ is maximal and not closed the same line of arguments can be used. One needs only replace the convexity of $\mathrm{dom}(\mathcal{H})$ as defined in the upcoming Definition \ref{def_dom_H} with star-shaped with center $0\in\mathbb{R}^n$ in the proof of \cite[Lem.\,3.5]{L2} to obtain the desired conclusion.
		\end{proof}
	\end{Lem}
	
	\begin{rem}\label{rem_motivation_quartic_homs_sing_at_inf}
		It has been shown in \cite[Prop.\,4.6]{L2} that homogeneous PSR manifolds are singular at infinity. More precisely, it was shown that every critical value of the term $P_3$ restricted to the Euclidean unit sphere $S^{n-1}$ in the corresponding defining polynomials in standard form turned out to be of absolute value $\frac{2}{3\sqrt{3}}$. In the quartic case we will prove a similar statement and find in Corollary \ref{cor_crit_values_quartic_homs} that there is not one, but two such possible values for the pair $(|P_3|,P_4)$.
	\end{rem}
	
	CC(G)PSR manifolds $\mathcal{H}\subset\{h=1\}\subset\mathbb{R}^{n+1}$ can be parametrised over any precompact section of their spanned coned $\mathbb{R}_{>0}\cdot\mathcal{H}$ with an affinely embedded $\mathbb{R}^n\hookrightarrow\mathbb{R}^{n+1}$ via $p\mapsto h^{-1/\tau}(p)p$, where $\tau=\deg(h)$. If $\mathcal{H}$ is in standard form, a natural choice for the affinely embedded $\mathbb{R}^n$ is $T_{\left(\begin{smallmatrix}1\\0\end{smallmatrix}\right)}\mathcal{H}$. We will use the so-obtain intersection with the cone spanned by $\mathcal{H}$ multiple times and thus give it a name. The following is also well defined for maximal incomplete connected (G)PSR manifolds.
	
	\begin{Def}\label{def_dom_H}
		Let $\mathcal{H}$ be a maximal connected (G)PSR manifold in standard form. We call the set
			\begin{equation*}
				\mathrm{dom}(\mathcal{H}):=\mathrm{pr}_{\mathbb{R}^n}\left(\left(\mathbb{R}_{>0}\cdot\mathcal{H}\right)\cap \left(\left(\begin{smallmatrix}1\\0\end{smallmatrix}\right)+T_{\left(\begin{smallmatrix}1\\0\end{smallmatrix}\right)}\mathcal{H}\right)\right)
			\end{equation*}
		the \textbf{domain} of $\mathcal{H}$, where $\mathrm{pr}_{\mathbb{R}^n}:\mathbb{R}^{n+1}\to\mathbb{R}^n$ denotes the projection onto the last $n$ coordinates.
	\end{Def}
	
	\begin{figure}[H]%
		\centering%
		\includegraphics[scale=0.2]{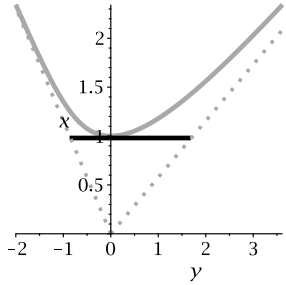}%
		\caption{A typical plot of a CC(G)PSR curve $\mathcal{H}$ in grey, with the boundary of its spanned cone marked with two dotted grey rays, and $\mathrm{dom}(\mathcal{H})$ in black.}
	\end{figure}
	
	Note that in general this does not define $\mathrm{dom}(\mathcal{H})$ uniquely as, in general, the standard form of $\mathcal{H}$ is not uniquely determined but might vary depending on both the reference point $p\in\mathcal{H}$ and a choice of coordinates. When parametrised over $\mathrm{dom}(\mathcal{H})$, the centro-affine metric of $\mathcal{H}$ can be written in a simple form.
	
	\begin{Lem}\label{lem_pullback_g}
		Let $\Phi:\mathrm{dom}(\mathcal{H})\to\mathcal{H}$ denote the central parametrisation of a given maximal connected (G)PSR manifold $\mathcal{H}\subset\{h=1\}$ in standard form with $\mathrm{deg}(h)=\tau\geq 3$. Let $f:= h\left(\left(\begin{smallmatrix}1\\y\end{smallmatrix}\right)\right)$. Then the pullback of the centro-affine fundamental form $g_\mathcal{H}$ of $\mathcal{H}$ is given by
			\begin{equation}
				\Phi^*g_\mathcal{H}=-\frac{\partial^2f}{\tau f} + \frac{(\tau-1)\D f^2}{\tau^2 f^2}.\label{eqn_pullback_metric_to_dom}
			\end{equation}
		\begin{proof}
			\cite[Cor.\,1.13]{CNS}.
		\end{proof}
	\end{Lem}
	
	The following holds for $\mathrm{dom}(\mathcal{H})$ independent of whether $\mathcal{H}$ is closed or not.
	
	\begin{Lem}
		$\mathrm{dom}(\mathcal{H})$ is precompact for all maximal connected (G)PSR manifolds $\mathcal{H}$.
		\begin{proof}
			For any given $\mathcal{H}\subset\mathbb{R}^{n+1}$ in standard form, the intersection of $\mathcal{H}$ with a plate $P\subset\mathbb{R}^{n+1}$ is a maximal connected (G)PSR curve. It thus suffices to show that for any hyperbolic homogeneous polynomial $h$ in two real variables that is in standard form, the connected component of $U:=\left\{y\in\mathbb{R}\ \left|\ \left(\begin{smallmatrix}1\\y\end{smallmatrix}\right)\text{ is a hyperbolic point of }h\right\}\right.$ containing $y=0$ is precompact in $\mathbb{R}$. If $h=0$ has a positive solution, which in particular includes all cases of closed connected (G)PSR curves, $U$ is precompact. Suppose that $h=0$ has no positive solution. Since $y=0$ is a local maximum of $h\left(\left(\begin{smallmatrix}1\\y\end{smallmatrix}\right)\right)$ and the latter is a polynomial, $h=0$ having no positive solution implies that $h\left(\left(\begin{smallmatrix}1\\y\end{smallmatrix}\right)\right)\to\infty$ as $y\to\infty$ and, hence, that $h\left(\left(\begin{smallmatrix}1\\y\end{smallmatrix}\right)\right)$ has a positive local minimum with $y>0$. At that point, $dh_y=0$ and \eqref{eqn_pullback_metric_to_dom} implies that $\left(\begin{smallmatrix}1\\y\end{smallmatrix}\right)$ is not a hyperbolic point of $h$. We conclude that $U$ is precompact in these cases as well.
		\end{proof}
	\end{Lem}
	
	We will need one more technical result that will play a role in our proof of Theorem \ref{thm_quartic_CCGPSR_curves_classification}.
	
	\begin{Lem}\label{lem_partial_x_h_not_vanishing}
		Let $\mathcal{H}\subset\{h=1\}$ be a maximal (G)PSR curve in standard form. Then $\partial_x h$ does not vanish on the cone $\mathbb{R}_{>0}\cdot\mathcal{H}$.
		\begin{proof}
			Let $\tau=\mathrm{deg}(h)$. Suppose there exists a counterexample $\mathcal{H}\subset\{h=1\}$ to our claim. By the homogeneity of $h$ we might assume that $\partial_x h|_p=0$ for some $p\in\mathcal{H}$. Assume without loss of generality that $y(p)>0$, which can always be achieved via a possible sign-flip in the $y$-coordinate. Then $\partial_y h|_p$ must be negative. This follows from the fact that $h\left((1-t)\left(\begin{smallmatrix}1\\0\end{smallmatrix}\right)+tp\right)>1$ for all $t\in(0,1)$ by the hyperbolicity of $h$. We thus obtain using $y(p)>0$
				\begin{equation*}
					dh_p(p)=\partial_xh|_p x(p)+\partial_yh|_p y(p)<0.
				\end{equation*}
			But Euler's theorem for homogeneous functions implies $dh_p(p)=\tau>0$, a contradiction.
		\end{proof}
	\end{Lem}

\section{Classification results}
	In this section we classify maximal quartic GPSR curves up to equivalence and study the linear automorphism groups of their defining polynomials. We start with the case of closed maximal quartic CCGPSR curves in Theorem \ref{thm_quartic_CCGPSR_curves_classification} and, after that, will deal whit the maximal non-closed cases in Theorem \ref{thm_incomplete_quartics}. In the following we assume that integral curves of vector fields are connected and parametrised over an open interval.
	
	\begin{Th}\label{thm_quartic_CCGPSR_curves_classification}
		Every quartic CCGPSR curve is equivalent to precisely one of the following quartic CCGPSR curves $\mathcal{H}\subset\{h=1\}$ in standard form:
			\begin{enumerate}[a)]
				\item $h=x^4-x^2y^2+\frac{1}{4}y^4$, $\{h=1\}\cap\{\text{hyperbolic points of }h\}$ has $4$ equivalent closed connected components,\label{eqn_qCCPSRcurves_class_a} and
					\begin{equation*}
						G^h\cong\mathbb{R}\ltimes(\mathbb{Z}_2\ltimes\mathbb{Z}_2),
					\end{equation*}
				where the $\mathbb{R}$-factor acts by hyperbolic rotations with respect to the Minkowski inner product $\D x^2-\tfrac{1}{2}\D y^2$, the first $\mathbb{Z}_2$-factor acts via $y\to-y$, and the second $\mathbb{Z}_2$-factor acts via $(x,y)\to\left(-\tfrac{1}{\sqrt{2}}y,\sqrt{2}x\right)$,
				\item $h=x^4-x^2y^2+\frac{2\sqrt{2}}{3\sqrt{3}}xy^3-\frac{1}{12} y^4=1$, $\{h=1\}\cap\{\text{hyperbolic points of }h\}$ has $2$ equivalent closed connected components,\label{eqn_qCCPSRcurves_class_b} and
					\begin{equation*}
						G^h\cong\mathbb{R}\times\mathbb{Z}_2,
					\end{equation*}
				where the $\mathbb{R}$-factor acts on $\mathcal{H}$ as described in equation \eqref{eqn_nontriv_R_action_b},
				\item $h=x^4-x^2y^2+\frac{2}{3\sqrt{3}}xy^3=1$, $\{h=1\}\cap\{\text{hyperbolic points of }h\}$ has $4$ equivalent closed connected components,\label{eqn_qCCPSRcurves_class_c} and
					\begin{equation*}
						G^h\cong\mathbb{Z}_2\times\mathbb{Z}_2,
					\end{equation*}
				\item $h=x^4-x^2y^2+Ky^4=1$ for exactly one $K<\frac{1}{4}$. The set $\{h=1\}\cap\{\text{hyperbolic points of }h\}$ has $4$ equivalent closed connected components for $0<K<\frac{1}{4}$\label{eqn_qCCPSRcurves_class_d} with
					\begin{equation*}
						G^h\cong\mathbb{Z}_2\ltimes\mathbb{Z}_2
					\end{equation*}
				where the first $\mathbb{Z}_2$-factor acts via $y\to-y$ and the second via $(x,y)\to\left(-\sqrt[4]{K}y,\sqrt[4]{K}^{-1}x\right)$, and $2$ equivalent closed connected components for $K\leq 0$ with
					\begin{equation*}
						G^h\cong\mathbb{Z}_2\times\mathbb{Z}_2,
					\end{equation*}
				where the first $\mathbb{Z}_2$-factor acts via $y\to-y$ and the second via $x\to -x$.
			\end{enumerate}
		The curves in \hyperref[eqn_qCCPSRcurves_class_a]{a)} and \hyperref[eqn_qCCPSRcurves_class_b]{b)} are homogeneous spaces, the curves \hyperref[eqn_qCCPSRcurves_class_c]{c)} and \hyperref[eqn_qCCPSRcurves_class_d]{d)} are inhomogeneous. All of the curves in \hyperref[eqn_qCCPSRcurves_class_a]{a)}, \hyperref[eqn_qCCPSRcurves_class_b]{b)}, and \hyperref[eqn_qCCPSRcurves_class_c]{c)} are singular at infinity, and each element in the one-parameter family \hyperref[eqn_qCCPSRcurves_class_d]{d)} is not singular at infinity.%\todo{check that this is shown!}
	\begin{proof}
		Let $\mathcal{H}$ be a maximal connected quartic GPSR curve. Proposition \ref{prop_standard_form} implies that we can assume without loss of generality that $\mathcal{H}=\mathcal{H}_{L,K}$ is the connected component of $\{h_{L,K}=1\}\cap\{\text{hyperbolic points of }h_{L,K}\}$, where
			\begin{equation}\label{eqn_quartic_curves_std_form}
				h_{L,K}=x^4-x^2y^2+Lxy^3+Ky^4
			\end{equation}
		for some $(L,K)^T\in\mathbb{R}^2$, that contains the point $\left(\begin{smallmatrix}
		x\\ y
		\end{smallmatrix}
		\right)=\left(\begin{smallmatrix}
		1\\ 0
		\end{smallmatrix}
		\right)$. This means that $P_3(y)=Ly^3$ and $P_4(y)=Ky^4$ in equation \eqref{eqn_h_standard_form}, leading to \eqref{eqn_quartic_curves_std_form}. We will say that the polynomial $h_{L,K}$ ``corresponds'' to the point $(L,K)^T\in\mathbb{R}^2$. On the other hand, note that for all $(L,K)^T\in\mathbb{R}^2$, the connected component of $\{h_{L,K}=1\}\cap\{\text{hyperbolic points of }h_{L,K}\}$ that contains $\left(\begin{smallmatrix}
		x\\ y
		\end{smallmatrix}
		\right)=\left(\begin{smallmatrix}
		1\\0
		\end{smallmatrix}\right)$ is by definition a maximal connected quartic GPSR curve. This proof primarily relies on the properties of $\delta P_3(y)$ and $\delta P_4(y)$, defined in equation \eqref{eqn_diff_standard_form_explicit} and specified for quartics in \eqref{eqn_deltaP3_quartics} and \eqref{eqn_deltaP4_quartics}. Since $\dim\left(\mathcal{H}_{L,K}\right)=1$, the term $b$ in equations \eqref{eqn_deltaP3_quartics} and \eqref{eqn_deltaP4_quartics} vanishes and we obtain
			\begin{align}
				\delta P_3(y)&=\left(\frac{9}{2}L^2+4K-1\right)y^3 \D y,\label{eqn_deltaP3_quarticcurves}\\
				\delta P_4(y)&=L\left(6K+\frac{1}{2}\right)y^4 \D y.\label{eqn_deltaP4_quarticcurves}
			\end{align}
		In the above formulas, $y$ denotes the induced coordinate of $\mathrm{dom}(\mathcal{H}_{L,K})$, cf. Definition \ref{def_dom_H}. This motivates the consideration of the vector field
			\begin{equation}\label{eqn_VFqCCGPSRc_def}
				\mathcal{V}=\mathcal{V}_{\left(\begin{smallmatrix}
				L\\ K
				\end{smallmatrix}\right)}:=\left(\frac{9}{2}L^2+4K-1\right)\partial_L+L\left(6K+\frac{1}{2}\right)\partial_K\in\mathfrak{X}(\mathbb{R}^2).
			\end{equation}
		See Figure \ref{fig_V_fieldplot} for a plot of $\mathcal{V}$.
			\begin{figure}[H]%
				\centering%
				\includegraphics[scale=0.3]{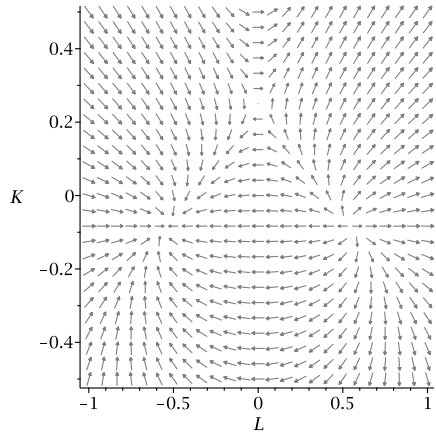}%
				\caption{A plot of $\mathcal{V}$ created in Maple using the option \texttt{fieldstrength=log[8]}.}\label{fig_V_fieldplot}%
			\end{figure}
		\noindent
		Observe that
			\begin{equation}
				\{\mathcal{V}=0\}=\left.\left\{\left(\begin{matrix}
				L\\ K
				\end{matrix}
				\right)\in\mathbb{R}^2\ \right|\ \mathcal{V}_{\left(\begin{smallmatrix}
				L\\ K
				\end{smallmatrix}
				\right)}=0\right\}=\left\{\left(\begin{matrix}
				0\\ \frac{1}{4}
				\end{matrix}
				\right),\left(\begin{matrix}
				\frac{2\sqrt{2}}{3\sqrt{3}}\\ -\frac{1}{12}
				\end{matrix}
				\right),\left(\begin{matrix}
				-\frac{2\sqrt{2}}{3\sqrt{3}}\\ -\frac{1}{12}
				\end{matrix}
				\right)\right\}.\label{eqn_mathcalV_zeros}
			\end{equation}
		$\mathcal{V}$ has the property that the polynomials
			\begin{equation}\label{eqn_quadpoly_integralcurves}
				h_{\gamma_L(t),\gamma_K(t)}=x^4-x^2y^2+\gamma_L(t)xy^3+\gamma_K(t)y^4
			\end{equation}
		associated to each integral curve and in particular each maximal integral curve $\gamma$ of the restricted vector field $\mathcal{V}|_{\mathbb{R}^2\setminus\{\mathcal{V}=0\}}$,
			\begin{equation*}
				t\mapsto \gamma(t)=\left(\begin{smallmatrix}
				\gamma_L(t)\\ \gamma_K(t)
				\end{smallmatrix}
				\right)
				\in\mathbb{R}^2,\quad \mathcal{V}_{\gamma}=\dot\gamma,\quad \gamma_L(0)=L,\ \gamma_K(0)=K,
			\end{equation*}
		for all initial values $(L,K)^T\in\mathbb{R}^2\setminus\{\mathcal{V}=0\}$, are equivalent to $h_{L,K}$. To see that this is true, we will use the techniques of Propositions \ref{prop_standard_form} and \ref{prop_A_explicit_form}. Firstly, we parametrise $\mathcal{H}_{L,K}$ over the intersection of its spanned cone with its tangent space at $\left(\begin{smallmatrix}1\\0\end{smallmatrix}\right)$ given by the diffeomorphism $\Phi:\mathrm{dom}(\mathcal{H}_{L,K})\to\mathcal{H}_{L,K}$, $T\mapsto \sqrt[4]{h_{L,K}\left(\left(\begin{smallmatrix}1\\ T\end{smallmatrix}\right)\right)}^{-1}\left(\begin{smallmatrix}1\\ T\end{smallmatrix}\right)$. We further define $\mathcal{A}:=A\circ \Phi$ with $A:\mathcal{H}_{L,K}\to\mathrm{GL}(2)$ as in \eqref{eqn_A_explicit}. We calculate $\mathcal{A}:\mathrm{dom}\left(\mathcal{H}_{L,K}\right)\to\mathrm{GL}(2)$ explicitly and find
			\begin{equation}
				\mathcal{A}(T)=\left(\begin{array}{c|c}
				\underset{}{\sqrt[4]{h_{L,K}\left(\left(\begin{smallmatrix}1\\ T\end{smallmatrix}\right)\right)}^{-1}} & \frac{2T-3LT^2-4KT^3}{4-2T^2+LT^3}r(L,K,T)\\ \hline
				\overset{}{\sqrt[4]{h_{L,K}\left(\left(\begin{smallmatrix}1\\ T\end{smallmatrix}\right)\right)}^{-1} T} & \overset{}{r(L,K,T)}
				\end{array}
				\right)\label{eqn_A_pulled_back_to_dom}
			\end{equation}
		with
			\begin{align*}
				r(L,K,T)&=\sqrt[4]{h_{L,K}\left(\left(\begin{smallmatrix}1\\ T\end{smallmatrix}\right)\right)}\left(4-2T^2+LT^3\right)\\
				&\ \cdot\left(16-48LT+(-8-96K)T^2+56LT^3+\left(-8-42L^2+128K\right)T^4\right.\\
				&\left.\quad+(16-144LK)T^5+\left(-14L^2-8K-96K^2\right)T^6\right.\\
				&\left.\quad+\left(6L^3+8LK\right)T^7+\left(6L^2K+16K^2\right)T^8\right)^{-\frac{1}{2}}.
			\end{align*}
		It follows from Lemma \ref{lem_partial_x_h_not_vanishing} that $\partial_x h_{L,K}|_{\Phi(T)}> 0$ for all $T\in \mathrm{dom}\left(\mathcal{H}_{L,K}\right)$ and, hence, that $\mathcal{A}(T)$ is well-defined for all $T\in\mathrm{dom}\left(\mathcal{H}_{L,K}\right)$ and not just in some neighbourhood of $0\in\mathrm{dom}\left(\mathcal{H}_{L,K}\right)$ as implied by Proposition \ref{prop_standard_form}. With
			\begin{equation*}
				h_{L(T),K(T)}=h_{L,K}\circ \mathcal{A}(T)=x^4-x^2y^2+L(T)xy^3+K(T)y^4
			\end{equation*}
		we then obtain
			\begin{align}
				L(T)&=-2\sqrt{2}\sqrt{h_{L,K}\left(\left(\begin{smallmatrix}1\\ T\end{smallmatrix}\right)\right)}\notag\\
				&\ \cdot\left(-8L+(8-32K)T-20LT^2+20L^2T^3+40LKT^4\right.\notag\\
				&\left.\quad+\left(-2L^2-8K+32K^2\right)T^5+\left(L^3+4LK\right)T^6\right)\notag\\
				&\ \cdot \left(8-24LT+(4-48K)T^2-4LT^3+\left(3L^2+8K\right)T^4\right)^{-1}\notag\\
				&\ \cdot \left(8-24LT+(-4-48K)T^2+28LT^3+\left(-4-21L^2+64K\right)T^4\right.\notag\\
				&\left.\quad+\left(8L-72LK\right)T^5+\left(-7L^2-4K-48K^2\right)T^6\right.\notag\\
				&\left.\quad+\left(3L^3+4LK\right)T^7+\left(3L^2K+8K^2\right)T^8\right)^{-\frac{1}{2}}\label{eqn_L_T_pmoving}
			\end{align}
		and
			\begin{align}
				K(T)&=\frac{1}{4}\left(256K+128LT+\left(-64-192L^2-256K\right)T^2+\left(128L-256LK\right)T^3 \right.\notag\\
				&\left.\quad\quad+\left(16-80L^2+384K-256K^2\right)T^4+\left(48L^3-32L-512LK\right)T^5 \right.\notag\\
				&\left.\quad\quad+\left(8L^2+352L^2K-64K-256K^2\right)T^6+\left(8L^3+64LK+512LK^2 \right)T^7\right.\notag\\
				&\left.\quad\quad+\left(-3L^4-16L^2K+256K^3\right)T^8\right)\notag\\
				&\ \cdot \left(8-24LT+(4-48K)T^2-4LT^3+\left(3L^2+8K\right)T^4 \right)^{-2},\label{eqn_K_T_pmoving}
			\end{align}
		where $L(0)=L$, $K(0)=K$. Note that
			\begin{equation}\label{eqn_delLTKT_at_0}
				\left.\frac{\partial}{\partial T}\left(\begin{matrix}
				L(T)\\ K(T)
				\end{matrix}
				\right)\right|_{T=0}=\mathcal{V}_{\left(\begin{smallmatrix}
				L\\ K
				\end{smallmatrix}
				\right)}
			\end{equation}
		as expected. By construction we know that for all $T\in\mathrm{dom}\left(\mathcal{H}_{L,K}\right)$, $h_{L(T),K(T)}$ and $h_{L,K}$ are equivalent. Since $\mathcal{A}$ depends smoothly on $T$, the connected components of $\{h_{L,K}=1\}\cap\{\text{hyperbolic points of }h_{L,K}\}$ that contain the point $\left(\begin{smallmatrix}
		x\\ y
		\end{smallmatrix}
		\right)=\left(\begin{smallmatrix}
		1\\0
		\end{smallmatrix}\right)$ are pairwise equivalent for all $T\in\mathrm{dom}\left(\mathcal{H}_{L,K}\right)$. The Euclidean velocities of the considered integral curve $\gamma$ of $\mathcal{V}|_{\mathbb{R}^2\setminus\{\mathcal{V}=0\}}$ and the curve $T\mapsto\left(\begin{smallmatrix}
		L(T)\\ K(T)
		\end{smallmatrix}
		\right)$ will in general not coincide, exceptions being the constant integral curves of $\mathcal{V}$, although these might not be the only exceptions. It however still holds true that the image of $\gamma$ is always contained in the image of $T\mapsto\left(\begin{smallmatrix}
		L(T)\\ K(T)
		\end{smallmatrix}
		\right)$, which we will show now. The constant integral curves of $\mathcal{V}$ are precisely those with initial values in $\{\mathcal{V}=0\}$ described in equation \eqref{eqn_mathcalV_zeros}. One now checks that for $\left(\begin{smallmatrix}
		L\\ K
		\end{smallmatrix}
		\right)\in\{\mathcal{V}=0\}$, $L(T)\equiv L(0)=L$ and $K(T)\equiv K(0)=K$. Hence, for constant integral curves the maps $T\mapsto L(T)$ in \eqref{eqn_L_T_pmoving} and $T\mapsto K(T)$ in \eqref{eqn_K_T_pmoving} are constant and thus the images of $\gamma$ and $(L(T),K(T))^T$ in $\mathbb{R}^2$ coincide. For all $\left(\begin{smallmatrix}
		L\\ K
		\end{smallmatrix}
		\right)\in\mathbb{R}^2\setminus\{\mathcal{V}=0\}$
		one can now verify that
			\begin{equation*}
				\D L\left(\mathcal{V}_{\left(\begin{smallmatrix}
				L(T)\\ K(T)
				\end{smallmatrix}
				\right)}\right)\cdot\partial_T K(T)=\D K\left(\mathcal{V}_{\left(\begin{smallmatrix}
				L(T)\\ K(T)
				\end{smallmatrix}
				\right)}\right)\cdot\partial_T L(T)
			\end{equation*}
		for all $T\in\mathrm{dom}\left(\mathcal{H}_{L,K}\right)$, and we find for all $\left(\begin{smallmatrix}
		L\\ K
		\end{smallmatrix}
		\right)\in\mathbb{R}^2\setminus\{\mathcal{V}=0\}$
		\begin{align}
			&\frac{dL\left(\mathcal{V}_{\left(\begin{smallmatrix}
			L(T)\\ K(T)
			\end{smallmatrix}
			\right)}\right)}{\partial_T L(T)} =\frac{dK\left(\mathcal{V}_{\left(\begin{smallmatrix}
			L(T)\\ K(T)
			\end{smallmatrix}
			\right)}\right)}{\partial_T K(T)}\notag\\ &=\frac{2\sqrt{2}h_{L,K}\left(\left(\begin{smallmatrix}1\\ T\end{smallmatrix}\right)\right)}{\sqrt{8-24LT+\left(4-48K\right)T^2-4LT^3+\left(3L^2+8K\right)T^4}}\label{eqn_ratios_VF_dLK}
		\end{align}
		which is well defined and positive for $T\in\mathrm{dom}\left(\mathcal{H}_{L,K}\right)$ small enough. Note that we might have zeros in the denominator of \eqref{eqn_ratios_VF_dLK}, hence the restriction for the Euclidean norm of $T$. Now suppose that there exists a maximal integral curve $\gamma=(\gamma_L,\gamma_K)^T:I\to\mathbb{R}^2\setminus\{\mathcal{V}=0\}$ of $\mathcal{V}|_{\mathbb{R}^2\setminus\{\mathcal{V}=0\}}$, such that at least two quartic GPSR curves, and thus also the corresponding polynomials, associated to two points in the image of $\gamma$ are not equivalent. Then for any fixed $w\in\gamma(I)$ there exists $\varepsilon>0$, such that all polynomials corresponding to elements in
			\begin{equation}\label{eqn_varepsilon_equal_pf}
				\gamma\left(\left(\gamma^{-1}(w)-\varepsilon,\gamma^{-1}(w)+\varepsilon\right)\right)
			\end{equation}
		are equivalent. This follows from
			\begin{equation*}
				\gamma(w)=:\left(\begin{smallmatrix}
				L_w\\ K_w
				\end{smallmatrix}
				\right)\not\in\{\mathcal{V}=0\}
			\end{equation*}
		and equation \eqref{eqn_ratios_VF_dLK} which shows that the described ratios are locally positive and boun\-ded, and thus implies that there exists an open interval
			\begin{equation*}
				I^{\mathcal{H}_{L_w,K_w}}_w\subset \mathrm{dom}\left(\mathcal{H}_{L_w,K_w}\right),\quad 0\in I^{\mathcal{H}_{L_w,K_w}}_w,
			\end{equation*}
		(recall that being a hyperbolic point is an open condition and, hence, $\mathrm{dom}\left(\mathcal{H}_{L_w,K_w}\right)$ is in all cases an open interval) such that with $L(0)=L_w$, $K(0)=K_w$,
			\begin{equation*}
				\gamma(w)\in\left.
				\left\{\left(\begin{matrix}
				L\left(T\right)\\ K\left(T\right)
				\end{matrix}
				\right)\ \right|\ T\in I^{\mathcal{H}_{L_w,K_w}}_w\right\}\subset \gamma(I).
			\end{equation*}
		Since the map $T\mapsto(L(T),K(T))^T$ is smooth and, by the assumption $(L_w,K_w)^T\not\in\{\mathcal{V}=0\}$, non-constant locally around $T=0$ it follows that the set
			\begin{equation*}
				\left.
				\left\{\left(\begin{matrix}
				L\left(T\right)\\ K\left(T\right)
				\end{matrix}
				\right)\ \right|\ T\in I^{\mathcal{H}_{L_w,K_w}}_w\right\}\subset\gamma(I)
			\end{equation*}
		contains an open neighbourhood of $\gamma (w)$ in the subspace topology of the submanifold $\gamma(I)\subset\mathbb{R}^2$. This is precisely the reason why we consider the restriction $\mathcal{V}|_{\mathbb{R}^2\setminus\{\mathcal{V}=0\}}$. We can thus choose $\varepsilon>0$ as in \eqref{eqn_varepsilon_equal_pf}. This in particular implies that we can choose a maximal open interval $I_w\subset I$, $w\in I_w$, such that all polynomials corresponding to points in $\gamma\left(I_w\right)$ are equivalent and for any $\overline{w}\in (\partial I_w)\cap I$, which is by assumption not empty, $h_{\gamma(w)}$ and $h_{\gamma\left(\overline{w}\right)}$ are not equivalent. Also by assumption we have
			\begin{equation*}
				\gamma\left(\overline{w}\right)=:\left(\begin{smallmatrix}
				L_{\overline{w}}\\ K_{\overline{w}}
				\end{smallmatrix}
				\right)\in\mathbb{R}^2\setminus\{\mathcal{V}=0\},
			\end{equation*}
		hence we can use the same procedure for $\overline{w}$ as we used for $w$ and find that there exists a maximal open interval $I_{\overline{w}}\subset I$, $\overline{w}\in I_{\overline{w}}$, such that all polynomials corresponding to elements in $\gamma\left(I_{\overline{w}}\right)$ are equivalent. The constructed intervals $I_w$ and $I_{\overline{w}}$ are both open, and since $\overline{w}\in (\partial I_w)\cap I$ it follows that $I_w\cap I_{\overline{w}}\ne \emptyset$. But this implies that the polynomials corresponding to $\gamma(w)$ and $\gamma(\overline{w})$ are equivalent, which is a contradiction. Summarising, this proves the claim that for all maximally extended integral curves $\gamma:I\to\mathbb{R}^2\setminus\{\mathcal{V}=0\}$ of the restricted vector field $\mathcal{V}|_{\mathbb{R}^2\setminus\{\mathcal{V}=0\}}$ the corresponding polynomials $h_{\gamma_L(t),\gamma_K(t)}$ defined in \eqref{eqn_quadpoly_integralcurves} and the corresponding maximal connected quartic GPSR curves $\mathcal{H}_{\gamma_L(t),\gamma_K(t)}$ are equivalent for all $t\in I$. By construction of $\gamma$, all transformations in the image of $\gamma$ must be of the form $A(p)$ \eqref{eqn_A_explicit} for some $p\in\mathcal{H}_{L,K}$. Since the image of $\gamma$ is connected and $\mathrm{dom}(\mathcal{H}_{L,K})$ and $\mathcal{H}_{L,K}$ are diffeomorphic via the central projection, we conclude that the image of $\gamma$ is indeed contained in the image of $T\mapsto\left(\begin{smallmatrix} L(T)\\ K(T) \end{smallmatrix}\right)$, $T\in\mathrm{dom}(\mathcal{H}_{L,K})$.
		
		Observing the complexity of the formulas \eqref{eqn_L_T_pmoving} and \eqref{eqn_K_T_pmoving}, the above discussion suggests that it might be easier to be concerned with properties of the vector field $\mathcal{V}|_{\mathbb{R}^2\setminus\{\mathcal{V}=0\}}$ and its integral curves in order to find the desired classification result instead of studying the equations \eqref{eqn_L_T_pmoving} and \eqref{eqn_K_T_pmoving} directly. This is precisely what we will do from this point on in the proof of this theorem. Note that we have not shown that the set of maximal integral curves of $\mathcal{V}|_{\mathbb{R}^2\setminus\{\mathcal{V}=0\}}$ is in one-to-one correspondence with equivalence classes of polynomials, but rather that for each $h_{L,K}$ with $(L,K)^T\in\mathbb{R}^2\setminus\{\mathcal{V}=0\}$ as in \eqref{eqn_quartic_curves_std_form} there exists at least one maximal integral curve $\gamma:I\to\mathbb{R}^2\setminus\{\mathcal{V}=0\}$ of $\mathcal{V}|_{\mathbb{R}^2\setminus\{\mathcal{V}=0\}}$, such that each polynomial corresponding to a point in $\gamma(I)$ is equivalent to $h_{L,K}$. Note that since we can assume that $h_{L,K}$ corresponds to a point in $\gamma(I)$ itself, we also get that the corresponding maximal connected quartic GPSR curves are equivalent. This leaves us with the task of checking which maximal integral curves of $\mathcal{V}|_{\mathbb{R}^2\setminus\{\mathcal{V}=0\}}$ do contain points corresponding to closed connected quartic GPSR curves, and then checking if pairwise different maximal integral curves might still contain points corresponding to equivalent closed connected quartic GPSR curves. The quartic GPSR curves corresponding to points in $\{\mathcal{V}=0\}$ need to be treated as well.
		
		Lemma \ref{lem_CC_implies_precompact} implies that for $\mathcal{H}_{L,K}$ to be closed, it is a necessary requirement that the function
			\begin{equation}\label{eqn_fLK}
				f_{L,K}(t):=h_{L,K}\left(\left(\begin{smallmatrix}
				1\\ t
				\end{smallmatrix}
				\right)\right)=1-t^2+Lt^3+Kt^4
			\end{equation}
		has at least one positive and one negative real zero in $t$. This follows from the fact that the set of these zeros must coincide with $\partial\mathrm{dom}\left(\mathcal{H}_{L,K}\right)$ since otherwise the connected component of $\{h_{L,K}=1\}$ that contains the point $\left(\begin{smallmatrix}
		x\\y
		\end{smallmatrix}
		\right)=\left(\begin{smallmatrix}
		1\\ 0
		\end{smallmatrix}
		\right)$ would not coincide with $\mathcal{H}_{L,K}$, which would in turn not be closed. Recall that in the cubic case, that is for CCPSR curves, it turned out that $h\left(\left(\begin{smallmatrix} 1\\ t \end{smallmatrix} \right)\right)=0$ having at least one negative and one positive solution was also a sufficient condition for a connected PSR manifold in standard form to be closed, cf. \cite[Thm.\,1.1, Cor.\,4.2]{L2}. This is not true in the quartic case as we will see in this proof. For quartic curves, Thm. \ref{thm_incomplete_quartics} \hyperref[eqn_incomplete_qGPSRcurves_class_c]{c)} and all curve in Thm. \ref{thm_incomplete_quartics} \hyperref[eqn_incomplete_qGPSRcurves_class_b]{b)} with $K\in\left(-\frac{25}{72},0\right)$ are precisely the cases, up to equivalence, where this statement does not hold true. For the following studies we will frequently need the formulas for $f_{L,K}(t)$ and its first and second derivatives
			\begin{align*}
				\dot{f}_{L,K}(t)&=-2t+3Lt^2+4Kt^3,\\
				\ddot{f}_{L,K}(t)&=-2+6Lt+12Kt^2.
			\end{align*}
		
		Consider first $L=0$. Then $f_{0,K}(t)=1-t^2+Kt^4$. For $K\ne 0$,
			\begin{equation}\label{eqn_f0K_realroots}
				f_{0,K}(t)=0\quad\Leftrightarrow\quad t^2=\frac{1}{2K}\pm\sqrt{\frac{1-4K}{4K^2}}.
			\end{equation}
		This shows that $f_{0,K}(t)$ has no real zeros for $K>\frac{1}{4}$. It follows that for all $K>\frac{1}{4}$, $f_{L,K}(t)$ has no positive real zero for $L>0$ and no negative real zero for $L<0$. Hence, $K\leq\frac{1}{4}$ (see Figure \ref{fig_underKonequarterarea}) is a necessary requirement for $\mathcal{H}_{L,K}$ to be a quartic CCGPSR curve.
			\begin{figure}[H]%
				\centering%
				\includegraphics[scale=0.2]{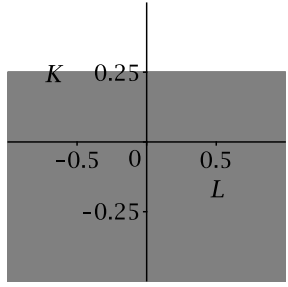}%
				\caption{$\left\{K\leq\frac{1}{4}\right\}\subset\mathbb{R}^2$ marked in grey.}\label{fig_underKonequarterarea}%
			\end{figure}
		
		Next, consider $K=\frac{1}{4}$. For that specific value of $K$, $f_{0,\frac{1}{4}}(t)=0$ if and only if $t=\pm \sqrt{2}$. Since
			\begin{align*}
				\dot{f}_{0,\frac{1}{4}}\left(\pm\sqrt{2}\right)=0,\quad
				\ddot{f}_{0,\frac{1}{4}}\left(\pm\sqrt{2}\right)=4>0,
			\end{align*}
		it follows that $t=\sqrt{2}$ and $t=-\sqrt{2}$ are both double zeros and local minima. Hence, for $L>0$ we have $f_{L,\frac{1}{4}}(t)>0$ for all $t>0$ and for $L<0$ we have $f_{L,\frac{1}{4}}(t)>0$ for all $t<0$. This shows that $f_{L,\frac{1}{4}}(t)$ has a positive and a negative real zero if and only if $L=0$.
		
		Now consider $K<\frac{1}{4}$. For $L=0$ we have shown above that $f_{0,K}(t)$ has at least one positive and one negative real zero. To analyse the cases $L\ne 0$ we will study the (possibly complex) zeros of $\dot{f}_{L,K}(t)$. We will without loss of generality assume that $L>0$, since $h_{L,K}$ and $h_{-L,K}$ are equivalent via $y\mapsto -y$. Then
			\begin{equation}\label{eqn_fLK_loc_extrama}
				\dot{f}_{L,K}(t)=0\quad\Leftrightarrow\quad t=0\ \text{ or }\left\{\begin{tabular}{rl}
				$t=\frac{2}{3L}$, & $K=0$,\\
				$t=-\frac{3L}{8K}\pm\frac{1}{8K}\sqrt{9L^2+32K}$, & $K\ne 0$.
				\end{tabular}\right.
			\end{equation}
		We want to stress that for $K\ne 0$ the latter two of the above zeros of $\dot{f}_{L,K}(t)$ might not be real. See the following figures \eqref{fig_fLKexample_1}, \eqref{fig_fLKexample_2}, and \eqref{fig_fLKexample_3} for plots of the functions $f_{L,K}(t)$ for specific values of $L$ and $K$.
			\begin{figure}[H]
				\centering
					\begin{subfigure}[h]{0.3\linewidth}
						\centering
						\includegraphics[scale=0.16]{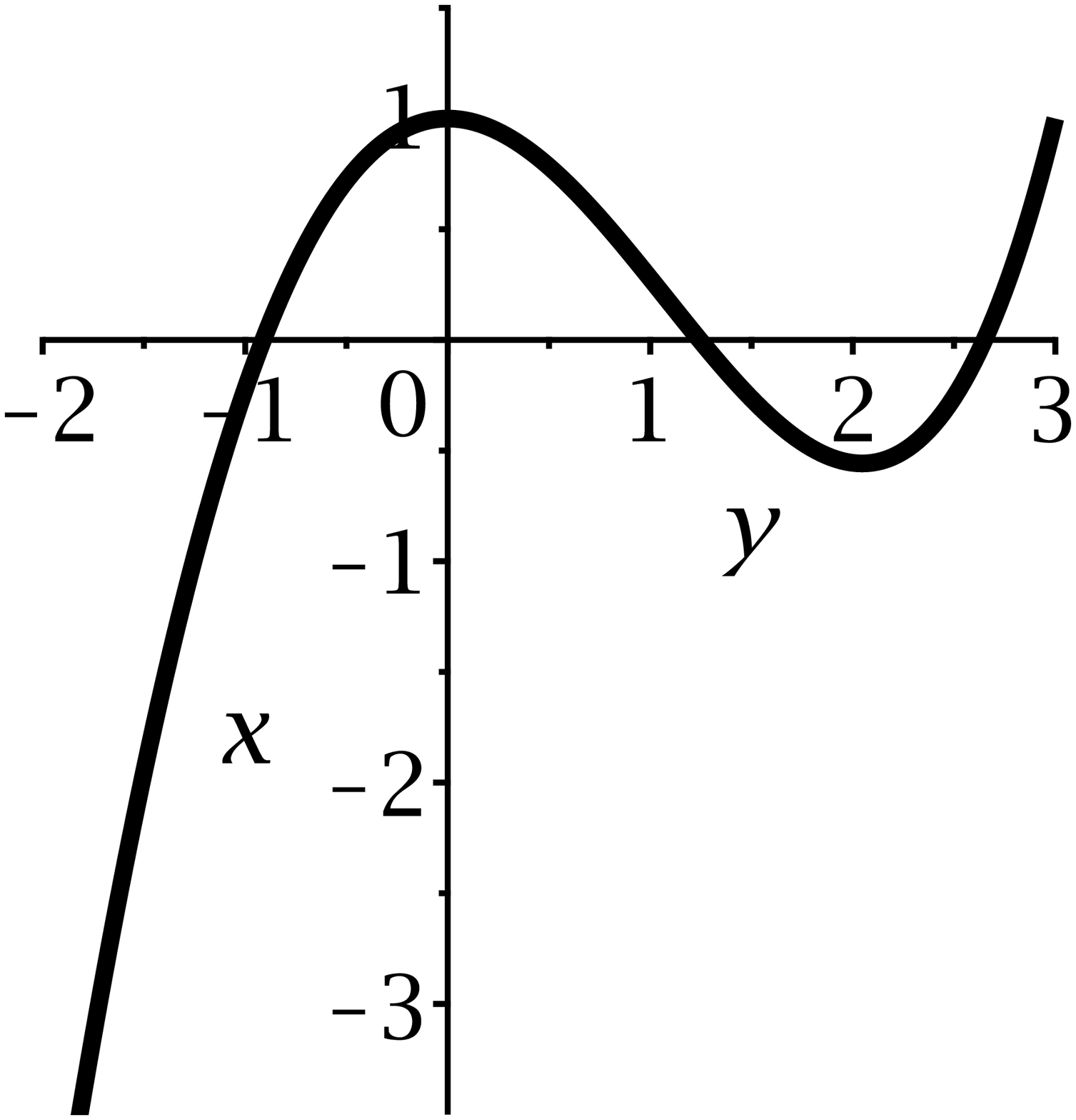}
						\caption{$f_{\frac{1}{4},\frac{1}{36}}(t),\ t\in\left(-\frac{3}{2},3\right)$.}\label{fig_fLKexample_1}
					\end{subfigure}
					\begin{subfigure}[h]{0.3\linewidth}
						\centering
						\includegraphics[scale=0.16]{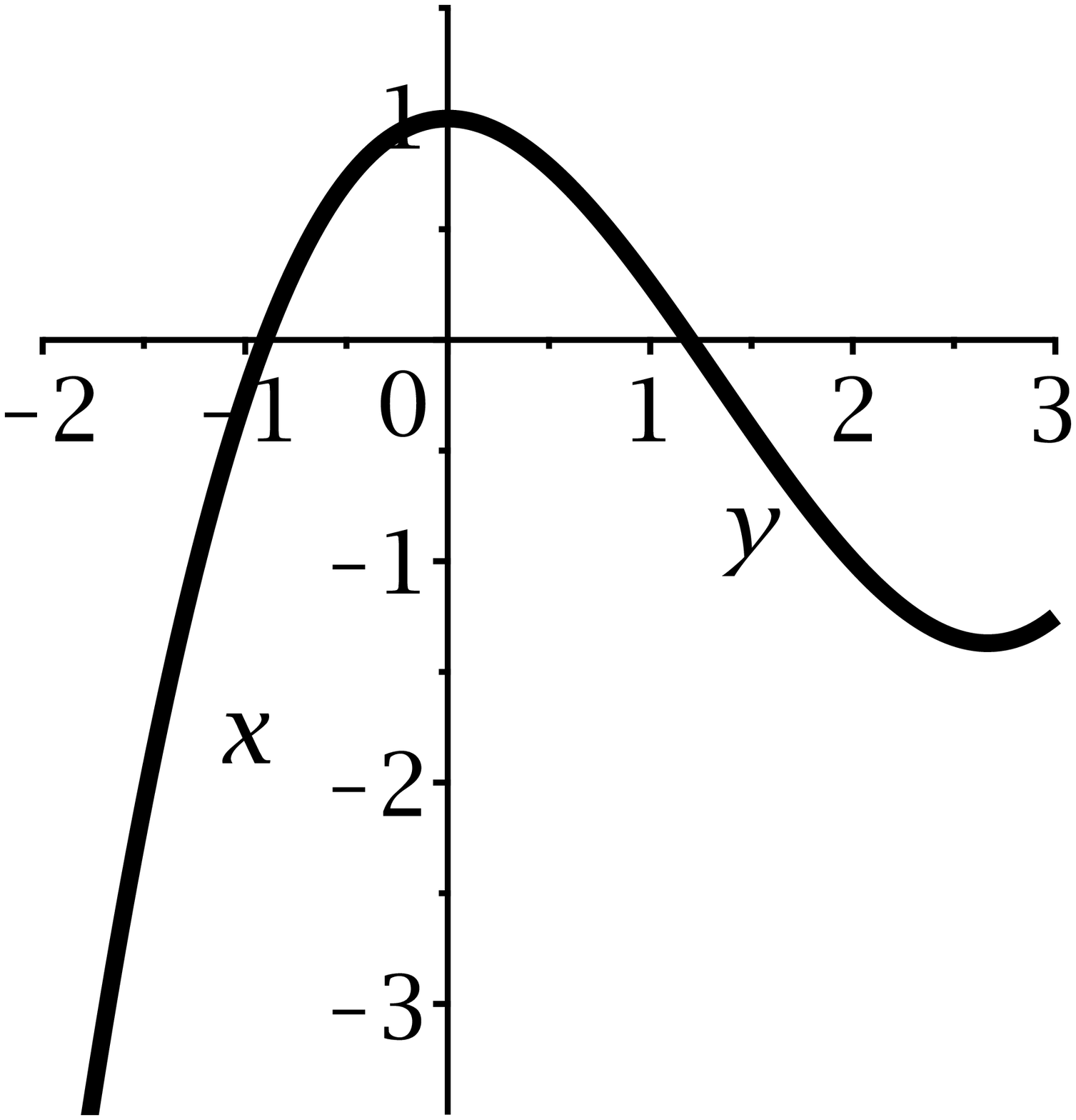}
						\caption{$f_{\frac{1}{4},0}(t),\ t\in\left(-\frac{3}{2},3\right)$.}\label{fig_fLKexample_2}
					\end{subfigure}
					\begin{subfigure}[h]{0.3\linewidth}
						\centering
						\includegraphics[scale=0.16]{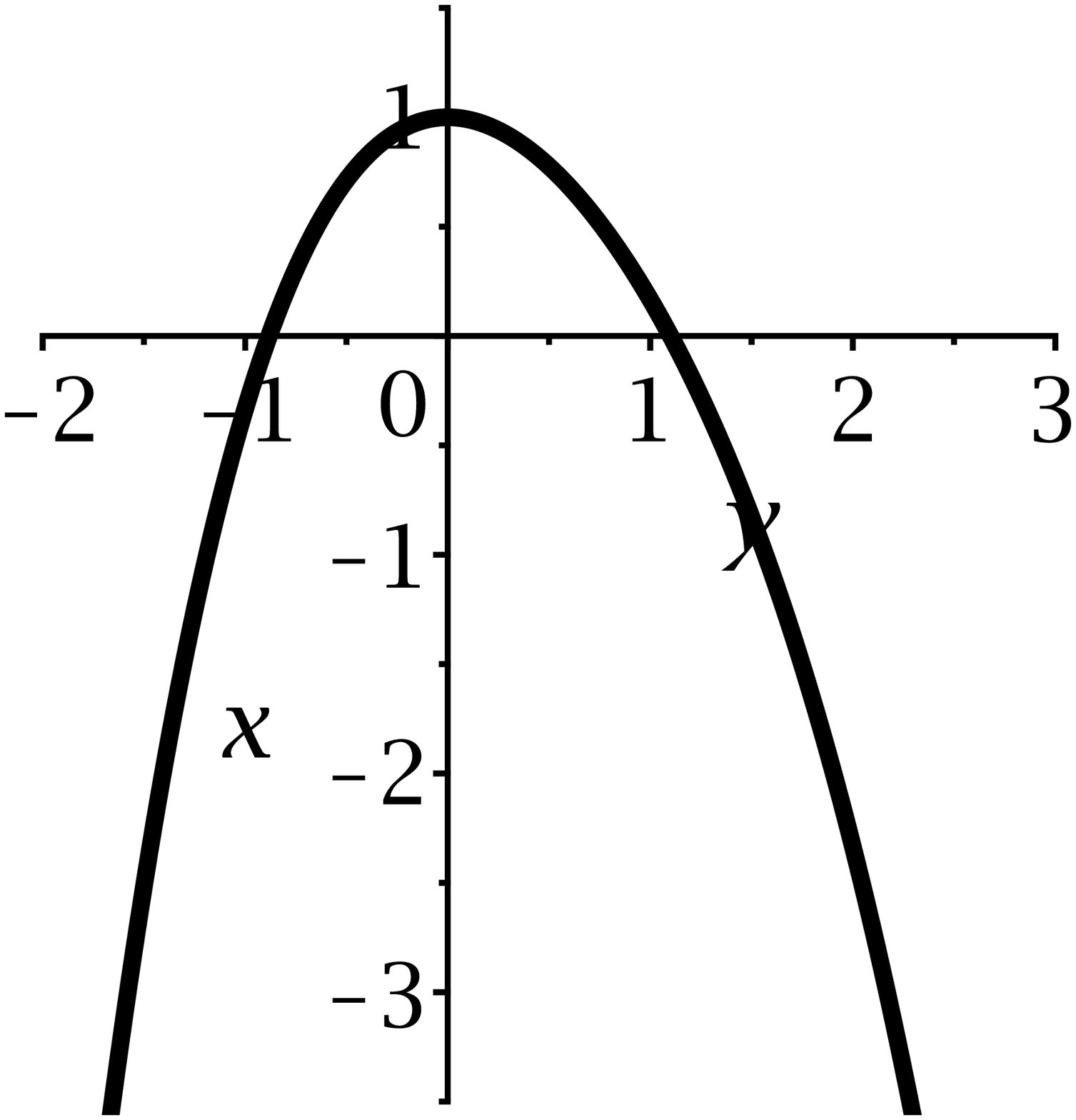}
						\caption{$f_{\frac{1}{4},-\frac{1}{12}}(t),\ t\in\left(-\frac{3}{2},3\right)$.}\label{fig_fLKexample_3}
					\end{subfigure}%
				\caption{Example plots of $f_{L,K}(t)$.}
			\end{figure}
		\noindent
		We now use Maple to symbolically solve the system of equations
			\begin{equation}\label{eqn_upper_bdr_system_eqns}
				f_{L,K}(t)=0,\quad \dot{f}_{L,K}(t)=0,
			\end{equation}
		for the variable $L$ and, as one of the solutions under the restriction
			\begin{equation}\label{eqn_t_m_def}
				t=t_m:=-\frac{3L}{8K}+\frac{1}{8K}\sqrt{9L^2+32K},
			\end{equation}
		which is one solution of $\dot{f}_{L,K}(t)=0$, we obtain for $K\ne 0$
			\begin{equation}\label{eqn_upper_bdr}
				L=\frac{\sqrt{2}}{3\sqrt{3}}\sqrt{1-36K+\sqrt{(1+12K)^3}}=:\mathbf{u}(K).
			\end{equation}
		We will consider $\mathbf{u}$ as a function on the interval $\left(-\frac{1}{12},\frac{1}{4}\right)$,
			\begin{equation}\label{eqn_bfu_def}
				\mathbf{u}:\left(-\frac{1}{12},\frac{1}{4}\right)\to\mathbb{R},
			\end{equation}
		where we note that for $K=0$, $L=\mathbf{u}(0)=\frac{2}{3\sqrt{3}}$ and $t=\frac{2}{3L}=\sqrt{3}$ solve \eqref{eqn_upper_bdr_system_eqns}.
		Before explaining the reason why we choose this particular lower bound for the interval $\left(-\frac{1}{12},\frac{1}{4}\right)$ (see equation \eqref{eqn_mathbfu_limits}), we will analyse $\mathbf{u}(K)$ further. We will show that $\mathbf{u}(K)>0$ for all $K\in\left(-\frac{1}{12},\frac{1}{4}\right)$. We obtain
			\begin{equation*}
				1-36K+\sqrt{(1+12K)^3}=0\quad\Rightarrow \quad (1+12K)^3-(1-36K)^2=0\quad\Rightarrow\quad K\in\left\{0,\frac{1}{4}\right\}.
			\end{equation*}
		So the only candidate for a solution of $\mathbf{u}(K)=0$ that is contained in the set $\left(-\frac{1}{12},\frac{1}{4}\right)$ is $K=0$. But $\mathbf{u}(0)=\frac{2}{3\sqrt{3}}$. We conclude that $\mathbf{u}(K)>0$ for all $K\in\left(-\frac{1}{12},\frac{1}{4}\right)$ as claimed. Note that for $K=0$, $\mathbf{u}(0)=\frac{2}{3\sqrt{3}}$ coincides with the unique real solution of $\dot{f}_{L,0}\left(t_m|_{K=0}\right)=0$ in $L$, where ${t_m}|_{K=0}:=\frac{2}{3L}$, cf. \eqref{eqn_fLK_loc_extrama}. Furthermore, for $K=\frac{1}{4}$, $\mathbf{u}\left(\frac{1}{4}\right)=0$, so the point $\left(\mathbf{u}\left(\frac{1}{4}\right),\frac{1}{4}\right)^T\in\mathbb{R}^2$ coincides with one of the fixed points of $\mathcal{V}$. Next we will show that the graph of $\mathbf{u}$ coincides with the image of a maximal integral curve of the vector field $\mathcal{V}|_{\mathbb{R}^2\setminus\{\mathcal{V}=0\}}$ defined in \eqref{eqn_VFqCCGPSRc_def} (see Figure \ref{fig_u_image} for a plot of the graph of $\mathbf{u}$). We know that $\mathcal{V}$ has no zeroes in the set
			\begin{equation}\label{eqn_image_upper_u}
				\left\{\left(\begin{smallmatrix}
				\mathbf{u}(K)\\ K
				\end{smallmatrix}
				\right)\ \left|\ K\in  \left(-\frac{1}{12},\frac{1}{4}\right)\right.\right\}\subset\mathbb{R}^2,
			\end{equation}
		see equation \eqref{eqn_mathcalV_zeros}. 
			\begin{figure}[H]%
				\centering%
				\includegraphics[scale=0.2]{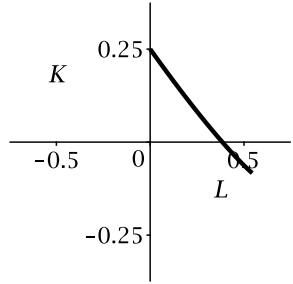}%
				\caption{The graph of $\mathbf{u}$ embedded in $\mathbb{R}^2$ as in \eqref{eqn_image_upper_u}.}\label{fig_u_image}%
			\end{figure}
		\noindent
		Furthermore $\D K(\mathcal{V})=L\left(6K+\frac{1}{2}\right)$ does not vanish for $L\ne 0$ and $K\ne-\frac{1}{12}$, so in particular it does not vanish on the graph of $\mathbf{u}$. We can extend $\mathbf{u}$ continuously to $\left[-\frac{1}{12},\frac{1}{4}\right]$, which shows that the graph of $\mathbf{u}$ is precompact in $\mathbb{R}^2$. One now verifies 
			\begin{equation}\label{eqn_partial_u}
				\partial_K\mathbf{u}(K)=\frac{6\left(-2+\sqrt{1+12K} \right)}{\sqrt{6-216K+6\sqrt{(1+12K)^3}}}= \left.\frac{\D L(\mathcal{V})}{\D K(\mathcal{V})}\right|_{\left(\begin{smallmatrix}
				L\\ K
				\end{smallmatrix}
				\right)=\left(\begin{smallmatrix}
				\mathbf{u}(K)\\ K
				\end{smallmatrix}
				\right)}.
			\end{equation}
		This shows that the image of $K\mapsto(\mathbf{u}\left(K\right),K)^T$, $K\in\left(-\frac{1}{12},\frac{1}{4}\right)$, is contained in the image of a maximal integral curve of $\mathcal{V}|_{\mathbb{R}^2\setminus\{\mathcal{V}=0\}}$. Since
			\begin{equation}\label{eqn_mathbfu_limits}
				\left(\begin{smallmatrix}
				\mathbf{u}\left(\frac{1}{4}\right)\\ \frac{1}{4}
				\end{smallmatrix}
				\right)=\left(\begin{smallmatrix}
				0\\ \frac{1}{4}
				\end{smallmatrix}
				\right)\in\{\mathcal{V}=0\}, \quad \left(\begin{smallmatrix}
				\mathbf{u}\left(-\frac{1}{12}\right)\\ -\frac{1}{12}
				\end{smallmatrix}
				\right)=\left(\begin{smallmatrix}
				\frac{2\sqrt{2}}{3\sqrt{3}}\\ -\frac{1}{12}
				\end{smallmatrix}
				\right)\in\{\mathcal{V}=0\},
			\end{equation}
		we conclude that said image coincides with a maximal integral curve of $\mathcal{V}|_{\mathbb{R}^2\setminus\{\mathcal{V}=0\}}$. Note that it contains in particular the point $\left(\frac{2}{3\sqrt{3}},0\right)^T\in \mathbb{R}^2$, which corresponds to the polynomial \hyperref[eqn_qCCPSRcurves_class_c]{c)}. Equation \eqref{eqn_mathbfu_limits} also explains the lower bound $-\frac{1}{12}$ of the domain of definition of $\mathbf{u}$. Further note that $\mathbf{u}$ is a convex function, which one can quickly verify by determining $\partial^2_K\mathbf{u}$.
		
		Next we will show that for all $K\in\left(-\frac{1}{12},\frac{1}{4}\right)$ and all $\ell>0$, $\mathcal{H}_{L,K}$ corresponding to
			\begin{equation*}
				\left(\begin{smallmatrix}
				L\\ K
				\end{smallmatrix}
				\right)=
				\left(\begin{smallmatrix}
				\mathbf{u}(K)+\ell\\ K
				\end{smallmatrix}
				\right)\in\mathbb{R}^2\setminus\{\mathcal{V}=0\}
			\end{equation*}
		is not closed, cf. Figure \ref{fig_area_rigth_of_u}.
			\begin{figure}[H]%
				\centering%
				\includegraphics[scale=0.2]{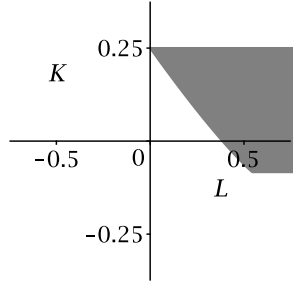}%
				\caption{Part of the set $\left\{(\mathbf{u}(K)+\ell,K)^T,\ K\in\left(-\frac{1}{12},\frac{1}{4}\right),\ \ell>0\right\}$ marked in grey.}\label{fig_area_rigth_of_u}%
			\end{figure}
		\noindent
		For $K=0$, we need to analyse the functions
			\begin{equation*}
				f_{\frac{2}{3\sqrt{3}}+\ell,0}(t)=1-t^2+\left(\frac{2}{3\sqrt{3}}+\ell\right)t^3
			\end{equation*}
		for $\ell>0$. Observe that $f_{\frac{2}{3\sqrt{3}},0}(t)=1-t^2+ \frac{2}{3\sqrt{3}} t^3$ has precisely one positive zero $t=\sqrt{3}$. Since $\dot{f}_{\frac{2}{3\sqrt{3}},0}\left(\sqrt{3}\right)=0$ and $\ddot{f}_{\frac{2}{3\sqrt{3}},0}\left(\sqrt{3}\right)=2$, $f_{\frac{2}{3\sqrt{3}},0}(t)$ has a local minimum at $t=\sqrt{3}$. Furthermore
			\begin{equation*}
				f_{\frac{2}{3\sqrt{3}}+\ell,0}(t)>f_{\frac{2}{3\sqrt{3}},0}(t)
			\end{equation*}
		for all $\ell>0$ and all $t>0$. We conclude that for all $\ell>0$, $f_{\frac{2}{3\sqrt{3}}+\ell,0}(t)$ has no positive real zero and, hence, $\mathcal{H}_{\frac{2}{3\sqrt{3}}+\ell,0}$ is not closed. See Figure \ref{fig_rightofu_example_plotf} for an example plot of a function of the form $f_{\frac{2}{3\sqrt{3}}+\ell,0}(t)$. 
			\begin{figure}[H]
				\centering
				\includegraphics[scale=0.15]{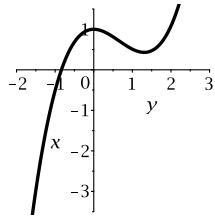}
				\caption{A plot of $f_{\mathbf{u}(K)+\ell,K}(t)$, $t\in(-2,4)$, for $K=0$ and $\ell=\frac{1}{8}$.}\label{fig_rightofu_example_plotf}
			\end{figure}
		
		For $K\ne 0$, $t_m$, defined in \eqref{eqn_t_m_def} above, is positive whenever $L>0$ and $9L^2+32K\geq 0$. For
			\begin{equation}\label{eqn_K_ne_0_uKl}
				\left(\begin{smallmatrix}
						L\\ K
						\end{smallmatrix}
						\right)=
						\left(\begin{smallmatrix}
						\mathbf{u}(K)+\ell\\ K
						\end{smallmatrix}
						\right),\quad K\in\left(-\frac{1}{12},\frac{1}{4}\right),\quad\ell>0,
			\end{equation}
		we have
			\begin{equation}\label{eqn_existence_two_extrama_GPSRcurves}
				9L^2+32K=\underbrace{18\ell  \mathbf{u}(K)+9\ell^2}_{>0}+\underbrace{\frac{2}{3}\left(1+12K\right)}_{>0}+\underbrace{\frac{2}{3}\sqrt{(1+12K)^3}}_{>0}>0.
			\end{equation}
		Hence, $t_m$ is real and positive along the considered points \eqref{eqn_K_ne_0_uKl} for all $K\in\left(-\frac{1}{12},\frac{1}{4}\right)\setminus\{0\}$ and all $\ell>0$.
		
		For all $K\in\left(0,\frac{1}{4}\right)$ and all $\ell>0$, the function $f_{\mathbf{u}(K)+\ell,K}(t)=1-t^2+\left(\mathbf{u}(K)+\ell\right)t^3+Kt^4$ has a local maximum at $t=0$ and diverges to $+\infty$ for $t\to\pm\infty$. Hence, $f_{\mathbf{u}(K)+\ell,K}(t)$ must have a local minimum at $t=t_m$. For $\ell=0$, $t_m$ is by construction also a zero of $f_{\mathbf{u}(K),K}(t)$. Since for all $t>0$
			\begin{equation*}
				f_{\mathbf{u}(K)+\ell,K}(t)>f_{\mathbf{u}(K),K}(t)
			\end{equation*}
		and, implied by the construction of $\mathbf{u}$, cf. \eqref{eqn_upper_bdr_system_eqns} and \eqref{eqn_t_m_def}, $f_{\mathbf{u}(K),K}(t)$ has precisely one positive real zero, we conclude that for all $\ell>0$, $f_{\mathbf{u}(K)+\ell,K}(t)$ has no positive real zero. Hence, $\mathcal{H} _{\mathbf{u}(K)+\ell,K}$ is not closed.
		
		Next consider $K\in\left(-\frac{1}{12},0\right)$. In these cases $f_{\mathbf{u}(K)+\ell,K}(t)$ does always have at least one negative and one positive real zero for all $\ell>0$ since its highest order monomial $t^4$ has a negative coefficient. Thus, in order to prove the claim that $\mathcal{H} _{\mathbf{u}(K)+\ell,K}$ is not closed for all $\ell>0$, we need to check that there exists at least one point $\widetilde{t}$ in the connected component containing $t=0$ of the set
			\begin{equation*}
				\left.\left\{t\in\mathbb{R}\ \right|\ f_{\mathbf{u}(K)+\ell,K}(t)>0\right\},
			\end{equation*}
		such that $\left(\begin{smallmatrix}
		1\\ \widetilde{t}
		\end{smallmatrix}
		\right)$ is not a hyperbolic point of the corresponding quartic polynomial $h_{\mathbf{u}(K)+\ell,K}$. To do so recall that we have already shown that along $L=\mathbf{u}(K)+\ell$ and $K\in\left(-\frac{1}{12},0\right)\subset\left(-\frac{1}{12},\frac{1}{4}\right)\setminus\{0\}$, the term $9L^2+32K$ is positive, see \eqref{eqn_existence_two_extrama_GPSRcurves}. Hence, \eqref{eqn_fLK_loc_extrama} implies that for all $\ell>0$, $f_{\mathbf{u}(K)+\ell,K}(t)$ has exactly three local extrema, namely at $t=0$, $t_m=-\frac{3\left(\mathbf{u}(K)+\ell\right)}{8K}+\frac{1}{8K}\sqrt{9\left(\mathbf{u}(K)+\ell\right)^2+32K}$, and with
			\begin{equation}\label{eqn_t_M_def}
				t_M:=-\frac{3L}{8K}-\frac{1}{8K}\sqrt{9L^2+32K},
			\end{equation}
		at $t_M=-\frac{3\left(\mathbf{u}(K)+\ell\right)}{8K}-\frac{1}{8K}\sqrt{9\left(\mathbf{u}(K)+\ell\right)^2+32K}$. 
		Since $K<0$, it follows that both $t_m$ and $t_M$ are positive and \eqref{eqn_existence_two_extrama_GPSRcurves} implies the strict inequality $t_m<t_M$. For all $K\in\left(-\frac{1}{12},\frac{1}{4}\right)$ and all $\ell>0$, the function $f_{\mathbf{u}(K)+\ell,K}(t)$ is a quartic polynomial in $t$, which implies that it has at most three distinct local extrema. Hence, we have shown that $f_{\mathbf{u}(K)+\ell,K}(t)$ has precisely three extrema at the distinct points $t=0$, $t_m$, and $t_M$, $f_{\mathbf{u}(K)+\ell,K}(t)\to -\infty$ for $t\to\pm\infty$, and furthermore that $f_{\mathbf{u}(K)+\ell,K}(t)$ has a local maximum at $t=0$. We deduce that $f_{\mathbf{u}(K)+\ell,K}(t)$ will always have a local minimum at $t_m$ and a local maximum at $t_M$. We will now show that for all $K\in\left(-\frac{1}{12},0\right)$ and all $\ell>0$, $t_m$ is contained in the connected component that contains $t=0$ of $\left.\left\{t\in\mathbb{R}\ \right|\ f_{\mathbf{u}(K)+\ell,K}(t)>0\right\}$ and that $\left(\begin{smallmatrix}
		1\\ t_m
		\end{smallmatrix}
		\right)$ is indeed not a hyperbolic point of $h_{\mathbf{u}(K)+\ell,K}$. To show the first statement, it suffices to show that $f_{\mathbf{u}(K)+\ell,K}\left(t_m\right)>0$, since we have seen that for all $K\in\left(-\frac{1}{12},0\right)$ and all $\ell>0$, $f_{\mathbf{u}(K)+\ell,K}(t)$ has a local maximum at $t=0$, a local minimum at $t_m>0$, and no extremal point at any $t\in\left(0,t_m\right)$. We view $f_{\mathbf{u}(K)+\ell,K}\left(t_m\right)$ as a function in the variables $K$ and $\ell$ and calculate
			\begin{align*}
				&\partial_\ell\left(f_{\mathbf{u}(K)+\ell,K}\left(t_m\right)\right)=\frac{1}{13824K^3}\left(-9\ell -\sqrt{6-216K+6\sqrt{(1+12K)^3}}\right.\\
				&\left.\quad +  \sqrt{6+72K+6\sqrt{(1+12K)^3}+18\ell\sqrt{6-216K+6\sqrt{(1+12K)^3}}+81\ell^2}\right)^3.
			\end{align*}
		At $\ell=0$, 
			\begin{align*}
				&\left.\partial_\ell\left(f_{\mathbf{u}(K)+\ell,K}\left(t_m\right)\right)\right|_{\ell=0}\\ 
				&=\frac{1}{13824K^3}\left(-\sqrt{6-216K+6\sqrt{(1+12K)^3}}+ \sqrt{6+72K+6\sqrt{(1+12K)^3}}\right)>0,
			\end{align*}
		which is easily seen for all $K\in\left(-\frac{1}{12},0\right)$. Suppose that there exists $K\in\left(-\frac{1}{12},0\right)$ and $\ell>0$, such that $\partial_\ell\left(f_{\mathbf{u}(K)+\ell,K}\left(t_m\right)\right)=0$. Then
			\begin{align*}
				&9\ell +\sqrt{6-216K+6\sqrt{(1+12K)^3}}\\
				&=\sqrt{6+72K+6\sqrt{(1+12K)^3}+18\ell\sqrt{6-216K+6\sqrt{(1+12K)^3}}+81\ell^2}\\
				\Rightarrow\quad & K=0,
			\end{align*}
		obtained by taking the square of both sides of the first equation. This contradicts $K\in\left(-\frac{1}{12},0\right)$. We conclude that for all $K\in\left(-\frac{1}{12},0\right)$ and all $\ell>0$, $\partial_\ell\left(f_{\mathbf{u}(K)+\ell,K}\left(t_m\right)\right)>0$. Since $f_{\mathbf{u}(K),K}\left(t_m\right)=0$ by construction, this shows that $t_m$ is in fact contained in the connected component that contains $t=0$ of $\left.\left\{t\in\mathbb{R}\ \right|\ f_{\mathbf{u}(K)+\ell,K}(t)>0\right\}$ as claimed. In order to show that $\left(\begin{smallmatrix}
		1\\ t_m
		\end{smallmatrix}
		\right)$ is not a hyperbolic point of $h_{\mathbf{u}(K)+\ell,K}$ we will use Lemma \ref{lem_pullback_g}. In the one-dimensional case, i.e. in our case where $\dim\left(\mathcal{H}_{L,K}\right)=1$, the function $f$ used in Lemma \ref{lem_pullback_g} and $f_{L,K}:\mathbb{R}\to\mathbb{R}$ coincide. Using equation \eqref{eqn_pullback_metric_to_dom} yields that the pullback of the centro-affine metric at $t=t_m\in\mathrm{dom}\left(\mathcal{H}_{\mathbf{u}(K)+\ell,K}\right)$ is given by
			\begin{equation}\label{eqn_pullback_g_quartic_curves_right_of_upper}
				\left(\Phi^*g_{\mathcal{H}_{\mathbf{u}(K)+\ell,K}}\right)_{t_m}=-\frac{\ddot{f}_{\mathbf{u}(K)+\ell,K}(t_m)}{4f_{\mathbf{u}(K)+\ell,K}(t_m)}\D t^2,
			\end{equation} 
		where $t$ denotes the coordinate of $\mathrm{dom}\left(\mathcal{H}_{\mathbf{u}(K)+\ell,K}\right)$. We have shown that $f_{\mathbf{u}(K)+\ell,K}(t)$ has a local minimum at $t_m$ for all $K\in\left(-\frac{1}{12},0\right)$ and all $\ell>0$, which implies that $-\frac{1}{4}\ddot{f}_{\mathbf{u}(K)+\ell,K}(t_m)\leq 0$. Hence, $\left(\begin{smallmatrix}
		1\\ t_m
		\end{smallmatrix}
		\right)$ is not a hyperbolic point of $h_{\mathbf{u}(K)+\ell,K}$ as claimed, and we deduce that for all $K\in\left(-\frac{1}{12},0\right)$ and for all $\ell>0$, $\mathcal{H}_{\mathbf{u}(K)+\ell,K}$ is not closed in $\mathbb{R}^2$.
		
		Summarising up to this point, we have determined for each $K>-\frac{1}{12}$ a non-negative lower bound for $L$, and by equivalence of $\mathcal{H}_{L,K}$ and $\mathcal{H}_{-L,K}$, also a non-positive upper bound for $L$, such that the maximal connected quartic GPSR curve $\mathcal{H}_{L,K}$ is not a
		CCGPSR curve. These points form the set
			\begin{equation}\label{eqn_set_covered_inbetweener}
				\left\{K>\frac{1}{4}\right\}\cup\left\{K=\frac{1}{4},\ |L|>0\right\}
				\cup\left\{K\in\left(-\frac{1}{12},\frac{1}{4}\right),\ |L|>\mathbf{u}(K)\right\},
			\end{equation}
		see Figure \ref{fig_sofarcovered}.
			\begin{figure}[H]
				\centering
				\includegraphics[scale=0.2]{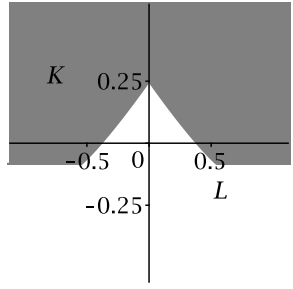}
				\caption{The set \eqref{eqn_set_covered_inbetweener} marked in grey.}\label{fig_sofarcovered}
			\end{figure}
		
		Next we will deal the case $K=-\frac{1}{12}$. It turns out that we can use the same strategy as for $K\in\left(-\frac{1}{12},0\right)$ by considering $\mathbf{u}(K)$ at the limit point $K=-\frac{1}{12}$, $\mathbf{u}\left(-\frac{1}{12}\right)=\frac{2\sqrt{2}}{3\sqrt{3}}$. For all $\ell>0$ and the corresponding function $f_{\frac{2\sqrt{2}}{3\sqrt{3}}+\ell,-\frac{1}{12}}(t)$, the points $t=0$, $t_m=\frac{3}{2}\left(\frac{2\sqrt{2}}{\sqrt{3}}+\ell-\sqrt{\frac{4\sqrt{2}}{3\sqrt{3}}\ell+\ell^2} \right)$ \eqref{eqn_t_m_def}, and $t_M=\frac{3}{2}\left(\frac{2\sqrt{2}}{\sqrt{3}}+\ell+\sqrt{\frac{4\sqrt{2}}{3\sqrt{3}}\ell+\ell^2} \right)$ \eqref{eqn_t_M_def}, still fulfil $0<t_m<t_M$ and are also still critical points of $f_{\frac{2\sqrt{2}}{3\sqrt{3}}+\ell,-\frac{1}{12}}(t)$. Also, we can show similarly as for the case $K\in\left(-\frac{1}{12},0\right)$ that
			\begin{equation*}
				\partial_\ell\left(f_{\frac{2\sqrt{2}}{3\sqrt{3}}+\ell,-\frac{1}{12}}(t_m)\right)>0
			\end{equation*}
		for all $\ell\geq 0$. Since $\left.f_{\frac{2\sqrt{2}}{3\sqrt{3}}+\ell,-\frac{1}{12}}(t_m)\right|_{\ell=0}=0$ and $t=0$, $t_m$, and $t_M$ are precisely the critical points of $f_{\frac{2\sqrt{2}}{3\sqrt{3}}+\ell,-\frac{1}{12}}(t)$ for all $\ell>0$, we conclude that for all $\ell>0$, $t_m$ is contained in the connected component of the set
			\begin{equation*}
				\left\{t\in\mathbb{R}\ \left|\ f_{\frac{2\sqrt{2}}{3\sqrt{3}}+\ell,-\frac{1}{12}}(t)>0\right\}\right.
			\end{equation*}
		that contains the point $t=0$. From equation \eqref{eqn_pullback_g_quartic_curves_right_of_upper} for the limit $K=-\frac{1}{12}$ it follows that $\left(\begin{smallmatrix}
		x\\ y
		\end{smallmatrix}
		\right)=\left(\begin{smallmatrix}
		1\\ t_m
		\end{smallmatrix}
		\right)$ is not a hyperbolic point of $h_{\frac{2\sqrt{2}}{3\sqrt{3}}+\ell,-\frac{1}{12}}$. Thus, for all $\ell>0$, $\mathcal{H}_{\frac{2\sqrt{2}}{3\sqrt{3}}+\ell,-\frac{1}{12}}$ is not closed in $\mathbb{R}^2$, i.e. not a quartic CCGPSR curve.
		
		Lastly in this sta ge of the proof we will study the case $K<-\frac{1}{12}$. For $K<0$ and $L>0$ the points $t_m$ \eqref{eqn_t_m_def} and $t_M$ \eqref{eqn_t_M_def} are both real numbers if and only if $9L^2+32K\geq 0$, and coincide if and only if
			\begin{equation*}
				9L^2+32K=0\quad \Leftrightarrow\quad L=\frac{4\sqrt{-2K}}{3}=:\mathbf{v}(K).
			\end{equation*}
		The function $\mathbf{v}$ is smooth and positive on $\{K<0\}$, but we only need to consider its restriction
			\begin{equation*}
				\mathbf{v}:\left(-\infty,-\frac{1}{12}\right)\to\mathbb{R}_{>0},
			\end{equation*}
		with smooth continuation to $\mathbf{v}\left(-\frac{1}{12}\right)=\frac{2\sqrt{2}}{3\sqrt{3}}$. Observe that $\mathbf{u}\left(-\frac{1}{12}\right)$ and $\mathbf{v}\left(-\frac{1}{12}\right)$ coincide (see Figure \ref{fig_u_v}),
			\begin{figure}[H]
				\centering
				\includegraphics[scale=0.2]{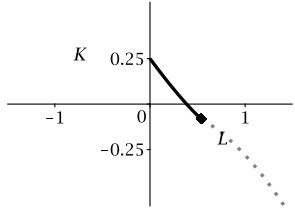}
				\caption{The graphs of $\mathbf{u}$ and $\mathbf{v}$ in $\mathbb{R}^2$, and the point $\left(L,K\right)^T=\left(\frac{2\sqrt{2}}{3\sqrt{3}},-\frac{1}{12}\right)^T$ marked with a black diamond.}\label{fig_u_v}
			\end{figure}
		\noindent
		and that the image of $(\mathbf{v}(K),K)^T$ for $K\in\left(-\infty,-\frac{1}{12}\right)$ is contained in $\mathbb{R}^2\setminus\{\mathcal{V}=0\}$. Along points of the form $(\mathbf{v}(K),K)^T\in\mathbb{R}^2$, $f_{\mathbf{v}(K),K}(t)$ has a saddle point at $t_m=t_M$. We will now show that $f_{\mathbf{v}(K),K}(t_m)>0$ for all $K<-\frac{1}{12}$. We view $t_m$ as a function depending on $L$ and $K$, and obtain along points of the form $(\mathbf{v}(K),K)^T\in\mathbb{R}^2$
			\begin{equation}\label{eqn_ftm_lowerboundary}
				f_{\mathbf{v}(K),K}(t_m)=\frac{1+12K}{12K}>0
			\end{equation}
		for all $K<-\frac{1}{12}$. Since $f_{\mathbf{v}(K),K}(t)$ is monotonously decreasing for $t>0$, this shows that $t_m$ is contained in the connected component of $\{f_{\mathbf{v}(K),K}(t)>0\}$ that contains the point $t=0$. But $t_m$ is a saddle point of $f_{\mathbf{v}(K),K}(t)$, which implies using \eqref{eqn_pullback_metric_to_dom} that
			\begin{equation*}
				\left(\Phi^*g_{\mathcal{H}_{\mathbf{v}(K),K}}\right)_{t_m}=0
			\end{equation*}
		for all $K<-\frac{1}{12}$. Hence, for all $K<-\frac{1}{12}$, $\mathcal{H}_{\mathbf{v}(K),K}$ is not closed. 
		
		Next we will show that for all $K<-\frac{1}{12}$ and all $\ell>0$, $\mathcal{H}_{\mathbf{v}(K)+\ell,K}$ is not closed. In these cases, $t_m$ \eqref{eqn_t_m_def} and $t_M$ \eqref{eqn_t_M_def} are both real numbers and $f_{\mathbf{v}(K)+\ell,K}(t)$ has a local minimum at $t_m>0$ and a local maximum at $t_M>t_m$. We will show that $f_{\mathbf{v}(K)+\ell,K}(t_m)>0$. Since $f_{\mathbf{v}(K)+\ell,K}(t)\to -\infty$ for $t\to \pm\infty$, this means that $t_m$ is an element of the connected component of $\{f_{\mathbf{v}(K)+\ell,K}>0\}$ that contains $t=0$. This implies that $\left(\begin{smallmatrix}
		1\\ t_m
		\end{smallmatrix}
		\right)$ is not a hyperbolic point of the corresponding $h_{\mathbf{v}(K)+\ell,K}$ and thus that $\mathcal{H}_{\mathbf{v}(K)+\ell,K}$ is not closed for all $K<-\frac{1}{12}$ and all $\ell>0$. We proceed similarly to the calculations for $f_{\mathbf{u}(K)+\ell,K}(t)$ with $K\in\left(-\frac{1}{12},0\right)$, $\ell>0$. Along $(L,K)^T=(\mathbf{v}(K)+\ell,K)^T$ we have
			\begin{align}
				t_m&=-\frac{3\left(\mathbf{v}(K)+\ell\right)}{8K}+\frac{1}{8K}\sqrt{9\left(\mathbf{v}(K)+\ell\right)^2+32K}\notag\\
				&=-\frac{4\sqrt{-2K}+3\ell}{8K}+\frac{1}{8K}\sqrt{24\sqrt{-2K}\ell+9\ell^2}\label{eqn_tm_rightoflowerboundary}
			\end{align}
		and
			\begin{equation*}
				\partial_\ell\left(f_{\mathbf{v}(K)+\ell,K}\left(t_m\right)\right)=\dot{f}_{\mathbf{v}(K)+\ell,K}\left(t_m\right)\cdot\partial_\ell t_m + t_m^3.
			\end{equation*}
		Hence, $t_m$ being positive, a critical point of $f_{\mathbf{v}(K)+\ell,K}\left(t\right)$, and smooth in $\ell$-dependence for $\ell>0$ implies that
			\begin{equation*}
				\partial_\ell\left(f_{\mathbf{v}(K)+\ell,K}\left(t_m\right)\right)=t_m^3>0
			\end{equation*}
		for all $\ell>0$. We have to be careful with the limit case $\ell=0$ since the first $\ell$-derivative of $t_m$ \eqref{eqn_tm_rightoflowerboundary} is easily seen to diverge as $\ell\searrow 0$ for all $K<-\frac{1}{12}$. However, \eqref{eqn_tm_rightoflowerboundary} also implies that $t_m$ can be continuously extended to $\ell=0$ for fixed $K<-\frac{1}{12}$, and hence it follows with $\left.f_{\mathbf{v}(K)+\ell,K}\left(t_m\right)\right|_{\ell=0}>0$ \eqref{eqn_ftm_lowerboundary} and $\partial_\ell\left(f_{\mathbf{v}(K)+\ell,K}\left(t_m\right)\right)>0$ for all $K<-\frac{1}{12}$ and all $\ell>0$ that
			\begin{equation*}
				f_{\mathbf{v}(K)+\ell,K}\left(t_m\right)>0.
			\end{equation*}
		We conclude that $f_{\mathbf{v}(K)+\ell,K}(t)$ has precisely one positive and one negative real zero for each pair $K<-\frac{1}{12}$, $\ell>0$, and thus that $t_m$ is indeed contained in the connected component of $\{f_{\mathbf{v}(K)+\ell,K}(t)>0\}$ that contains $t=0$. Hence, $\left(\begin{smallmatrix}
		1\\ t_m
		\end{smallmatrix}\right)$ is by $t_m$ being a local minimum of $f_{\mathbf{v}(K)+\ell,K}(t)$ never a hyperbolic point of the corresponding cubic polynomial $h_{\mathbf{v}(K)+\ell,K}$ and $\mathcal{H}_{\mathbf{v}(K)+\ell,K}$ is not closed.
		
		Up to this point we have shown for all
			\begin{equation*}
				\left(\begin{matrix}
				L\\ K
				\end{matrix}\right)\in\left\{K>\frac{1}{4}\right\}\cup  \left\{K\in\left[-\frac{1}{12},\frac{1}{4}\right],\ |L|>\mathbf{u}(K)\right\}\cup\left\{K<-\frac{1}{12},\ |L|\geq \mathbf{v}(K)\right\}\subset\mathbb{R}^2,
			\end{equation*}
		where $\mathbf{u}:\left[-\frac{1}{12},\frac{1}{4}\right]\to\mathbb{R}$ is the continuous extension of $\mathbf{u}$ \eqref{eqn_bfu_def}, that the maximal connected quartic GPSR curve $\mathcal{H}_{L,K}$ is not closed (see Figure \ref{fig_areaover_hardbell}) .
			\begin{figure}[H]
				\centering
				\includegraphics[scale=0.2]{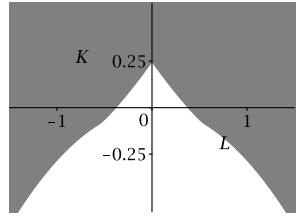}
				\caption{The set $\left\{K>\frac{1}{4}\right\}\cup \left\{K\in\left[-\frac{1}{12},\frac{1}{4}\right],\ |L|>\mathbf{u}(K)\right\}\cup\left\{K<-\frac{1}{12},\ |L|\geq \mathbf{v}(K)\right\}\subset\mathbb{R}^2$ marked in grey.}\label{fig_areaover_hardbell}
			\end{figure}
		The next step might look a bit arbitrary at first sight, but its motivation will quickly become clear. We define
			\begin{equation}\label{eqn_sharp_bdr}
				\mathbf{w}:\left(-\infty,-\frac{1}{12}\right)\to\mathbb{R},\quad \mathbf{w}(K)=\frac{\sqrt{6-216K}}{9}.
			\end{equation}
		The function $\mathbf{w}$ is positive and is smoothly extended to $K=-\frac{1}{12}$,
			\begin{equation*}
				\mathbf{w}\left(-\frac{1}{12}\right)=\frac{2\sqrt{2}}{3\sqrt{3}}=\mathbf{u}\left(-\frac{1}{12}\right)=\mathbf{v}\left(-\frac{1}{12}\right),
			\end{equation*}
		see Figure \ref{fig_u_v_w}.
			\begin{figure}[H]
				\centering
				\includegraphics[scale=0.25]{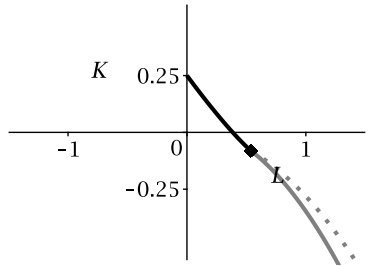}
				\caption{The graphs of $\mathbf{u}$, $\mathbf{v}$, and $\mathbf{w}$ in $\mathbb{R}^2$. The graph of $\mathbf{u}$ is the black line, the graph of $\mathbf{w}$ is the grey line, and the graph of $\mathbf{v}$ is the dotted grey line. The point $\left(L,K\right)^T=\left(\frac{2\sqrt{2}}{3\sqrt{3}},-\frac{1}{12}\right)^T$ is marked with a black diamond.}\label{fig_u_v_w}
			\end{figure}
		\noindent
		The definition of $\mathbf{w}$ is motivated by considering symbolic solutions in $L$ and $K$ of the system of equations
			\begin{equation}\label{eqn_sharpbdr_base_equations}
				\left(\Phi^*g_{\mathcal{H}_{L,K}}\right)_{t}=0,\quad \frac{\partial}{\partial t}\left(\left(\Phi^*g_{\mathcal{H}_{L,K}}\right)_{t}\right)=0,
			\end{equation}
		which can be obtained with the help of a computer algebra system like Maple. It turns out that the graph of $\mathbf{w}$ embedded in $\mathbb{R}^2$ via $K\mapsto\left(\begin{smallmatrix}
		\mathbf{w}(K)\\ K
		\end{smallmatrix}\right)$ consists of solutions in $L$ and $K$ of \eqref{eqn_sharpbdr_base_equations}. Observe that
			\begin{equation*}
				\mathbf{v}(K)=\mathbf{w}(K)\quad\Leftrightarrow\quad K=-\frac{1}{12}
			\end{equation*}
		when considering their smooth extensions to $K=-\frac{1}{12}$, and that $\mathbf{v}\left(-\frac{1}{2}\right)=\frac{4}{3}>\frac{\sqrt{114}}{9}=\mathbf{w}\left(-\frac{1}{2}\right)$, which implies
			\begin{equation}\label{eqn_vgreaterw}
				\mathbf{v}(K)>\mathbf{w}(K)
			\end{equation}
		for all $K\in\left(-\infty,-\frac{1}{12}\right)$. We will now show that for each $K<-\frac{1}{12}$ and all $L\in\left[\mathbf{w}(K),\mathbf{v}(K)\right)$, $\mathcal{H}_{L,K}$ is not closed, see Figure \ref{fig_area_between_v_w}.
			\begin{figure}[H]
				\centering
				\includegraphics[scale=0.2]{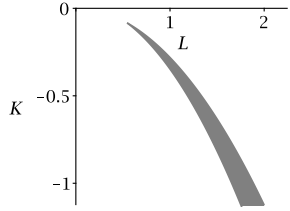}
				\caption{The area $\left.\left\{(L,K)^T\in\mathbb{R}^2\ \right|\ K<-\frac{1}{12},\ L\in[\mathbf{w}(K),\mathbf{v}(K))\right\}$ between the graphs of $\mathbf{w}$ and $\mathbf{v}$ marked in grey.}\label{fig_area_between_v_w}
			\end{figure}
		\noindent
		The inequality \eqref{eqn_vgreaterw} implies that for each such pair $(L,K)^T\in\mathbb{R}^2$ the corresponding function $f_{L,K}(t)$ has precisely one negative and one positive real zero, and furthermore only one critical value at $t=0$, since the corresponding points $t_m$ \eqref{eqn_t_m_def} and $t_M$ \eqref{eqn_t_M_def} are not real-valued in these cases. We start by showing that the graph of $\mathbf{w}$ coincides with the image of a maximal integral curve of $\mathcal{V}|_{\mathbb{R}^2\setminus\{\mathcal{V}=0\}}$. Firstly note that the graph of $\mathbf{w}$,
			\begin{equation*}
				\left\{\left(\begin{smallmatrix}\mathbf{w}(K)\\ K \end{smallmatrix}\right)\ \left|\ K\in\left(-\infty,-\frac{1}{12}\right)\right\}\right.\subset\mathbb{R}^2,
			\end{equation*}
		is contained in $\mathbb{R}^2\setminus\{\mathcal{V}=0\}$. Similarly to the consideration of the graph of $\mathbf{u}$, we note that $dK(\mathcal{V})=L\left(6K+\frac{1}{2}\right)$ does not vanish if $L\ne 0$ and $K\ne-\frac{1}{12}$, hence it does not vanish along the graph of $\mathbf{w}$. Also note that
			\begin{equation*}
				\partial \left\{\left(\begin{smallmatrix}\mathbf{w}(K)\\ K  \end{smallmatrix}\right)\ \left|\ K\in\left(-\infty,-\frac{1}{12}\right)\right\}\right.=\left\{\left(\begin{smallmatrix}
				L\\ K 
				\end{smallmatrix}
				\right)=\left(\begin{smallmatrix}
				\frac{2\sqrt{2}}{3\sqrt{3}}\\ -\frac{1}{12}
				\end{smallmatrix}
				\right)\right\}\subset\{\mathcal{V}=0\}.
			\end{equation*}
		We check that
			\begin{equation}\label{eqn_partial_w}
				\partial_K\mathbf{w}(K)=-\frac{12}{\sqrt{6-216K}}=\left.\frac{\D L(\mathcal{V})}{\D K(\mathcal{V})} \right|_{\left(\begin{smallmatrix}
				L\\ K
				\end{smallmatrix}
				\right)=\left(\begin{smallmatrix}
				\mathbf{w}(K)\\ K
				\end{smallmatrix}
				\right)},
			\end{equation}
		which shows that the graph of $\mathbf{w}$ does coincide with a maximal integral curve of $\mathcal{V}|_{\mathbb{R}^2\setminus\{\mathcal{V}=0\}}$. In order to show that each $\mathcal{H}_{\mathbf{w}(K),K}$, $K<-\frac{1}{12}$, is not closed, it thus suffices to verify the latter for one arbitrary $K<-\frac{1}{12}$. This follows from the fact that each pair of maximal connected quartic GPSR curves $\mathcal{H}_{\mathbf{w}(K),K}$, $K<-\frac{1}{12}$, is equivalent as we have shown in the beginning of this proof. We choose $K=-\frac{1}{6}$ with $\mathbf{w}\left(-\frac{1}{6}\right)=\frac{\sqrt{42}}{9}$. We now solve
			\begin{equation}\label{eqn_1_6_K_metriccond}
				\left(\Phi^*g_{\mathcal{H}_{\mathbf{w}\left(-\frac{1}{6}\right),-\frac{1}{6}}}\right)_{t}=0
			\end{equation}
		for $t$ and obtain as one solution $t_0=\frac{\sqrt{21}-3}{\sqrt{2}}$, and we have
			\begin{equation}\label{eqn_1_6_K_initial}
				f_{\mathbf{w}\left(-\frac{1}{6}\right),-\frac{1}{6}}\left(\frac{\sqrt{21}-3}{\sqrt{2}}\right)=6\sqrt{21}-27>0.
			\end{equation}
		Hence, $t_0$ is contained in the connected component of $\left\{f_{\mathbf{w}\left(-\frac{1}{6}\right),-\frac{1}{6}}(t)>0\right\}$ that contains $t=0$. But $\left(\Phi^*g_{\mathcal{H}_{\mathbf{w}\left(-\frac{1}{6}\right),-\frac{1}{6}}}\right)_{t_0}=0$ implies that $\left(\begin{smallmatrix}
		1\\ t_0
		\end{smallmatrix}
		\right)$ is not a hyperbolic point of $h_{\mathbf{w}\left(-\frac{1}{6}\right),-\frac{1}{6}}$. Hence, $\mathcal{H}_{\mathbf{w}\left(-\frac{1}{6}\right),-\frac{1}{6}}$ is not closed. We conclude that for all $K<-\frac{1}{12}$, $\mathcal{H}_{\mathbf{w}(K),K}$ is not closed. Note that it turns out that $t_0$ is the unique point in the connected component of $\left\{f_{\mathbf{w}\left(-\frac{1}{6}\right),-\frac{1}{6}}(t)>0\right\}$ containing $t=0$ that solves \eqref{eqn_sharp_bdr_g_root}, meaning that the curve $\{h=1\}$ has a flat point, cf. the related discussion in the proof of Theorem \ref{thm_incomplete_quartics} for the curve Thm. \ref{thm_incomplete_quartics} \hyperref[eqn_incomplete_qGPSRcurves_class_c]{c)}.
		
		Next, we will consider $K<-\frac{1}{12}$ and $L\in(\mathbf{w}(K),\mathbf{v}(K))$. With the help of Maple or any other modern computer algebra system we can solve the equation
			\begin{equation}\label{eqn_sharp_bdr_g_root}
				\left(\Phi^*g_{\mathcal{H}_{\mathbf{w}(K),K}}\right)_{t}=0
			\end{equation}
		for $t$ explicitly and obtain as one solution
			\begin{equation*}
				t_0=t_0(K):=\frac{\sqrt{6-216K}-3\sqrt{-2-24K}}{2}.
			\end{equation*}
		The point $t_0$ is real and positive for all $K<-\frac{1}{12}$ as the equation $t_0=0$ has no solutions for $K<-\frac{1}{12}$ and one can check that $t_0\left(-\frac{1}{6}\right)=\frac{\sqrt{21}-3}{\sqrt{2}}>0$. Thus, at $K=-\frac{1}{6}$, $t_0$ coincides with the point used in equation \eqref{eqn_1_6_K_initial}. Since $\mathbf{w}(K)<\mathbf{v}(K)$ for all $K<-\frac{1}{12}$, the corresponding functions $f_{\mathbf{w}(K),K}(t)$ have precisely one negative and one positive real zero. We calculate
			\begin{align*}
				f_{\mathbf{w}(K),K}(t_0)=162K+3888K^2+23328K^3+\left(\frac{3}{2}+45K+324K^2\right) \sqrt{6-216K}\sqrt{-2-24K}
			\end{align*}
		and obtain, again by using a computer algebra system, $f_{\mathbf{w}(K),K}(t_0)=0$ if and only if $K=-\frac{1}{12}$. Hence, we obtain with \eqref{eqn_1_6_K_initial} that 
			\begin{equation}\label{eqn_f_sign_sharp_bdr_g_root}
				f_{\mathbf{w}(K),K}(t_0)>0
			\end{equation}
		for all $K<-\frac{1}{12}$, and together with the uniqueness of the positive and negative real zeros of $f_{\mathbf{w}(K),K}(t)$ we conclude that for all $K<-\frac{1}{12}$, the point $t_0$ is contained in the connected component of $\{f_{\mathbf{w}(K),K}(t)>0\}$ that contains $t=0$. This motivates studying the expression
			\begin{equation*}
				f^2_{\mathbf{w}(K)+\ell,K}(t_0)\cdot \left(\Phi^*g_{\mathcal{H}_{\mathbf{w}(K)+\ell,K}}\right)_{t_0}
			\end{equation*}
		as a function of $K<-\frac{1}{12}$ and, depending on $K$, $\ell\in\left(0,\mathbf{v}(K)-\mathbf{w}(K)\right)$ via canonically identifying sections in $\mathrm{Sym}^2\left(\mathbb{R}^*\right)\to\mathbb{R}$ with smooth functions on $\mathbb{R}$. We obtain with the latter identification
			\begin{align*}
				f^2_{\mathbf{w}(K)+\ell,K}(t_0)\cdot  \left(\Phi^*g_{\mathcal{H}_{\mathbf{w}(K)+\ell,K}}\right)_{t_0}&=\left(-\frac{27}{8}+243K+4374K^2\right)\ell^2\\
				&+\left(\frac{27}{16}+\frac{243}{4}K\right)\sqrt{6-216K}\sqrt{-2-24K}\ell^2\\
				&+\left(81K+972K^2\right) \sqrt{6-216K}\ell\\
				&+\left(\frac{9}{2}-81K-2916K^2\right) \sqrt{-2-24K}\ell.
			\end{align*}
		Note that for all $K<-\frac{1}{12}$, $f^2_{\mathbf{w}(K)+\ell,K}(t_0)\cdot \left(\Phi^*g_{\mathcal{H}_{\mathbf{w}(K)+\ell,K}}\right)_{t_0}$ is also smooth when considered for all $\ell\in\mathbb{R}$. In order to show that $\mathcal{H}_{\mathbf{w}(K)+\ell,K}$ is not closed for any $K<-\frac{1}{12}$ and, depending on $K$, $\ell\in(0,\mathbf{v}(K)-\mathbf{w}(K))$ it thus suffices to show that $f^2_{\mathbf{w}(K)+\ell,K}(t_0)\cdot \left(\Phi^*g_{\mathcal{H}_{\mathbf{w}(K)+\ell,K}}\right)_{t_0}<0$ for these points. We will make use of
			 \begin{align*}
				\partial_\ell\left(f^2_{\mathbf{w}(K)+\ell,K}(t_0)\cdot  \left(\Phi^*g_{\mathcal{H}_{\mathbf{w}(K)+\ell,K}} \right)_{t_0}\right)&=\left(-\frac{27}{4}+486K+8748K^2\right)\ell\\
				&+\left(\frac{27}{8}+\frac{243}{2}K\right)\sqrt{6-216K}\sqrt{-2-24K}\ell\\
				&+\left(81K+972K^2\right)\sqrt{6-216K}\\
				&+\left(\frac{9}{2}-81K-2916K^2\right)\sqrt{-2-24K}.
			 \end{align*}
		Solving $\partial_\ell\left(f^2_{\mathbf{w}(K)+\ell,K}(t_0)\cdot \left(\Phi^*g_{\mathcal{H}_{\mathbf{w}(K)+\ell,K}}\right)_{t_0}\right)=0$ for $\ell$, we obtain
			\begin{align}
				&\partial_\ell\left(f^2_{\mathbf{w}(K)+\ell,K}(t_0)\cdot \left(\Phi^*g_{\mathcal{H}_{\mathbf{w}(K)+\ell,K}}\right)_{t_0}\right)=0\notag\\
				\Leftrightarrow\quad & \ell=\ell_0:=\frac{4}{3}\cdot\frac{\left(-18K-216K^2\right)\sqrt{6-216K}+\left(-1+18K+648K^2\right)\sqrt{-2-24K}}{-2+144K+2592K^2+\left(1+36K\right)\sqrt{6-216K}\sqrt{-2-24K}}.\label{eqn_ell0_def}
			\end{align}
		Furthermore, we get
			\begin{equation}
				\ell_0=0\quad\Leftrightarrow\quad K\in\left\{-\frac{1}{12},\frac{1}{36}\right\}\label{eqn_ell0_vanishing_condition}
			\end{equation}
		and
			\begin{equation*}
				\ell_0=\mathbf{v}(K)-\mathbf{w}(K)\quad\Leftrightarrow\quad K=-\frac{1}{12}.
			\end{equation*}
		This shows that $\ell_0$ is not contained in the boundary of the open interval $(0,\mathbf{v}(K)-\mathbf{w}(K))$ for any $K<-\frac{1}{12}$. We check that
			\begin{equation*}
				\ell_0|_{K=-\frac{1}{6}}=\frac{2\sqrt{2}\left(3\sqrt{21}-14\right)}{15\sqrt{21}-69}
			\end{equation*}
		and further calculate
			\begin{equation*}
				\mathbf{v}\left(-\frac{1}{6}\right)-\mathbf{w}\left(-\frac{1}{6}\right)- \ell_0|_{K=-\frac{1}{6}}=\frac{-21\sqrt{2}-92\sqrt{3}+60\sqrt{7}+5\sqrt{42}}{45\sqrt{21}-207}<0.
			\end{equation*}
		It follows that for all $K<-\frac{1}{12}$, $\ell_0\not\in\left(0,\mathbf{v}(K)-\mathbf{w}(K)\right)$. Since
			\begin{equation*}
				\partial_\ell\left(f^2_{\mathbf{w}(K)+\ell,K}(t_0)\cdot \left(\Phi^*g_{\mathcal{H}_{\mathbf{w}(K)+\ell,K}}\right)_{t_0}\right)
			\end{equation*}
		can be smoothly extended to $\ell\in\mathbb{R}$ as mentioned before, we now consider the term
			\begin{equation*}
				\left.\partial_\ell\left(f^2_{\mathbf{w}(K)+\ell,K}(t_0)\cdot \left(\Phi^*g_{\mathcal{H}_{\mathbf{w}(K)+\ell,K}}\right)_{t_0}\right)\right|_{\ell=0}.
			\end{equation*}
		We have seen in \eqref{eqn_ell0_def} and \eqref{eqn_ell0_vanishing_condition} that the sign of $\left.\partial_\ell\left(f^2_{\mathbf{w}(K)+\ell,K}(t_0)\cdot \left(\Phi^*g_{\mathcal{H}_{\mathbf{w}(K)+\ell,K}}\right)_{t_0}\right)\right|_{\ell=0}$ is constant for $K<-\frac{1}{12}$, and it thus coincides with
			\begin{equation*}
				\mathrm{sgn}\left(\left.\partial_\ell\left(f^2_{\mathbf{w}(K)+\ell,K}(t_0)\cdot \left(\Phi^*g_{\mathcal{H}_{\mathbf{w}(K)+\ell,K}}\right)_{t_0}\right)\right|_{\ell=0,\, K=-\frac{1}{6}}\right)=-1.
			\end{equation*}
		We deduce that $\partial_\ell\left(f^2_{\mathbf{w}(K)+\ell,K}(t_0)\cdot \left(\Phi^*g_{\mathcal{H}_{\mathbf{w}(K)+\ell,K}}\right)_{t_0}\right)<0$ for all $K<-\frac{1}{12}$ and correspondingly for all $\ell\in\left(0,\mathbf{v}(K)-\mathbf{w}(K)\right)$. Since by construction \eqref{eqn_sharp_bdr_g_root}
			\begin{equation*}
				\left.f^2_{\mathbf{w}(K)+\ell,K}(t_0)\cdot \left(\Phi^*g_{\mathcal{H}_{\mathbf{w}(K)+\ell,K}}\right)_{t_0}\right|_{\ell=0}= 0
			\end{equation*}
		for all $K<-\frac{1}{12}$, we conclude with \eqref{eqn_f_sign_sharp_bdr_g_root} that
			\begin{equation*}
				\left(\Phi^*g_{\mathcal{H}_{\mathbf{w}(K)+\ell,K}}\right)_{t_0}<0
			\end{equation*}
		for all $K<-\frac{1}{12}$ and correspondingly all $\ell\in(0,\mathbf{v}(K)-\mathbf{w}(K))$. This finally implies that for all such $K$ and $\ell$, $\mathcal{H}_{\mathbf{w}(K)+\ell,K}$ is not closed and thereby not a quartic CCGPSR curve.
		
		We have now, as it will turn out, identified all $(L,K)^T\in\mathbb{R}^2$, such that the corresponding maximal connected quartic GPSR curve $\mathcal{H}_{L,K}$ is not closed, namely
			\begin{equation}
				\left(\begin{matrix}
				L\\ K
				\end{matrix}\right)\in\left\{K>\frac{1}{4}\right\}\cup  \left\{K\in\left[-\frac{1}{12},\frac{1}{4}\right],\ |L|>\mathbf{u}(K)\right\}\cup\left\{K<-\frac{1}{12},\ |L|\geq \mathbf{w}(K)\right\},\label{ean_non_closed_quartic_curves_set}
			\end{equation}
		see Figure \ref{fig_final_area_nonclosed}.
			\begin{figure}[H]
				\centering
				\includegraphics[scale=0.2]{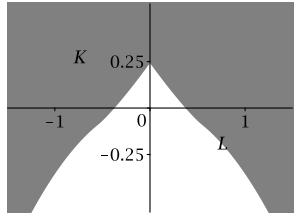}
				\caption{The set \eqref{ean_non_closed_quartic_curves_set} marked in grey.}\label{fig_final_area_nonclosed}
			\end{figure}
		\noindent
		We will now show that every point $(L,K)^T$ not contained in the set \eqref{ean_non_closed_quartic_curves_set} does indeed define a quartic CCGPSR curve $\mathcal{H}_{L,K}$, and we will show that we can choose for each equivalence class of such a curve a representative which is either \hyperref[eqn_qCCPSRcurves_class_a]{a)}, \hyperref[eqn_qCCPSRcurves_class_b]{b)}, \hyperref[eqn_qCCPSRcurves_class_c]{c)}, or contained in the one-parameter family of quartic CCGPSR curves \hyperref[eqn_qCCPSRcurves_class_d]{d)}.
		
		We start with the points contained in the image of $(\mathbf{u}(K),K)^T$, $K\in\left(-\frac{1}{12},\frac{1}{4}\right)$, cf. \eqref{eqn_image_upper_u}. We have shown that this set coincides with the image of a maximal integral curve of $\mathcal{V}|_{\mathbb{R}^2\setminus\{\mathcal{V}=0\}}$. Thus it suffices to show that one point in that set defines a quartic CCGPSR curve to conclude that all points have that property, and furthermore that the corresponding quartic CCGPSR curves are all equivalent. We choose to check
			\begin{equation*}
				\left(\begin{matrix}
				L\\ K
				\end{matrix}
				\right)=\left(\begin{matrix}
				\mathbf{u}(0)\\ 0
				\end{matrix}
				\right)=\left(\begin{matrix}
				\frac{2}{3\sqrt{3}}\\ 0
				\end{matrix}
				\right).
			\end{equation*}
		In this case, the corresponding function $f$ as defined in Lemma \ref{lem_pullback_g} coincides with $f_{\frac{2}{3\sqrt{3}},0}(t)$ and has the form $f=1-t^2+\frac{2}{3\sqrt{3}}t^3$. Equation \eqref{eqn_pullback_metric_to_dom} implies that
			\begin{equation*}
				\Phi^*g_{\mathcal{H}_{\frac{2}{3\sqrt{3}},0}}=-\frac{\partial^2f}{4f}+\frac{3\D f^2}{16f^2}.
			\end{equation*}
		Note that the connected component of $\left\{f_{\frac{2}{3\sqrt{3}},0}(t)>0\right\}$ that contains $t=0$ is given by $\left(-\frac{\sqrt{3}}{2},\sqrt{3}\right)$. The form of $f$ and \cite[Thm.\,1.1]{L2} motivate considering the cubic polynomial
			\begin{equation*}
				\widetilde{h}=x^3-xy^2+\frac{2}{3\sqrt{3}}y^3
			\end{equation*}
		with corresponding CCPSR curve in standard form $\widetilde{\mathcal{H}}\subset\{\widetilde{h}=1\}$. Its centro-affine metric is given by
			\begin{equation*}
				\left(\Phi^*g_{\widetilde{\mathcal{H}}}\right)_t=-\frac{\partial^2f}{3f}+\frac{2\D f^2}{9f^2}
			\end{equation*}
		on $\mathrm{dom} (\widetilde{\mathcal{H}} )=\left(-\frac{\sqrt{3}}{2},\sqrt{3}\right)$. Now we see that
			\begin{equation*}
				\left(\Phi^*g_{\mathcal{H}_{\frac{2}{3\sqrt{3}},0}}\right)_t=\frac{3}{4}\cdot \left(\Phi^*g_{\widetilde{\mathcal{H}}}\right)_t+\underbrace{\frac{\D f_t^2}{48 f^2(t)}}_{\geq 0}>0
			\end{equation*}
		for all $t\in\mathrm{dom} (\widetilde{\mathcal{H}} )$. Since the connected PSR curve $\left(\widetilde{\mathcal{H}},g_{\widetilde{\mathcal{H}}}\right)$ is equivalent to the curve \cite[Thm.\,8\,a)]{CHM} and is thereby a closed PSR curve, this shows that $\mathrm{dom}\left(\mathcal{H}_{\frac{2}{3\sqrt{3}},0}\right)=\mathrm{dom} (\widetilde{\mathcal{H}} )=\left(-\frac{\sqrt{3}}{2},\sqrt{3}\right)$, which coincides with the connected component of $\left\{f_{\frac{2}{3\sqrt{3}},0}(t)>0\right\}$ that contains $t=0$ and, hence, proves that $\mathcal{H}_{\frac{2}{3\sqrt{3}},0}$ is indeed also closed. Thus $\mathcal{H}_{\frac{2}{3\sqrt{3}},0}$ is in fact a quartic CCGPSR curve, which proves the claim that all points in the set described in \eqref{eqn_image_upper_u} define quartic CCGPSR curves, each equivalent to $\mathcal{H}_{\frac{2}{3\sqrt{3}},0}$. This is precisely the quartic CCGPSR curve \hyperref[eqn_qCCPSRcurves_class_c]{$c)$}. It remains to determine the connected components of $\left\{h_{\frac{2}{3\sqrt{3}},0}=x^4-x^2y^2+\frac{2}{3\sqrt{3}}xy^3=1\right\}$ and show that they are equivalent. To do so we will determine the connected components of $\left\{h_{\frac{2}{3\sqrt{3}}}>0\right\}\subset\mathbb{R}^2$. Since the quartic polynomial $h_{\frac{2}{3\sqrt{3}},0}$ is homogeneous, it suffices to study $\left\{h_{\frac{2}{3\sqrt{3}},0}\left(\left(\begin{smallmatrix}
		1\\ y
		\end{smallmatrix}
		\right)\right)>0\right\}$ and $\left\{h_{\frac{2}{3\sqrt{3}},0}\left(\left(\begin{smallmatrix}
		x\\ 1
		\end{smallmatrix}
		\right)\right)>0\right\}$. We obtain
			\begin{equation*}
				\left\{h_{\frac{2}{3\sqrt{3}},0}\left(\left(\begin{smallmatrix}
				1\\ y
				\end{smallmatrix}
				\right)\right)>0\right\}=\left(-\frac{\sqrt{3}}{2},\sqrt{3}\right)\dot\cup\left(\sqrt{3},\infty\right),
			\end{equation*}
		and, using $h_{\frac{2}{3\sqrt{3}},0}\left(\left(\begin{smallmatrix}
		x\\ 1
		\end{smallmatrix}
		\right)\right)=x\left(x+\frac{2}{\sqrt{3}}\right)\left(x-\frac{1}{\sqrt{3}}\right)^2$,
			\begin{equation*}
				\left\{h_{\frac{2}{3\sqrt{3}}}\left(\left(\begin{smallmatrix}
				x\\ 1
				\end{smallmatrix}
				\right)\right)>0\right\}=\left(-\infty,-\frac{2}{\sqrt{3}}\right)\dot\cup\left(0,\frac{1}{\sqrt{3}}\right)\dot\cup\left(\frac{1}{\sqrt{3}},\infty\right).
			\end{equation*}
		Hence,
			\begin{align}\label{eqn_h_23sqrt3_concomp}
				\left\{h_{\frac{2}{3\sqrt{3}},0}>0\right\}&=\mathbb{R}_{>0}\cdot \left\{\left(\begin{smallmatrix}
				1\\ y
				\end{smallmatrix}
				\right)\in\mathbb{R}^2\ \left|\ y\in\left(-\tfrac{\sqrt{3}}{2},\sqrt{3}\right)\right\}\right. \notag\\
				&\ \dot\cup\ \mathbb{R}_{>0}\cdot \left\{\left(\begin{smallmatrix}
				x\\ 1
				\end{smallmatrix}
				\right)\in\mathbb{R}^2\ \left|\  x\in\left(0,\tfrac{1}{\sqrt{3}}\right)\right\}\right.\notag\\
				&\ \dot\cup\ \mathbb{R}_{>0}\cdot\left\{-\left(\begin{smallmatrix}
				1\\ y
				\end{smallmatrix}
				\right)\in\mathbb{R}^2\ \left|\  y\in\left(-\tfrac{\sqrt{3}}{2},\sqrt{3}\right)\right\}\right.\notag\\
				&\ \dot\cup\ \mathbb{R}_{>0}\cdot \left\{-\left(\begin{smallmatrix}
				x\\ 1
				\end{smallmatrix}
				\right)\in\mathbb{R}^2\ \left|\  x\in\left(0,\tfrac{1}{\sqrt{3}}\right)\right\}\right..
			\end{align}
		The quartic CCGPSR curve \hyperref[eqn_qCCPSRcurves_class_c]{$c)$} is contained in the set $\mathbb{R}_{>0}\cdot\left\{\left(\begin{smallmatrix}
		1\\ y
		\end{smallmatrix}
		\right)\in\mathbb{R}^2\ \left|\ y\in\left(-\frac{\sqrt{3}}{2},\sqrt{3}\right)\right\}\right.$, and we see that it is equivalent to the unique quartic CCGPSR curve contained in the set
		$\mathbb{R}_{>0}\cdot\left\{-\left(\begin{smallmatrix}
		1\\ y
		\end{smallmatrix}
		\right)\in\mathbb{R}^2\ \left|\ y\in\left(-\frac{\sqrt{3}}{2},\sqrt{3}\right)\right\}\right.$
		via $\left(\begin{smallmatrix}
		x\\ y
		\end{smallmatrix}
		\right)\to -\left(\begin{smallmatrix}
		x\\ y
		\end{smallmatrix}
		\right).$
		We will now show that the remaining connected components of $\left\{h_{\frac{2}{3\sqrt{3}},0}>0\right\}$ also contain a (unique) quartic CCGPSR curve that is equivalent to \hyperref[eqn_qCCPSRcurves_class_c]{$c)$}. It suffices to consider the set $\mathbb{R}_{>0}\cdot \left\{\left(\begin{smallmatrix}
		x\\ 1
		\end{smallmatrix}
		\right)\in\mathbb{R}^2\ \left|\ x\in\left(0,\frac{1}{\sqrt{3}}\right)\right\}\right.$. One can easily check that the point $\left(\begin{smallmatrix}
		\frac{1}{\sqrt[4]{5}}\\ \frac{2\sqrt{3}}{\sqrt[4]{5}}\end{smallmatrix}\right)\in\left\{h_{\frac{2}{3\sqrt{3}},0}=1\right\}$ is a hyperbolic point of $h$. Consider the linear transformation of the form \eqref{eqn_A_explicit}
			\begin{equation*}
				A\left(\left(\begin{smallmatrix}
				\frac{1}{\sqrt[4]{5}}\\ \frac{2\sqrt{3}}{\sqrt[4]{5}}\end{smallmatrix}\right)\right)
				=\left(\begin{tabular}{c|c}
				$\underset{}{\vphantom{|}}\frac{1}{\sqrt[4]{5}}$ & $\frac{1}{\sqrt[4]{5}}$\\ \hline
				$\overset{}{\frac{2\sqrt{3}}{\sqrt[4]{5}}}$ & $\frac{1}{\sqrt{3}\sqrt[4]{5}}$
				\end{tabular}
				\right),
			\end{equation*}
		mapping $\left(\begin{smallmatrix}
		1\\ 0
		\end{smallmatrix}
		\right)\in\mathcal{H}_{\frac{2}{3\sqrt{3}},0}$ to said point. Then
			\begin{equation*}
				A\left(\left(\begin{smallmatrix}
				\frac{1}{\sqrt[4]{5}}\\ \frac{2\sqrt{3}}{\sqrt[4]{5}}\end{smallmatrix}\right)\right)^*h_{\frac{2}{3\sqrt{3}},0}=x^4-x^2y^2+\frac{4}{27}xy^3+\frac{4}{27}y^4.
			\end{equation*}
		We find that $\mathbf{u}\left(\frac{4}{27}\right)=\frac{4}{27}$ and deduce that the connected component of $\left\{h_{\frac{2}{3\sqrt{3}},0}=1\right\}\cap\left\{\text{hyperbolic points of }h_{\frac{2}{3\sqrt{3}},0}\right\}$ which contains the point $p=\left(\begin{smallmatrix}
		\frac{1}{\sqrt[4]{5}}\\ \frac{2\sqrt{3}}{\sqrt[4]{5}}\end{smallmatrix}\right)$ is as a GPSR curve equivalent to the quartic CCGPSR curve \hyperref[eqn_qCCPSRcurves_class_c]{$c)$}. Further note that equation \eqref{eqn_upper_bdr_system_eqns} together with Eulers theorem for homogeneous functions implies that the quartic CCGPSR curve \hyperref[eqn_qCCPSRcurves_class_c]{$c)$} is, in fact, singular at infinity. Summarising, we have shown that $\left\{h_{\frac{2}{3\sqrt{3}},0}=1\right\}$ has 4 equivalent closed connected hyperbolic components, one of which is given by the quartic CCGPSR curve $\mathcal{H}_{\frac{2}{3\sqrt{3}},0}$ which is the quartic CCGPSR curve $\mathcal{H}$ in \hyperref[eqn_qCCPSRcurves_class_c]{$c)$}. In order to determine $G^{h_{\frac{2}{3\sqrt{3}},0}}$, we need to show that the only linear map $A\in\mathrm{GL}(2)$ mapping $\mathcal{H}_{\frac{2}{3\sqrt{3}},0}$ to itself, that is $A\left(\mathcal{H}_{\frac{2}{3\sqrt{3}},0}\right)=\mathcal{H}_{\frac{2}{3\sqrt{3}},0}$, is the identity transformation $A=\mathbbm{1}\in\mathrm{GL}(2)$. Using the condition $A\cdot\left(\begin{smallmatrix}
		1\\ 0
		\end{smallmatrix}\right)\in\mathcal{H}_{\frac{2}{3\sqrt{3}},0}$, one can check that $A$ must be of the form $A=A(p)$ as in \eqref{eqn_A_explicit}, where we view $E(p)$ as an element of $\mathbb{R}\setminus\{0\}$. But then, independent of the sign $\pm$ of $E(p)$, we find with $p=\sqrt[4]{h_{\frac{2}{3\sqrt{3}},0}\left(\left(\begin{smallmatrix}1\\T\end{smallmatrix}\right)\right)}^{-1}\left(\begin{smallmatrix}
		1\\ T
		\end{smallmatrix}
		\right)$, $T\in\mathrm{dom}\left(\mathcal{H}_{\frac{2}{3\sqrt{3}},0}\right)=\left(-\frac{\sqrt{3}}{2},\sqrt{3}\right)$, and with the formulas \eqref{eqn_L_T_pmoving} and \eqref{eqn_K_T_pmoving} and the notation
			\begin{equation*}
				A(p)^*h_{\frac{2}{3\sqrt{3}},0}=x^4-x^2y^2+L(T)xy^3+K(T)y^4
			\end{equation*}
		that
			\begin{equation}\label{eqn_KT_c}
				K(T)=\frac{T\left(-T^7+4\sqrt{3}T^6+18T^5-84\sqrt{3}T^4+63T^3+432\sqrt{3}T^2-1404T+432\sqrt{3}\right)}{12\left(-T^4+2\sqrt{3}T^3-9T^2+12\sqrt{3}T-18\right)^2}.
			\end{equation}
		Using a computer algebra system like Maple, one can show that the denominator of $K(T)$ in the above formula \eqref{eqn_KT_c} does not have any zeros in $\mathrm{dom}\left(\mathcal{H}_{\frac{2}{3\sqrt{3}},0}\right)=\left(-\frac{\sqrt{3}}{2},\sqrt{3}\right)$, and the numerator has only one zero in $\mathrm{dom}\left(\mathcal{H}_{\frac{2}{3\sqrt{3}},0}\right)$, namely $T=0$. We conclude that the only linear transformations of the form \eqref{eqn_A_explicit} can either be the identity, or, corresponding to a possible minus sign of $E\left(\left(\begin{smallmatrix}
		1\\ 0
		\end{smallmatrix}
		\right)\right)$, the linear transformation $\widetilde{A}=\left(\begin{matrix}
		1 & \\
		 & -1
		\end{matrix}
		\right)$. But we can quickly check that $\widetilde{A}^*h_{\frac{2}{3\sqrt{3}},0}=x^4-x^2y^2-\frac{2}{3\sqrt{3}}xy^3\ne h_{\frac{2}{3\sqrt{3}},0}$. In order to show that
			\begin{equation*}
				G^{h_{\frac{2}{3\sqrt{3}},0}}\cong\mathbb{Z}_2\times\mathbb{Z}_2
			\end{equation*}
		and to find an explicit description in coordinates, we refer the reader to the proof of Theorem \ref{thm_incomplete_quartics} in which we show that $h_{\frac{2}{3\sqrt{3}},0}$ is equivalent to $\widehat{h}=xy(x^2-2xy+y^2)$ via \eqref{eqn_u_graph_A_trafo}. One can now quickly check that the connected components of $\{\widehat{h}=1\}$ are related by compositions of $(x,y)\to (y,x)$ and $(x,y)\to (-x,-y)$.
		
		Next we will study all maximal connected quartic GPSR curves of the form $\mathcal{H}_{0,K}$, $K\leq\frac{1}{4}$, that correspond to points of the form
			\begin{equation}\label{eqn_L0_general_case}
				\left(\begin{smallmatrix}
				L\\ K
				\end{smallmatrix}
				\right)\in\left\{K\leq \tfrac{1}{4},\ L=0\right\}.
			\end{equation}
		We will prove that each of these curves is a quartic CCGPSR curve, and furthermore that they are pairwise inequivalent.

		We have seen in \eqref{eqn_f0K_realroots} that the connected component of $\left\{f_{0,K}(t)>0\right\}$ that contains the point $t=0$ is precompact for all $K\leq \frac{1}{4}$. For each $K\leq \frac{1}{4}$, the pullback of the centro-affine metric $g_{\mathcal{H}_{0,K}}$ to $\mathrm{dom}\left(\mathcal{H}_{0,K}\right)$ fulfils by Lemma \ref{lem_pullback_g}
			\begin{align*}
				f^2_{0,K}(t)\left(\Phi^*g_{\mathcal{H}_{0,K}}\right)_t&=\left(-\frac{f_{0,K}(t)\ddot{f}_{0,K}(t)}{4}+\frac{3\dot{f}^2_{0,K}(t)}{16}\right)\D t^2\\
				&=\left(\frac{K}{2}t^4+\left(\frac{1}{4}-3K\right)t^2+\frac{1}{2}\right) \D t^2=:\widetilde{g}_K(t)\D t^2.
			\end{align*}
		First consider $K\leq 0$. Then $\ddot{f}_{0,K}(t)=-2+12Kt^2<0$ for all $t\in\mathbb{R}$. This immediately shows that $\left(\Phi^*g_{\mathcal{H}_{0,K}}\right)_t>0$ for all $t$ contained in the connected component of $\left\{f_{0,K}(t)>0\right\}$ that contains $t=0$, which thus coincides with $\mathrm{dom}\left(\mathcal{H}_{0,K}\right)$. We deduce that for all $K\leq 0$, $\mathcal{H}_{0,K}$ is closed and, hence, a quartic CCGPSR curve. Next, we will show that for all $0<K<\frac{1}{4}$ the smooth function $\widetilde{g}_K:\mathbb{R}\to\mathbb{R}$ has no real zeros and is positive. For $K=\frac{1}{4}$, we will show that $\widetilde{g}_K(t)=0$ if and only if $f_{0,K}(t)=0$. This will then imply that for all $0<K\leq\frac{1}{4}$ the set $\mathrm{dom}\left(\mathcal{H}_{0,K}\right)$ and the connected component containing $t=0$ of $\left\{f_{0,K}(t)>0\right\}$ coincide, which shows that $\mathcal{H}_{0,K}$ is closed and thus a quartic CCGPSR curve. We obtain for $0<K<\frac{1}{4}$ the symbolic equivalence
			\begin{equation*}
				\widetilde{g}_K(t)=0\quad \Leftrightarrow\quad t^2=\frac{-1+12K\pm\sqrt{144K^2-40K+1}}{4K}.
			\end{equation*}
		For $t$ to be real, one of the two possible terms $-1+12K\pm\sqrt{144K^2-40K+1}$ must be real and non-negative. Since $0<K\leq\frac{1}{4}$, this might only hold for $-1+12K+\sqrt{144K^2-40K+1}$. We find that
			\begin{equation*}
				144K^2-40K+1\geq 0\text{ and }K\in\left(0,\frac{1}{4}\right)\quad \Leftrightarrow \quad K\in\left(0,\frac{1}{36}\right],
			\end{equation*}
		and we observe that $\frac{1}{36}<\frac{1}{4}$ is a zero of $144K^2-40K+1$. This shows that $\widetilde{g}_K(t)$ might only have real zeros in the considered domain $\left(0,\frac{1}{4}\right)$ for $K$ if $K\in\left(0,\frac{1}{36}\right]$. One now verifies that $-1+12K-\sqrt{144K^2-40K+1}$ has no real zeros. Since $144K^2-40K+1\geq 0$ restricts $K$ to be an element of $\left(0,\frac{1}{36}\right]$, we evaluate 
			\begin{equation*}
				\left.-1+12K-\sqrt{144K^2-40K+1}\right|_{K=\frac{1}{72}}=-\frac{5+\sqrt{17}}{6}<0.
			\end{equation*}
		We further obtain
			\begin{equation*}
				-1+12K+\sqrt{144K^2-40K+1}=0\quad \Leftrightarrow\quad K=0
			\end{equation*}
		and deduce that the sign of the term $\frac{-1+12K+\sqrt{144K^2-40K+1}}{4K}$ is constant for $K\in\left(0,\frac{1}{36}\right]$. We further evaluate
			\begin{equation*}
				\left.-1+12K + \sqrt{144K^2-40K+1}\right|_{K=\frac{1}{72}}=-\frac{5-\sqrt{17}}{6}<0.
			\end{equation*}
		We conclude that there exists no $K\in\left(0,\frac{1}{4}\right)$, such that either term $\frac{-1+12K\mp\sqrt{144K^2-40K+1}}{4K}$ is positive. This and $\widetilde{g}_K(0)=\frac{1}{2}$ prove the claim that for all $K\in\left(0,\frac{1}{4}\right)$ the function $\widetilde{g}_K(t)$ is positive on $\mathbb{R}$ and, hence, that each $\mathcal{H}_{0,K}$ is a quartic CCGPSR curve.
		
		Now consider the case $K=\frac{1}{4}$ and note that $\left(\begin{smallmatrix}
		L\\ K
		\end{smallmatrix}\right)=\left(\begin{smallmatrix}
		0\\ \frac{1}{4}
		\end{smallmatrix}
		\right)\in\{\mathcal{V}=0\}$.
		We calculate
			\begin{equation*}
				\widetilde{g}_{\frac{1}{4}}(t)=0\quad\Leftrightarrow\quad t=\pm\sqrt{2},
			\end{equation*}
		which are precisely the zeros of $f_{0,\frac{1}{4}}(t)$. Similar to the previous case $K\in\left(0,\frac{1}{4}\right)$, $\widetilde{g}_{\frac{1}{4}}(0)=\frac{1}{2}$ implies that $\widetilde{g}_{\frac{1}{4}}$, restricted to the connected component $\left\{f_{0,\frac{1}{4}}(t)>0\right\}$ that contains the point $t=0$, is positive. We conclude that the maximal connected quartic GPSR curve $\mathcal{H}_{0,\frac{1}{4}}$ is closed and, hence, a quartic CCGPSR curve. The case $K=\frac{1}{4}$ and the cases $K<\frac{1}{4}$ correspond to the polynomial \hyperref[eqn_qCCPSRcurves_class_a]{a)} and the one-parameter family of polynomials \hyperref[eqn_qCCPSRcurves_class_d]{d)}, respectively. The quartic CCGPSR curve $\mathcal{H}_{0,\frac{1}{4}}$ is furthermore a homogeneous space under the action of the corresponding Lie group $G^{h_{0,\frac{1}{4}}}_0$. This follows from \cite[Prop.\,3.12]{L2} since $\left(0,\frac{1}{4}\right)^T\in\{\mathcal{V}=0\}$ \eqref{eqn_mathcalV_zeros}. Note that, using \cite[Prop.\,1.8]{CNS}, the homogeneity of $\mathcal{H}_{0,\frac{1}{4}}$ would also have been sufficient to prove that $\mathcal{H}_{0,\frac{1}{4}}$ is closed as a subset of $\mathbb{R}^2$, since Riemannian homogeneous spaces are always complete.
		
		It remains to prove the claim that the quartic CCGPSR curves $\mathcal{H}_{0,K}$ for $K\leq \frac{1}{4}$ are pairwise inequivalent. While proving this statement we will also determine the connected components of $\left\{h_{0,K}=1\right\}\cap\{\text{hyperbolic points of }h_{0,K}\}$ and show that they are equivalent for each fixed  $K\leq\frac{1}{4}$. It will turn out that for each $K\leq \frac{1}{4}$, every connected component of $\left\{h_{0,K}=1\right\}$ consists only of hyperbolic points.
		
		Since $\left(0,\frac{1}{4}\right)^T\in\{\mathcal{V}=0\}$ and $(0,K)^T\in\mathbb{R}^2\setminus\{\mathcal{V}=0\}$ for all $K<\frac{1}{4}$, we can use \cite[Prop.\,3.12]{L2} which implies that the connected component containing the identity of the automorphism group of $h_{0,K}$ acts transitively on $\mathcal{H}_{0,K}$ if and only if $K=\frac{1}{4}$. In particular this shows that $\mathcal{H}_{0,\frac{1}{4}}$ is not equivalent to $\mathcal{H}_{0,K}$ for any $K<\frac{1}{4}$. It remains to show that for $K<\frac{1}{4}$ the quartic CCGPSR curves $\mathcal{H}_{0,K}$ are pairwise inequivalent. For fixed $K< \frac{1}{4}$ we want to determine every $A\in \mathrm{GL}(2)$,
			\begin{equation*}
				A=\left(\begin{matrix}
				a_{11} & a_{12} \\ a_{21} & a_{22}
				\end{matrix}
				\right),
			\end{equation*}
		such that
			\begin{equation}\label{eqn_KKtilda_equiv}
				A^*h_{0,K} =x^4-x^2y^2+\widetilde{K}y^4=h_{0,\widetilde{K}}
			\end{equation}
		for some $\widetilde{K}< \frac{1}{4}$ with the additional restriction that the quartic CCGPSR curves $\mathcal{H}_{0,K}$ and $\mathcal{H}_{0,\widetilde{K}}$ are required to be equivalent via $A:\mathcal{H}_{0,\widetilde{K}}\to\mathcal{H}_{0,K}$. A necessary requirement for $A$ to fulfil \eqref{eqn_KKtilda_equiv} is that $h_{0,K}\left(\left(\begin{smallmatrix}
		a_{11}\\ a_{21}
		\end{smallmatrix}
		\right)\right)=1$ and furthermore $\left(\begin{smallmatrix}
		a_{11}\\ a_{21}
		\end{smallmatrix}
		\right)$ is required to be a hyperbolic point of $h_{0,K}$. We will treat the two cases $a_{11}\ne 0$ and $a_{11}=0$ separately. We start with assuming that $a_{11}\ne 0$. Then
			\begin{equation*}
				A^*h_{0,K} =\left.\D h_{0,K}\right|_{\left(\begin{smallmatrix}
				a_{11}\\ a_{21}
				\end{smallmatrix}
				\right)}\left(\left(\begin{smallmatrix}
				a_{12}\\ a_{22}
				\end{smallmatrix}
				\right)\right)
				x^3y + \text{remainder},
			\end{equation*}
		where the remainder in the above equation does not contain a non-zero multiple of the monomial $x^3y$. Hence, up to a scaling factor $r\in\mathbb{R}\setminus\{0\}$, 
			\begin{equation}\label{eqn_KKtilda_r_eqn}
				\left(\begin{matrix}
				a_{12}\\ a_{22}
				\end{matrix}
				\right)=r\left(\begin{matrix}
				\left.-\frac{\partial h}{\partial y}\right|_{\left(\begin{smallmatrix}
				a_{11}\\ a_{21}
				\end{smallmatrix}
				\right)}\\ \left. \frac{\partial h}{\partial x}\right|_{\left(\begin{smallmatrix}
				a_{11}\\ a_{21}
				\end{smallmatrix}
				\right)}
				\end{matrix}
				\right).
			\end{equation}
		Then, for $a_{21}=0$,
			\begin{equation*}
				A^*h_{0,K} =a_{11}^4x^4-16r^2a_{11}^8x^2y^2+256Kr^4a_{11}^12y^4.
			\end{equation*}
		In these cases we thus immediately obtain with \eqref{eqn_KKtilda_r_eqn} that $A\in\mathrm{GL}(2)$ needs to fulfil precisely
			\begin{equation*}
				A\in\left\{\mathbbm{1},-\mathbbm{1},\left(\begin{matrix}
				-1 & \\
				 & 1
				\end{matrix} 
				\right),\left(\begin{matrix}
				1 & \\ 
				 & -1
				\end{matrix} 
				\right)\right\}.
			\end{equation*}
		to solve \eqref{eqn_KKtilda_equiv}. Note that in these cases $\widetilde{K}=K$. For $a_{21}\ne 0$ consider
			\begin{align*}
				A^*h_{0,K}
				&=64r^3a_{11}a_{21}\\
				&\quad\cdot\left(-4K^3a_{21}^8+K^2\left(4a_{11}^2a_{21}^6+a_{21}^8\right)\right.\\ & \left.\quad\quad +K\left(4a_{11}^8-4a_{11}^6a_{21}^2-a_{11}^2a_{21}^6\right)-a_{11}^8+a_{11}^6a_{21}^2\right)xy^3\\
				&+\text{remainder},
			\end{align*}
		where the remainder in the above equation does not contain a non-zero multiple of the monomial $xy^3$. In order for \eqref{eqn_KKtilda_equiv} to be fulfilled, we thus need that (since by assumption $a_{11}\ne 0$ and $a_{21}\ne 0$)
			\begin{equation}\label{eqn_KKtilda_pmovingalike2}
				-4K^3a_{21}^8+K^2\left(4a_{11}^2a_{21}^6+a_{21}^8\right)+K\left(4a_{11}^8-4a_{11}^6a_{21}^2-a_{11}^2a_{21}^6\right)-a_{11}^8+a_{11}^6a_{21}^2=0.
			\end{equation}
		Solving equation \eqref{eqn_KKtilda_pmovingalike2} for $K$, we obtain
			\begin{equation*}
				\eqref{eqn_KKtilda_pmovingalike2}\quad \Leftrightarrow\quad K\in\left\{\frac{1}{4},-\frac{a_{11}^2\left(a_{11}^2-a_{21}^2\right)}{a_{21}^4},\frac{a_{11}^4}{a_{21}^4}\right\}.
			\end{equation*}
		The value $K=\frac{1}{4}$ has already been excluded. For $K=-\frac{a_{11}^2\left(a_{11}^2-a_{21}^2\right)}{a_{21}^4}$ we get
			\begin{equation*}
				A^*h_{0,K}
				\equiv 0,
			\end{equation*}
		hence we can also exclude this solution for $K$. For $K=\frac{a_{11}^4}{a_{21}^4}$, consider the condition
			\begin{equation*}
				h_{0,\frac{a_{11}^4}{a_{21}^4}}\left(\left(\begin{smallmatrix}
				a_{11}\\ a_{21}
				\end{smallmatrix}
				\right)\right)=1\quad\Leftrightarrow\quad a_{11}^2\left(2a_{11}^2-a_{21}^2\right)=1.
			\end{equation*}
		Hence, $a_{21}^2=2a_{11}^2-a_{11}^{-2}$ and, consequently,
			\begin{equation*}
				K=\frac{a_{11}^4}{a_{21}^4}=\frac{a_{11}^8}{(2a_{11}^4-1)^2}.
			\end{equation*}
		However
			\begin{equation*}
				\frac{a_{11}^8}{(2a_{11}^4-1)^2}>\frac{1}{4}
			\end{equation*}
		for all $a_{11}\ne 0$. We deduce that this solution for $K$ also does not fulfil our requirements, in this case the requirement $K<\frac{1}{4}$. Summarising, we have shown that for $a_{11}\ne 0$ the only linear transformations $A\in\mathrm{GL}(2)$ that fulfil \eqref{eqn_KKtilda_equiv} for $K<\frac{1}{4}$ are given by
			\begin{equation}\label{eqn_A_KKtilda_leq0}
				A\in\left\{\mathbbm{1},-\mathbbm{1},\left(\begin{matrix}
				-1 & \\ 
				 & 1
				\end{matrix} 
				\right),\left(\begin{matrix}
				1 & \\ 
				 & -1
				\end{matrix} 
				\right)\right\},
			\end{equation}
		and that in each case $\widetilde{K}=K$. Next consider the case $a_{11}=0$. In this case, for \eqref{eqn_KKtilda_equiv} to be true, $h_{0,K}\left(\left(\begin{smallmatrix}
		0\\ a_{21}
		\end{smallmatrix}
		\right)\right)=Ka_{21}^4$ must fulfil $Ka_{21}^4=1$. We see that this can only be the case for positive $K$, and thus can already say that for all $K\leq 0$ the transformation $A$ solves \eqref{eqn_KKtilda_equiv} if and only if $A\in\left\{\mathbbm{1},-\mathbbm{1},\left(\begin{matrix}
		-1 & \\ 
		 & 1
		\end{matrix} 
		\right),\left(\begin{matrix}
		1 & \\ 
		 & -1
		\end{matrix} 
		\right)\right\}$ and that in these cases $\widetilde{K}=K$. For $0<K<\frac{1}{4}$ we obtain
			\begin{equation*}
				Ka_{21}^4=1\quad\Leftrightarrow\quad a_{21}=\pm\frac{1}{\sqrt[4]{K}}.
			\end{equation*}
		Under the assumption that $a_{11}=0$ and $a_{21}=\pm\frac{1}{\sqrt[4]{K}}$ we find that
			\begin{equation*}
				A^*h_{0,K} =4Ka_{21}^3a_{22}x^3y+\text{remainder}= \pm 4\sqrt[4]{K}a_{22}x^3y+\text{remainder},
			\end{equation*}
		where the above remainder does not contain any $x^3y$-term. Since the $x^3y$-part of $(A^*h_{0,K})\left(\left(\begin{smallmatrix}
		x\\ y
		\end{smallmatrix}
		\right)\right)$ must vanish for \eqref{eqn_KKtilda_equiv} to be fulfilled, we deduce that $a_{22}=0$. So $A$ must be of the form $A=\left(\begin{matrix}
		0 & a_{12}\\
		a_{21} & 0
		\end{matrix}
		\right)$. Now for the final step we check that
			\begin{equation*}
				A^*h_{0,K} =-a_{21}^2a_{12}^2x^2y^2+\text{remainder}= -\frac{a_{21}^2}{\sqrt{K}}x^2y^2+\text{remainder},
			\end{equation*}
		the above remainder not containing any $x^2y^2$-term, which shows that $a_{21}=\pm\sqrt[4]{K}$. Summarising, we have shown that the possible $A\in\mathrm{GL}(2)$ that fulfil \eqref{eqn_KKtilda_equiv} for $0<K<\frac{1}{4}$ can only be of the form
			\begin{align}\label{eqn_A_KKtilda_geq0}
				A&\in\left\{\mathbbm{1},-\mathbbm{1},\left(\begin{matrix}
				-1 & \\ 
				 & 1
				\end{matrix} 
				\right),\left(\begin{matrix}
				1 & \\ 
				 & -1
				\end{matrix} 
				\right),\right.\notag\\
				&\quad\left.\left(\begin{matrix}
				  & \sqrt[4]{K} \\ 
				 \frac{1}{\sqrt[4]{K}} & 
				\end{matrix} 
				\right),\left(\begin{matrix}
				  & -\sqrt[4]{K} \\ 
				 \frac{1}{\sqrt[4]{K}} & 
				\end{matrix} 
				\right),\left(\begin{matrix}
				  & \sqrt[4]{K} \\ 
				 -\frac{1}{\sqrt[4]{K}} & 
				\end{matrix} 
				\right),\left(\begin{matrix}
				  & -\sqrt[4]{K} \\ 
				 -\frac{1}{\sqrt[4]{K}} & 
				\end{matrix} 
				\right)\right\},
			\end{align}
		and one can easily check that these matrices actually solve \eqref{eqn_KKtilda_equiv} with $\widetilde{K}=K$. We now conclude that the quartic CCGPSR curves $\mathcal{H}_{0,K}$ for $K\leq\frac{1}{4}$ are pairwise inequivalent. Now, similar to the quartic CCGPSR curve \hyperref[eqn_qCCPSRcurves_class_c]{$c)$}, equation \eqref{eqn_h_23sqrt3_concomp}, we will determine for each $K\leq \frac{1}{4}$ the connected components of $\left\{h_{0,K}>0\right\}\subset\mathbb{R}^2$. We have
			\begin{align*}
				h_{0,K}\left(\left(\begin{smallmatrix}
				1\\y
				\end{smallmatrix}\right)\right)=1-y^2+Ky^4,\quad  h_{0,K}\left(\left(\begin{smallmatrix}
				x\\1
				\end{smallmatrix}\right)\right)=x^4-x^2+K,
			\end{align*}
		and obtain
			\begin{align*}
				\underline{\text{for }K=\tfrac{1}{4}:}\left\{h_{0,\frac{1}{4}}>0\right\}&=\mathbb{R}_{>0}\cdot \left.\left\{\left(\begin{smallmatrix}
				1\\ y
				\end{smallmatrix}
				\right)\in\mathbb{R}^2\ \right|\ y\in\left(-\sqrt{2},\sqrt{2}\right)\right\}\notag\\
				&\ \dot\cup\ \mathbb{R}_{>0}\cdot \left\{\left(\begin{smallmatrix}
				x\\ 1
				\end{smallmatrix}
				\right)\in\mathbb{R}^2\ \left|\  x\in\left(-\tfrac{1}{\sqrt{2}},\tfrac{1}{\sqrt{2}}\right)\right\}\right.\notag\\
				&\ \dot\cup\ \mathbb{R}_{>0}\cdot\left.\left\{-\left(\begin{smallmatrix}
				1\\ y
				\end{smallmatrix}
				\right)\in\mathbb{R}^2\ \right|\ y\in\left(-\sqrt{2},\sqrt{2}\right) \right\} \notag\\
				&\ \dot\cup\ \mathbb{R}_{>0}\cdot \left\{-\left(\begin{smallmatrix}
				x\\ 1
				\end{smallmatrix}
				\right)\in\mathbb{R}^2\ \left|\  x\in\left(-\tfrac{1}{\sqrt{2}},\tfrac{1}{\sqrt{2}}\right)\right\}\right.,
			\end{align*}
		where we observe that for $K=\frac{1}{4}$, $h_{0,\frac{1}{4}}\leq 0$ if and only if $h_{0,\frac{1}{4}}=0$,
			\begin{align*}
				\underline{\text{for }K\in\left(0,\tfrac{1}{4}\right):}&\\
				\left\{h_{0,K}>0\right\}&=\mathbb{R}_{>0}\cdot\left\{\left(\begin{smallmatrix}
				1\\ y
				\end{smallmatrix}
				\right)\in\mathbb{R}^2\ \left|\  y\in\left(-\sqrt{\tfrac{2}{K}}\sqrt{1-\sqrt{1-\tfrac{K}{4}}},\sqrt{\tfrac{2}{K}}\sqrt{1-\sqrt{1-\tfrac{K}{4}}}\right)\right\}\right.\notag\\
				&\ \dot\cup\ \mathbb{R}_{>0}\cdot \left\{\left(\begin{smallmatrix}
				x\\ 1
				\end{smallmatrix}
				\right)\in\mathbb{R}^2\ \left|\  x\in\left(-\sqrt{\tfrac{1-\sqrt{1-4K}}{2}},\sqrt{\tfrac{1-\sqrt{1-4K}}{2}}\right)\right\}\right.\notag\\
				&\ \dot\cup\ \mathbb{R}_{>0}\cdot\left\{-\left(\begin{smallmatrix}
				1\\ y
				\end{smallmatrix}
				\right)\in\mathbb{R}^2\ \left|\  y\in\left(-\sqrt{\tfrac{2}{K}}\sqrt{1-\sqrt{1-\tfrac{K}{4}}},\sqrt{\tfrac{2}{K}}\sqrt{1-\sqrt{1-\tfrac{K}{4}}}\right)\right\}\right.\notag\\
				&\ \dot\cup\ \mathbb{R}_{>0}\cdot \left\{-\left(\begin{smallmatrix}
				x\\ 1
				\end{smallmatrix}
				\right)\in\mathbb{R}^2\ \left|\  x\in\left(-\sqrt{\tfrac{1-\sqrt{1-4K}}{2}},\sqrt{\tfrac{1-\sqrt{1-4K}}{2}}\right)\right\}\right.,
			\end{align*}
			\begin{align*}
				\underline{\text{for }K=0:}\quad  \left\{h_{0,0}>0\right\}&=\mathbb{R}_{>0}\cdot\left.\left\{\left(\begin{smallmatrix}
				1\\ y
				\end{smallmatrix}
				\right)\in\mathbb{R}^2\ \right|\ y\in\left(-1,1\right)\right\}\notag\\
				&\ \dot\cup\ \mathbb{R}_{>0}\cdot\left.\left\{-\left(\begin{smallmatrix}
				1\\ y
				\end{smallmatrix}
				\right)\in\mathbb{R}^2\ \right|\ y\in\left(-1,1\right)\right\},
			\end{align*}
		and
			\begin{align*}
				\underline{\text{for }K<0:}&\notag\\ 
				\left\{h_{0,K}>0\right\}
				&=\mathbb{R}_{>0}\cdot\left\{\left(\begin{smallmatrix}
				1\\ y
				\end{smallmatrix}
				\right)\in\mathbb{R}^2\ \left|\  y\in\left(-\sqrt{\tfrac{2}{K}}\sqrt{1+\sqrt{1-\tfrac{K}{4}}},\sqrt{\tfrac{2}{K}}\sqrt{1+\sqrt{1-\tfrac{K}{4}}}\right)\right\}\right.\notag\\
				&\ \dot\cup\ \mathbb{R}_{>0}\cdot \left\{-\left(\begin{smallmatrix}
				1\\ y
				\end{smallmatrix}
				\right)\in\mathbb{R}^2\ \left|\  y\in\left(-\sqrt{\tfrac{2}{K}}\sqrt{1+\sqrt{1-\tfrac{K}{4}}},\sqrt{\tfrac{2}{K}}\sqrt{1+\sqrt{1-\tfrac{K}{4}}}\right)\right\}\right..
			\end{align*}
		For $K\leq0$ we see that $\left\{h_{0,K}>0\right\}$ has exactly $2$ connected components, and the corresponding unique contained quartic CCGPSR curves are equivalent via $A=-\mathbbm{1}$, cf. \eqref{eqn_A_KKtilda_leq0}. For $0<K<\frac{1}{4}$, $\left\{h_{0,K}>0\right\}$ has exactly $4$ connected components, and the corresponding quartic CCGPSR curves are equivalent via compositions of $A=\left(\begin{matrix}
		1 & \\ & -1
		\end{matrix}\right)$ and $A=\left(\begin{matrix}
		 & -\sqrt[4]{K}\\
		 \frac{1}{\sqrt[4]{K}} & 
		\end{matrix}
		\right)$, cf. \eqref{eqn_A_KKtilda_geq0}. In the case $K=\frac{1}{4}$, $\left\{h_{0,\frac{1}{4}}>0\right\}$ has exactly $4$ connected components, and one can check that transformations $A\in\mathrm{GL}(2)$ of the form \eqref{eqn_A_KKtilda_geq0} are also automorphisms of $h_{0,\frac{1}{4}}$ as for $0<K<\frac{1}{4}$ and map the corresponding quartic CCGPSR curves bijectively to each other. Now, for the automorphism groups $G^{h_{0,K}}$ of $h_{0,K}$ for $K\leq \frac{1}{4}$, \eqref{eqn_A_KKtilda_leq0} and \eqref{eqn_A_KKtilda_geq0} imply that
			\begin{equation*}
				G^{h_{0,K}}\cong \mathbb{Z}_2\times\mathbb{Z}_2
			\end{equation*}
		for all $K\leq 0$, and
			\begin{equation*}
				G^{h_{0,K}}\cong \mathbb{Z}_2\ltimes\mathbb{Z}_2.
			\end{equation*}
		for all $K\in\left(0,\frac{1}{4}\right)$. For a geometric interpretation of the action of $G^{h_{0,\frac{1}{4}}}_0$, observe that $h_{0,\frac{1}{4}}=\left(x^2-\tfrac{1}{2}y^2\right)^2$. Hence, $G^{h_{0,\frac{1}{4}}}_0$ acts precisely via hyperbolic rotations of the Lorentz vector space $(\mathbb{R}^2,\D x^2-\tfrac{1}{2}\D y^2)$ which are given by
			\begin{align}
				\left(\begin{matrix}
				\cosh(t) & \frac{\sinh(t)}{\sqrt{2}}\\
				\sqrt{2}\sinh(t) & \cosh(t)
				\end{matrix}
				\right),\quad t\in\mathbb{R}.\label{eqn_G014_identification}
			\end{align}
		Summarising, $\mathcal{H}_{0,\frac{1}{4}}$ is a homogeneous space and
			\begin{equation*}
				G^{h_{0,\frac{1}{4}}}\cong \mathbb{R}\ltimes(\mathbb{Z}_2\ltimes\mathbb{Z}_2).%\label{eqn_G014_identification}
			\end{equation*}
		Lastly, we need to show that for all $K<\frac{1}{4}$, $\mathcal{H}_{0,K}$ is not singular at infinity. This is equivalent to showing that $\dot f_{0,K}(t)$ does not vanish at the smallest positive and biggest negative zero of $f_{0,K}(t)$ and we leave the calculations as a small exercise for the reader.
		
		Now we will consider points
			\begin{equation*}
				\left(\begin{matrix}
				L\\ K
				\end{matrix}
				\right)\in\left\{K<-\frac{1}{12},\ |L|<\mathbf{w}(K)\right\}\subset\mathbb{R}^2
			\end{equation*}
		and the corresponding maximal connected quartic GPSR curves $\mathcal{H}_{L,K}$. We will proceed as follows. We will show that for any such $(L,K)^T\in\left\{K<-\frac{1}{12},\ |L|<\mathbf{w}(K)\right\}$, the image of the maximal integral curve $\gamma$ of $\mathcal{V}|_{\mathbb{R}^2\setminus\{\mathcal{V}=0\}}$ with $\gamma(0)=(L,K)^T$ contains a point of the form $\left(0,\widetilde{K}\right)^T\in\mathbb{R}^2$, $\widetilde{K}<-\frac{1}{12}$. This will thus imply that $\mathcal{H}_{L,K}$ is equivalent to $\mathcal{H}_{0,\widetilde{K}}$ and, since we have already seen that $\mathcal{H}_{0,\widetilde{K}}$ is a quartic CCGPSR curve for all $\widetilde{K}\leq\frac{1}{4}$, that $\mathcal{H}_{L,K}$ is a quartic CCGPSR curve. We will without loss of generality assume that $0<L<\mathbf{w}(K)$ for $K<-\frac{1}{12}$ fixed, as we have already dealt with the case $L=0$ above \eqref{eqn_L0_general_case} and since $\mathcal{H}_{L,K}$ and $\mathcal{H}_{-L,K}$ are equivalent for all $L$. Instead of checking the maximal integral curves of $\mathcal{V}|_{\mathbb{R}^2\setminus\{\mathcal{V}=0\}}$, respectively their restriction to the set $\left\{K<-\frac{1}{12},\ 0<L<\mathbf{w}(K)\right\}$,
		directly, we will first transform this set using a suitable diffeomorphism. Recall that $\mathbf{w}(K)=\frac{\sqrt{6-216K}}{9}>0$ for all $K<-\frac{1}{12}$ and consider the smooth map
			\begin{align}
				&F:\left\{K<-\frac{1}{12},\ 0<L<\mathbf{w}(K)\right\}\to \left\{\psi<-\frac{1}{12},\ 0<\varphi<1\right\},\notag\\
				&F:\left(\begin{matrix}
				L\\ K
				\end{matrix}
				\right)\mapsto \left(\begin{matrix}
				\frac{L}{\mathbf{w}(K)}\\ K
				\end{matrix}
				\right),\label{eqn_first_F_trafo}
			\end{align}
		where $\left(\varphi,\psi\right)$ denote the canonical coordinates of $\left\{\psi<-\frac{1}{12},\ 0<\varphi<1\right\}\subset\mathbb{R}^2$ (see Figure \ref{fig_F_domain_codomain}).
			\begin{figure}[H]
				\centering
				\begin{subfigure}[h]{0.3\linewidth}
					\includegraphics[scale=0.2]{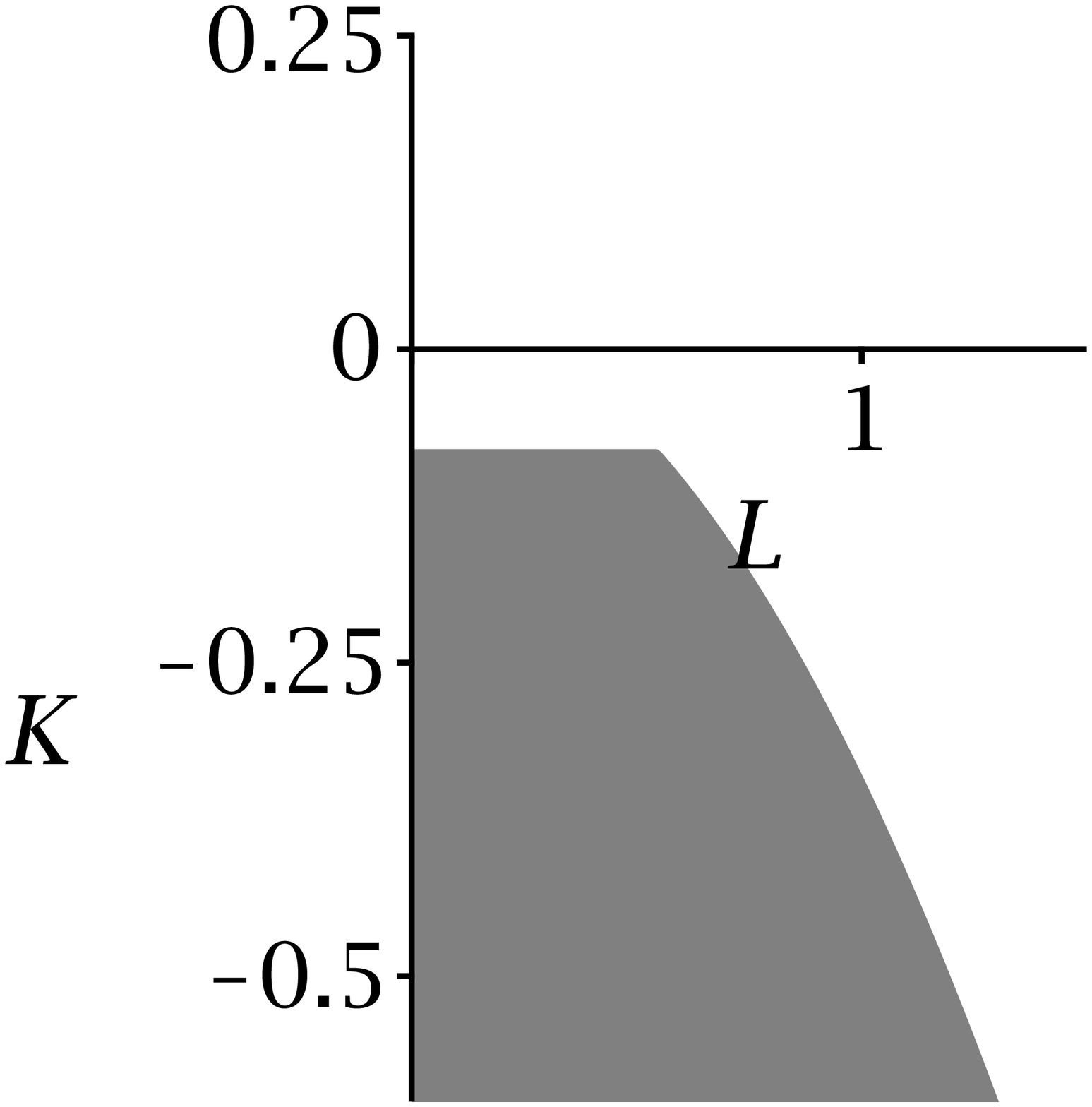}
				\end{subfigure}
				\begin{subfigure}[h]{0.3\linewidth}
					\includegraphics[scale=0.2]{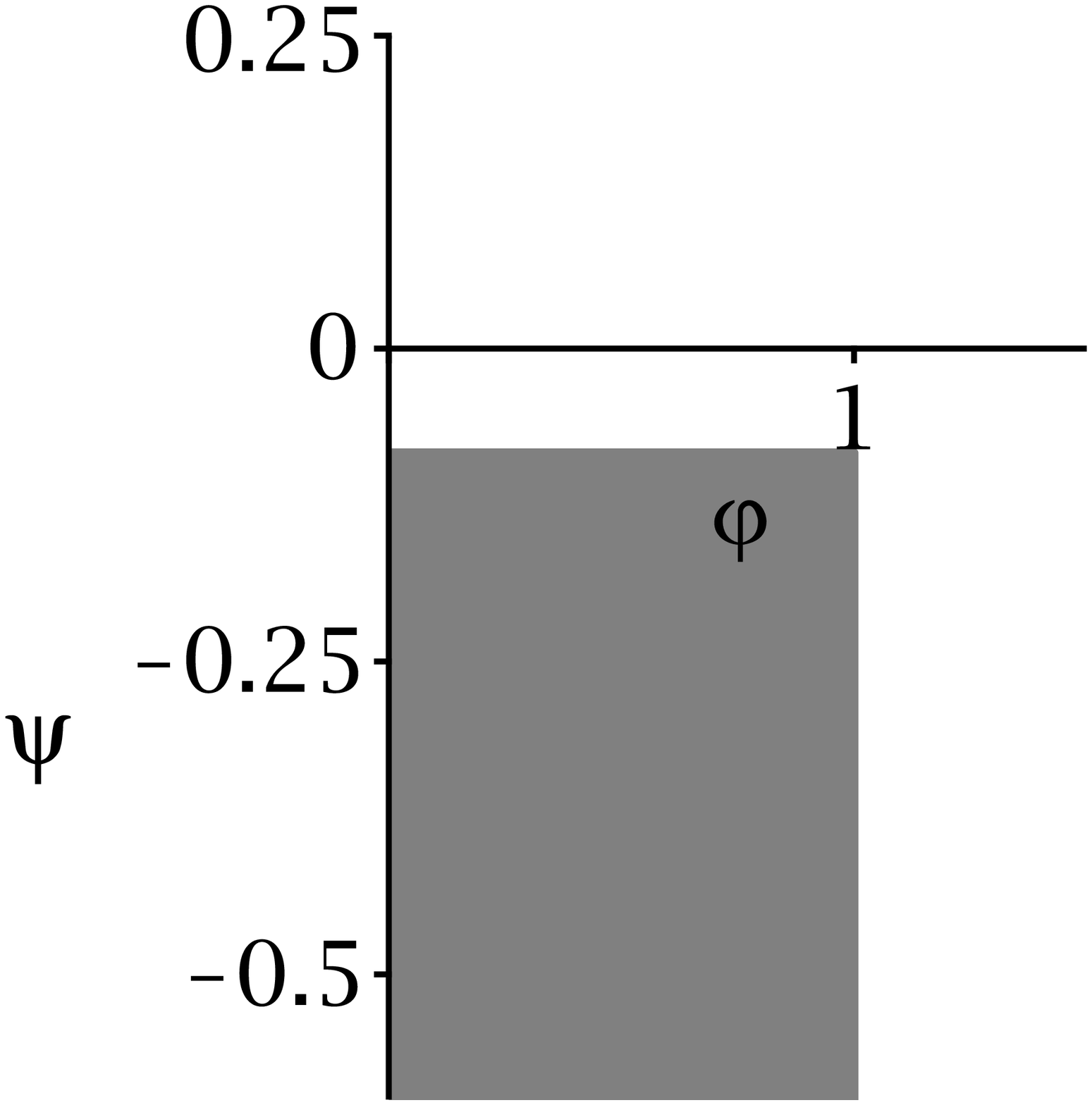}
				\end{subfigure}
				\caption{Domain and codomain of $F$ \eqref{eqn_first_F_trafo}.}\label{fig_F_domain_codomain}
			\end{figure}
		\noindent
		The differential of $F$ is given by
			\begin{equation*}
				\D F=\left(\begin{matrix}
				\frac{1}{\mathbf{w}(K)} & \frac{162L}{(1-36K)\sqrt{6-216K}}\\
				0 & 1
				\end{matrix}
				\right)
			\end{equation*}
		and we see that $F$ is, in fact, a diffeomorphism. We obtain for the inverse of $F$
			\begin{equation*}
				F^{-1}\left(\left(\begin{smallmatrix}
				\varphi\\ \psi
				\end{smallmatrix}
				\right)\right)=\left(\begin{smallmatrix}
				\varphi\mathbf{w}(\psi)\\ \psi
				\end{smallmatrix}
				\right)
			\end{equation*}
		and for the push-forward of the vector field $\mathcal{V}$ restricted to ${\{K<-1/12,\, 0<L<\mathbf{w}(K)\}}$
			\begin{equation}
				F_*\mathcal{V}=\left(F_*\mathcal{V}\right)_{\left(\begin{smallmatrix}
				\varphi\\ \psi
				\end{smallmatrix}
				\right)}= \frac{9\left(1-\varphi^2\right)(-1+4\psi)}{\sqrt{6-216\psi}}  \partial_{\varphi}+ \frac{\varphi(1+12\psi)\sqrt{6-216\psi}}{18}\partial_{\psi}.\label{eqn_first_F_V_pushforward}
			\end{equation}
		Since $\D\psi\left(F_*\mathcal{V}\right)<0$ for all $\left(\begin{smallmatrix}
		\varphi\\ \psi
		\end{smallmatrix}
		\right)\in\left\{\psi<-\frac{1}{12},\ 0<\varphi<1\right\}$, we can define the smooth function
			\begin{align*}
				&\mathcal{R}:\left\{\psi<-\frac{1}{12},\ 0<\varphi<1\right\}\to\mathbb{R},\\
				&\mathcal{R}\left(\left(\begin{smallmatrix}
				\varphi\\ \psi
				\end{smallmatrix}
				\right)\right):=\frac{d\varphi\left(F_*\mathcal{V}\right)}{d\psi\left(F_*\mathcal{V}\right)}=\frac{1-\varphi^2}{\varphi}\cdot\frac{27(1-4\psi)}{(1+12\psi)(-1+36\psi)}.
			\end{align*}
		Instead of studying the images of the maximal integral curves of the vector field
			\begin{equation*}
				F_*\mathcal{V}|_{\left\{\psi<-\frac{1}{12},\ 0<\varphi<1\right\}},
			\end{equation*}
		we can now study the images of the maximal integral curves of the vector field
			\begin{equation}\label{eqn_mathcalX_def}
				\mathcal{X}:=\mathcal{R}\,\partial_\varphi+\partial_\psi.
			\end{equation}
		defined on the set $\left\{\psi<-\frac{1}{12},\ 0<\varphi<1\right\}=F\left(\left\{K<-\frac{1}{12},\ 0<L<\mathbf{w}(K)\right\}\right)\subset\mathbb{R}^2$, since there is a one-to-one correspondence between them, which follows from $d\psi\left(F_*\mathcal{V}\right)\ne 0$ on said set. It turns out that we can, in fact, find the general solution of the equation for integral curves of $\mathcal{X}$. For $t<-\frac{1}{12}$ and some $a<-\frac{1}{12}$, consider with $\gamma:\left(a,-\frac{1}{12}\right)\to\left\{\psi<-\frac{1}{12},\ 0<\varphi<1\right\}$, $\gamma(t)=\left(\begin{smallmatrix}
		\varphi(t)\\ t
		\end{smallmatrix}
		\right)$,
			\begin{equation*}
				\mathcal{X}_\gamma=\dot{\gamma}\quad
				\Leftrightarrow\quad \frac{\varphi\dot{\varphi}}{1-\varphi^2}=\frac{27(1-4t)}{(-1+36t)(1+12t)}\quad
				\Leftrightarrow\quad \varphi(t)=\sqrt{1-\frac{c\sqrt{-(1+12t)^3}}{1-36t}},
			\end{equation*}
		where $c\in\mathbb{R}_{>0}$ is chosen in such a way so that the initial condition
			\begin{equation*}
				\gamma(t_0)=\left(\begin{smallmatrix}
				\varphi(t_0)\\ t_0
				\end{smallmatrix}\right)\in\left\{\psi<-\frac{1}{12},\ 0<\varphi<1\right\}
			\end{equation*}
		is fulfilled, see Figure \ref{fig_gamma_example} for an example of such a curve $\gamma$ (note: in our construction, the initial time $t_0$ fulfils $t_0<-\frac{1}{12}$).
			\begin{figure}[H]
				\centering
				\includegraphics[scale=0.2]{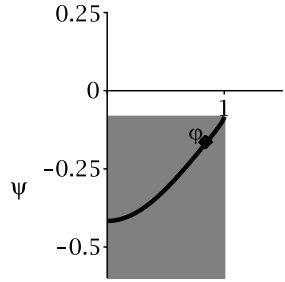}
				\caption{Example of a curve $\gamma$ (in black) fulfilling $\mathcal{X}_\gamma=\dot{\gamma}$ with initial condition $\varphi\left(-\frac{1}{6}\right)=\frac{\sqrt{6}}{\sqrt{7}}$ marked with a diamond.}\label{fig_gamma_example}
			\end{figure}
		\noindent
		Observe that for all $t<-\frac{1}{12}$ in the domain of $\varphi(t)$ for any fixed $c$, we have $\varphi(t)<1$. This holds in particular arbitrarily close to $t=-\frac{1}{12}$. We will now show that $\varphi(t)$ cannot converge to the value $1$ in finite negative time. Solving $\varphi(t)=1$, we obtain as the unique negative solution $t=-\frac{1}{12}$, but this is the upper bound of the domain of definition of $\gamma(t)$ and, hence, $\varphi(t)$. This shows that for all $t<-\frac{1}{12}$ for which $\varphi(t)$ is defined we indeed have $\varphi(t)<1$. Thus, if we can prove that each such curve $\gamma$ independent of the initial condition $(t_0,\gamma(t_0))$ converges to a point in the set $\left\{\psi<-\frac{1}{12},\ \varphi=0\right\}$ in finite negative time $t<-\frac{1}{12}$, then we will have shown using the fact $\mathcal{V}|_{\{K<-1/12,\, |L|<\mathbf{w}(K)\}}\ne 0$ that each maximal integral curve of $\mathcal{V}|_{\{K<-1/12,\, |L|<\mathbf{w}(K)\}}$ meets the set $\left\{K<-\frac{1}{12},\ L=0\right\}$ in either finite positive or finite negative time.
		Note that all possible integral curves $\gamma$ that fulfil $\varphi(t_0)>0$ move toward $\left\{\psi=-\frac{1}{12}\right\}$ in positive time-direction.
		This means that we have to solve
			\begin{align*}
				&\sqrt{1-\frac{c\sqrt{-(1+12t)^3}}{ 1-36t }}=0,\ t<-\frac{1}{12}\\
				\Leftrightarrow\quad & 1728c^2t^3+\left(432c^2+1296\right)+\left(36c^2-72\right)t^2+1=0,\ t<-\frac{1}{12}
			\end{align*}
		with the restriction that we are interested in the biggest possible negative solution in $t$. Substituting $t=s-\frac{1}{12}$ and dividing by $16$, we observe that
			\begin{align}
				&1728c^2t^3+\left(432c^2+1296\right)+\left(36c^2-72\right)t^2+1=0,\ t<-\frac{1}{12}\notag\\
				\Leftrightarrow\quad & 108 c^2s^3+81s^2-18s+1=0,\ s<0.\label{eqn_lower_area_s_root}
			\end{align}
		We see that equation \eqref{eqn_lower_area_s_root} always has a negative, and in particular uniquely determined biggest negative, solution in $s$ since it takes the value $1$ for $s=0$. Note that this solution coincides with the minimal possible value $a$ that was used to denote the domain of definition $\left(a,-\frac{1}{12}\right)\subset\mathbb{R}$ of $\gamma(t)$, respectively $\varphi(t)$. We deduce that each maximal integral curve of $\mathcal{X}$ does meet the set $\left\{\psi<-\frac{1}{12},\ \varphi=0\right\}$ in finite negative time and, hence, that each maximal integral curve $\gamma$ of $\mathcal{V}|_{\mathbb{R}^2\setminus\{\mathcal{V}=0\}}$ with initial condition $\gamma(0)\in\left\{K<-\frac{1}{12},\ |L|<\mathbf{w}(K)\right\}$ meets the set $\left\{K<-\frac{1}{12},\ L=0\right\}$ in either finite negative or finite positive time. Hence, each maximal connected quartic GPSR curve $\mathcal{H}_{L,K}$, $(L,K)^T\in\left\{K<-\frac{1}{12},\ |L|<\mathbf{w}(K)\right\}$, is equivalent to a maximally extended quartic GPSR curve of the form $\mathcal{H}_{0,\widetilde{K}}$ with $\widetilde{K}<-\frac{1}{12}$. We have already shown that the quartic CCGPSR curves $\mathcal{H}_{0,K}$, $K\leq\frac{1}{4}$, are pairwise inequivalent. Hence, the value for $\widetilde{K}$ in dependence of the initial $(L,K)$ is unique. We deduce that $\mathcal{H}_{L,K}$ is closed in $\mathbb{R}^2$ and thus a quartic CCGPSR curve as claimed.
		
		Now consider the set $\left\{K=-\frac{1}{12},\ |L|<\mathbf{w}\left(-\frac{1}{12}\right)=\frac{2\sqrt{2}}{3\sqrt{3}}\right\}$ and the restriction of $\mathcal{V}$ to it. It turns out that this set coincides with the image of a maximal integral curve of $\mathcal{V}|_{\mathbb{R}^2\setminus\{\mathcal{V}=0\}}$. This follows from the fact that $\left\{K=-\frac{1}{12},\ |L|<\frac{2\sqrt{2}}{3\sqrt{3}}\right\}\subset\mathbb{R}^2\setminus\{\mathcal{V}=0\}$, and that $\mathcal{V}$ is parallel to $\left\{K=-\frac{1}{12},\ |L|<\frac{2\sqrt{2}}{3\sqrt{3}}\right\}$ in the sense that 
			\begin{equation*}
				\D K\left(\mathcal{V}|_{\left\{K=-\frac{1}{12},\ |L|<\frac{2\sqrt{2}}{3\sqrt{3}}\right\}}\right)\equiv 0.
			\end{equation*}
		Hence, every maximal connected quartic GPSR curve of the form $\mathcal{H}_{L,-\frac{1}{12}}$ with $|L|<\frac{2\sqrt{2}}{3\sqrt{3}}$ is equivalent to $\mathcal{H}_{0,-\frac{1}{12}}$ and thus closed.
			
		Lastly, we have to consider the maximal connected quartic GPSR curves $\mathcal{H}_{L,K}$ with
			\begin{equation}
				\left(\begin{matrix}
				L\\ K
				\end{matrix}
				\right)\in\left\{-\frac{1}{12}<K<\frac{1}{4},\  |L|<\mathbf{u}(K)\right\}\subset\mathbb{R}^2,\label{eqn_set_upper_bell}
			\end{equation}
		respectively the restriction of $\mathcal{V}$ to the above set. Recall that
			\begin{equation*}
				\mathbf{u}(K)=\frac{\sqrt{2}}{3\sqrt{3}}\sqrt{1-36K+\sqrt{(1+12K)^3}},
			\end{equation*}
		cf. \eqref{eqn_upper_bdr}. We proceed similarly to the case where we considered points of the form $\left(\begin{smallmatrix}
		L\\ K
		\end{smallmatrix}
		\right)\in\left\{K<-\frac{1}{12},\ |L|<\mathbf{w}(K)\right\}$. We will show that every maximal integral curve of the restricted vector field $\mathcal{V}|_{\left\{-\frac{1}{12}<K<\frac{1}{4},\ |L|<\mathbf{u}(K)\right\}}$ contains a point of the form
			\begin{equation*}
				\left(\begin{matrix}
				0\\ \widetilde{K}
				\end{matrix}\right)\in\left\{-\frac{1}{12}<K<\frac{1}{4},\  |L|<\mathbf{u}(K)\right\}.
			\end{equation*}
		To do so, it suffices to consider points in $\left\{-\frac{1}{12}<K<\frac{1}{4},\ 0<L<\mathbf{u}(K)\right\}$ which provide the initial value for said integral curves. We define
			\begin{align}
				&\widetilde{F}:\left\{-\frac{1}{12}<K<\frac{1}{4},\  0<L<\mathbf{u}(K)\right\}\to\left\{-\frac{1}{12}<\psi<\frac{1}{4},\ 0<\varphi<1\right\},\notag\\
				&\widetilde{F}:\left(\begin{matrix}
				L\\ K
				\end{matrix}
				\right)\mapsto \left(\begin{matrix}
				\frac{L}{\mathbf{u}(K)}\\ K
				\end{matrix}
				\right),\label{eqn_first_F_tilda_trafo}
			\end{align}
		see Figure \ref{fig_Ftilda_domain_codomain}.
				\begin{figure}[H]
					\centering
					\begin{subfigure}[h]{0.3\linewidth}
						\includegraphics[scale=0.2]{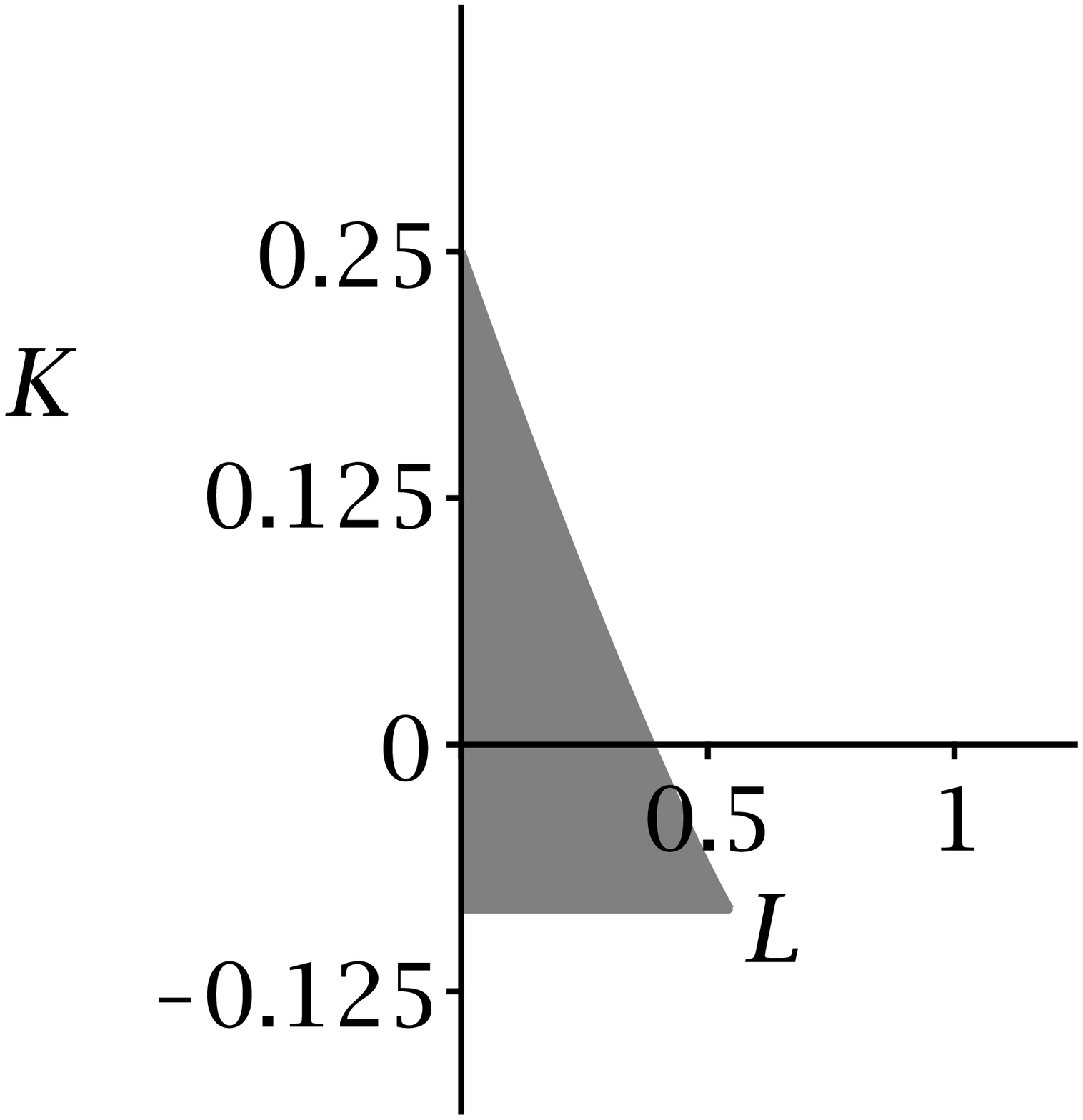}
					\end{subfigure}
					\begin{subfigure}[h]{0.3\linewidth}
						\includegraphics[scale=0.2]{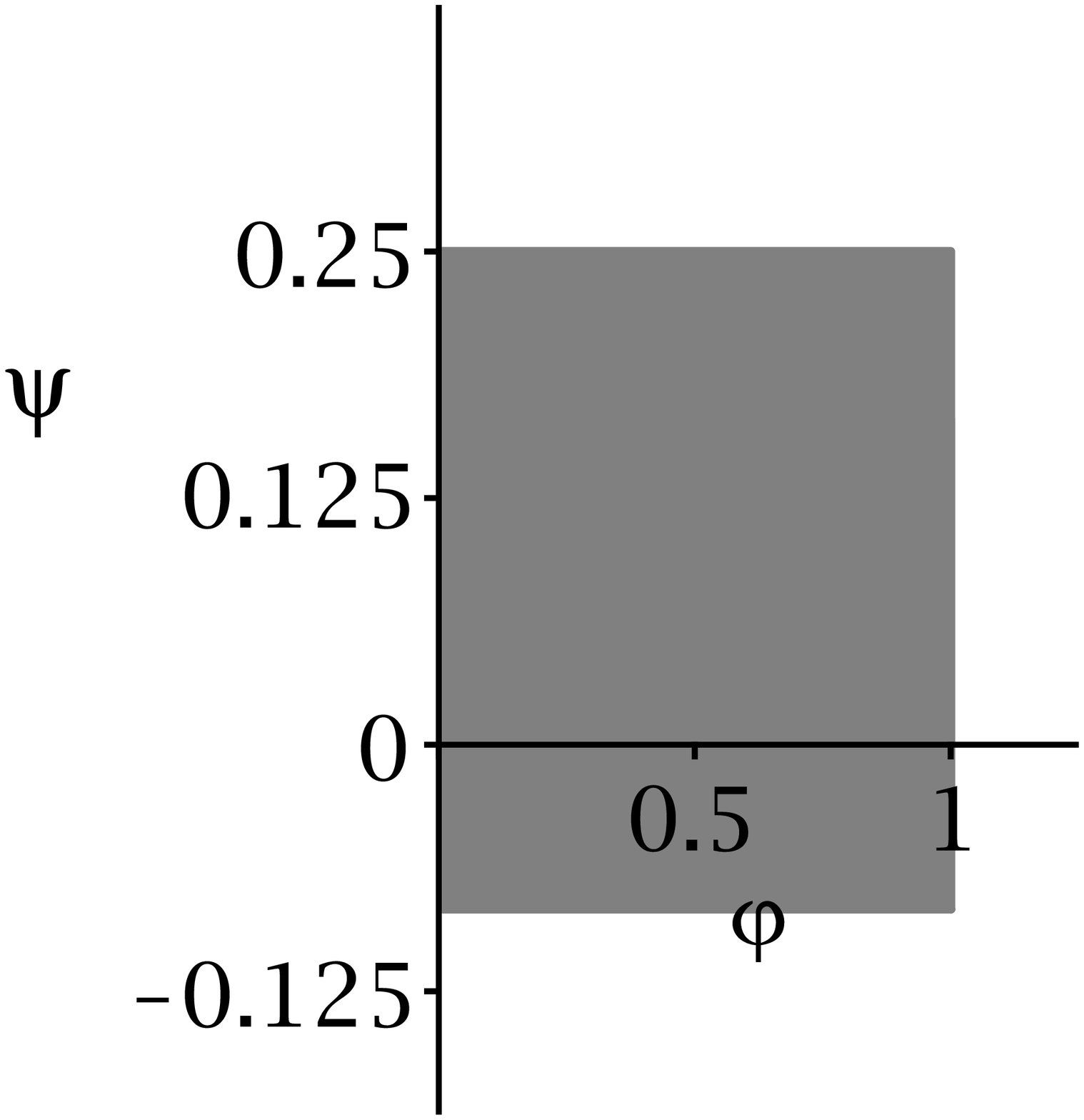}
					\end{subfigure}
					\caption{Domain and codomain of $\widetilde{F}$ \eqref{eqn_first_F_tilda_trafo}.}\label{fig_Ftilda_domain_codomain}
				\end{figure}
		\noindent
		The map $\widetilde{F}$ is a diffeomorphism with 
			\begin{equation}
				\D\widetilde{F}=\left(\begin{matrix}
				\frac{1}{\mathbf{u}(K)} &  \frac{27\sqrt{3}L\left(2-\sqrt{1+12K}\right)}{\sqrt{2\left(1-36K+\sqrt{(1+12K)^3}\right)^3}}\\
				0 & 1
				\end{matrix}
				\right)\label{eqn_Ftilda_differential}
			\end{equation}
		and
			\begin{equation}
				\widetilde{F}^{-1}\left(\left(\begin{smallmatrix}
				\varphi\\ \psi
				\end{smallmatrix}
				\right)\right)=\left(\begin{smallmatrix}
				\varphi \mathbf{u}(K)\\ \psi
				\end{smallmatrix}
				\right).\label{eqn_Ftilda_inverse}
			\end{equation}
		Note at this point that the term $1-36\psi+\sqrt{(1+12\psi)^3}$ is positive for all $\psi\in\left(-\frac{1}{12},\frac{1}{4}\right)$, and in fact it vanishes for $K>-\frac{1}{12}$ if and only if $K=\frac{1}{4}$. The push-forward $\widetilde{F}_*\mathcal{V}$ is of the form
			\begin{align}
				\widetilde{F}_*\mathcal{V}=\left(\widetilde{F}_*\mathcal{V}\right)_{\left(\begin{smallmatrix}
				\varphi\\ \psi
				\end{smallmatrix}
				\right)}&=\frac{3\sqrt{3}\left(1-\varphi^2\right)\left(-1+40\psi-144\psi^2+\left(-1-8\psi+48\psi^2\right)\sqrt{1+12\psi}\right)}{\sqrt{2\left(1-36\psi+\sqrt{(1+12\psi)^3}\right)^3}} \partial_{\varphi}\notag\\
				&+\frac{1}{3\sqrt{6}}\varphi(1+12\psi)\sqrt{1-36\psi+\sqrt{(1+12\psi)^3}}\partial_{\psi}.\label{eqn_Ftilda_V_pushforward}
			\end{align}
		The term $\D\psi\left(\widetilde{F}_*\mathcal{V}\right)$ is positive for all $\left(\begin{smallmatrix}
		\varphi\\ \psi
		\end{smallmatrix}
		\right)\in\left\{-\frac{1}{12}<\psi<\frac{1}{4},\ 0<\varphi<1\right\}$. Thus, the smooth function
			\begin{align}
				&\widetilde{\mathcal{R}}:\left\{-\frac{1}{12}<\psi<\frac{1}{4},\ 0<\varphi<1\right\}\to\mathbb{R},\notag\\ 
				&\widetilde{\mathcal{R}}\left(\left(\begin{smallmatrix}
				\varphi\\ \psi
				\end{smallmatrix}
				\right)\right):=\frac{\D\varphi\left(\widetilde{F}_*\mathcal{V}\right)}{\D\psi\left(\widetilde{F}_*\mathcal{V}\right)}=\frac{1-\varphi^2}{\varphi}\cdot\frac{-27(1-4\psi)}{(1+12\psi)\left(1-36\psi+\sqrt{(1+12\psi)^3}\right)},\label{eqn_Rtilda_def}
			\end{align}
		is well defined. Note that $\widetilde{\mathcal{R}}$ is positive on its domain of definition. Similarly to the definition of $\mathcal{X}$ \eqref{eqn_mathcalX_def} we define the vector field $\mathcal{Y}$ on the set $\left\{-\frac{1}{12}<\psi<\frac{1}{4},\ 0<\varphi<1\right\}\subset\mathbb{R}^2$ as
			\begin{equation}
				\mathcal{Y}:=\widetilde{\mathcal{R}}\,\partial_\varphi+\partial_\psi.\label{eqn_Y_def}
			\end{equation}
		Since $d\psi\left(\widetilde{F}_*\mathcal{V}\right)$ and, hence, $\widetilde{F}_*\mathcal{V}$ do not vanish on the set $\left\{-\frac{1}{12}<\psi<\frac{1}{4},\ 0<\varphi<1\right\}$, it follows that the images of the maximal integral curves of $\mathcal{Y}$ and $\widetilde{F}_*\mathcal{V}$ are in one-to-one correspondence. As for the vector field $\mathcal{X}$ defined in \eqref{eqn_mathcalX_def} we can find a formula for the integral curves of $\mathcal{Y}$. For any open interval $(a,b)\subset\left(-\frac{1}{12},\frac{1}{4}\right)$,
			\begin{align}
				&\mathcal{Y}_\gamma=\dot{\gamma},\quad \gamma(t)=\left(\begin{smallmatrix}
				\varphi(t)\\ t
				\end{smallmatrix}
				\right),\quad \gamma:\left(a,b\right)\to\left\{-\frac{1}{12}<\psi< \frac{1}{4},\  0<\varphi<1\right\},\ \gamma(t_0)=\varphi_0\notag\\
				\Leftrightarrow\quad & \frac{\varphi\dot{\varphi}}{1-\varphi^2}=\underbrace{\frac{-27(1-4t)}{(1+12t)\left(1-36t+\sqrt{(1+12t)^3}\right)}}_{=:J(t)}\label{eqn_J_def}\\ 
				\Leftrightarrow\quad & \varphi(t)=\sqrt{1-c \exp\left(-2\int\limits_{t_0}^tJ(s) \D s\right)},\label{eqn_gamma_t_in_Ftildastuff}
			\end{align}
		where $c>0$ is chosen so that the initial condition $\gamma(t_0)=(\varphi(t_0),t_0)^T$ is met. For an example of such a curve $\gamma$ see Figure \ref{fig_gammatilda_example}.
			\begin{figure}[H]
				\centering
				\includegraphics[scale=0.2]{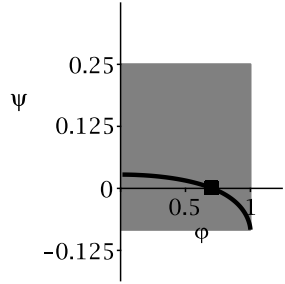}
				\caption{Example of a curve $\gamma$ fulfilling $\mathcal{Y}_\gamma=\dot{\gamma}$ with initial condition $\varphi\left(\frac{1}{36}\right)=\frac{1}{\sqrt{2}}$ marked with a box.}\label{fig_gammatilda_example}
			\end{figure}
		\noindent
		Observe that $\varphi(t)<1$ for all $t\in(a,b)$, and that
			\begin{align}
				J:\left(-\frac{1}{12},\frac{1}{4}\right)\to\mathbb{R},\quad
				J(t)=\frac{-27(1-4t)}{(1+12t)\left(1-36t+\sqrt{(1+12t)^3}\right)},\label{eqn_J_formula}
			\end{align}
		is a well-defined negative function. Note that we can actually find an explicit formula for $\gamma(t)$ as in equation \eqref{eqn_gamma_t_in_Ftildastuff}. For $t_0=\frac{1}{36}$ one can show that $\int\limits_{t_0}^t J(s)ds=-\frac{3}{4}\ln(1+12t)+\frac{1}{2}\ln\left(1-36t+\sqrt{(1+12t)^3}\right)$ and for other values of $t_0$ we might split up the integral to obtain a more general formula. We will however not need the explicit formula for the purpose of this proof. We will now show that $\gamma$ being maximal implies $b<\frac{1}{4}$ and that $\lim\limits_{t\nearrow b}\varphi(t)=0$. This will imply that the corresponding maximal integral curve of $\mathcal{V}|_{\left\{-\frac{1}{12}<K<\frac{1}{4},\, 0<L<\mathbf{u}(K)\right\}}$ has a limit point in $\left\{-\frac{1}{12}<K<\frac{1}{4},\ L=0\right\}$ bounded away from $K=\frac{1}{4}$ and, hence, that each maximal integral curve of $\mathcal{V}|_{\left\{-\frac{1}{12}<K<\frac{1}{4},\, |L|<\mathbf{u}(K)\right\}}$ meets $\left\{-\frac{1}{12}<K<\frac{1}{4},\ L=0\right\}$ in one point. Recall for this step that $\left\{-\frac{1}{12}<K<\frac{1}{4},\, |L|<\mathbf{u}(K)\right\}$ is contained in $\mathbb{R}^2\setminus\{\mathcal{V}=0\}$. To prove $b<\frac{1}{4}$ it suffices to show that for all $c>0$,
			\begin{equation}\label{eqn_varphi_upper_triangle_0}
				\varphi(t)=\sqrt{1-c \exp\left(-2\int\limits_{t_0}^tJ(s)\D s\right)}=0
			\end{equation}
		is fulfilled for some $t\in\left(a,\frac{1}{4}\right)$. Since $J(s)<0$ for all $s\in\left(-\frac{1}{12},\frac{1}{4}\right)$ and the term $-2\int\limits_{t_0}^tJ(s)\D s$ is thus strictly monotonously increasing in $t$, it is sufficient to show that for all $t_0\in\left(-\frac{1}{12},\frac{1}{4}\right)$
			\begin{equation}\label{eqn_J_int_infinity}
				\lim\limits_{t\nearrow1/4}\int\limits_{t_0}^tJ(s)\D s=-\infty.
			\end{equation}
		Since $J$ is smooth on $\left(-\frac{1}{12},\frac{1}{4}\right)$, we can without loss of generality assume $t_0=0$. Observe that $(1+12t)\in\left(1,4\right)$ for $t\in\left[0,\frac{1}{4}\right]$ implies that \eqref{eqn_J_int_infinity} is equivalent to
			\begin{equation}\label{eqn_J_int_1stequiv_infinity}
				\lim\limits_{t\nearrow1/4}\int\limits_{0}^t\frac{1-4s}{1-36s+\sqrt{(1+12s)^3}}\D s=\infty.
			\end{equation}
		We substitute $s$ by $-r+\frac{1}{4}$ in \eqref{eqn_J_int_1stequiv_infinity} and see that it is equivalent to
			\begin{equation}\label{eqn_J_int_2ndequiv_infinity}
				\lim\limits_{t\searrow 0}\int\limits_{t}^{\frac{1}{4}}\frac{r}{9r-2+(-6r+2)\sqrt{1-3r}}\D r=\infty.
			\end{equation}
		Note that both the numerator and denominator of the integrand in \eqref{eqn_J_int_2ndequiv_infinity} are positive and smooth on the interval $\left(0,\frac{1}{4}\right)$, and both converge to $0$ as $r\to 0$. To prove \eqref{eqn_J_int_2ndequiv_infinity} it is enough to show that there exists $\varepsilon\in\left(0,\frac{1}{4}\right)$ and $A>0$, such that for all $r\in (0,\varepsilon)$
			\begin{equation}\label{eqn_J_int_3rdequiv_infinity}
				9r-2+(-6r+2)\sqrt{1-3r}\leq Ar^2.
			\end{equation}
		The condition \eqref{eqn_J_int_3rdequiv_infinity} on the other hand can be proven by showing that
			\begin{equation}\label{eqn_J_int_4thequiv_infinity}
				\lim\limits_{r\searrow 0}\frac{ 9r-2+(-6r+2)\sqrt{1-3r}}{r^2}
			\end{equation}
		exists and is positive. Using L'H\^opitals rule for limits yields
			\begin{equation*}
				\lim\limits_{r\searrow 0}\frac{ 9r-2+(-6r+2)\sqrt{1-3r}}{r^2}=\lim\limits_{r\searrow 0}\frac{9\left(1-\sqrt{1-3r}\right)}{2r}=\frac{27}{4}.
			\end{equation*}
		Hence, \eqref{eqn_J_int_4thequiv_infinity} holds true, and since \eqref{eqn_J_int_4thequiv_infinity} $\Rightarrow$ \eqref{eqn_J_int_3rdequiv_infinity} $\Rightarrow$ \eqref{eqn_J_int_2ndequiv_infinity} $\Rightarrow$ \eqref{eqn_J_int_1stequiv_infinity} $\Rightarrow$ \eqref{eqn_J_int_infinity}, we have proven that for all initial values $\varphi_0\in(0,1)$ and corresponding $c>0$, there exists $t=\widetilde{t}\in\left(a,b\right)$, $\widetilde{t}>t_0$, such that equation \eqref{eqn_varphi_upper_triangle_0} is fulfilled. Summarising, we have shown that each maximal integral curve of $\mathcal{V}|_{\mathbb{R}^2\setminus\{\mathcal{V}=0\}}$ starting in $\left\{-\frac{1}{12}<K<\frac{1}{4},\ |L|<\mathbf{u}(K)\right\}$ meets the set $\left\{-\frac{1}{12}<K<\frac{1}{4},\ L=0\right\}$ in one point. We have already shown that the quartic CCGPSR curves $\mathcal{H}_{0,K}$ for $K\leq \frac{1}{4}$ are pairwise inequivalent and can thus deduce that this point is unique. This proves that every maximal connected quartic GPSR curve $\mathcal{H}_{L,K}$ with $(L,K)^T\in\left\{-\frac{1}{12}<K<\frac{1}{4},\ |L|<\mathbf{u}(K)\right\}$ is equivalent to a uniquely determined quartic CCGPSR curve $\mathcal{H}_{0,\widetilde{K}}$, $\widetilde{K}\in\left(-\frac{1}{12},\frac{1}{4}\right)$, and thus in particular itself a quartic CCGPSR curve.
		
		We have shown up to this point that the maximal connected quartic GPSR curves cor\-res\-pon\-ding to \hyperref[eqn_qCCPSRcurves_class_a]{$a)$}, \hyperref[eqn_qCCPSRcurves_class_c]{$c)$}, and \hyperref[eqn_qCCPSRcurves_class_d]{$d)$} are closed and, hence, quartic CCGPSR curves. The remaining case \hyperref[eqn_qCCPSRcurves_class_b]{$b)$} corresponds to the point
			\begin{equation*}
				\left(\begin{matrix}
				L\\ K
				\end{matrix}
				\right)=\left(\begin{matrix}
				\frac{2\sqrt{2}}{3\sqrt{3}}\\ -\frac{1}{12}
				\end{matrix}
				\right)\in\{\mathcal{V}=0\},
			\end{equation*}
		and we will now show that the corresponding maximal connected quartic GPSR curve $\mathcal{H}_{\frac{2\sqrt{2}}{3\sqrt{3}}, -\frac{1}{12}}$ is closed. Note that $\mathcal{H}_{\frac{2\sqrt{2}}{3\sqrt{3}}, -\frac{1}{12}}$ and $\mathcal{H}_{-\frac{2\sqrt{2}}{3\sqrt{3}}, -\frac{1}{12}}$ are equivalent, and so this is indeed the last remaining maximal connected quartic GPSR curve we have to study in this proof. To do so, it suffices to show that $\mathcal{H}_{\frac{2\sqrt{2}}{3\sqrt{3}}, -\frac{1}{12}}$ is a homogeneous space. Riemannian homogeneous spaces are automatically geodesically complete and we can then use \cite[Prop.\,1.8]{CNS} to conclude that $\mathcal{H}_{\frac{2\sqrt{2}}{3\sqrt{3}}, -\frac{1}{12}}\subset\mathbb{R}^2$ is closed and thus a quartic CCGPSR curve. In fact, $\left(\frac{2\sqrt{2}}{3\sqrt{3}},-\frac{1}{12}\right)^T\in\{\mathcal{V}=0\}$ and the formulas \eqref{eqn_VFqCCGPSRc_def}, \eqref{eqn_deltaP3_quarticcurves}, and \eqref{eqn_deltaP4_quarticcurves} show that for $\mathcal{H}_{\frac{2\sqrt{2}}{3\sqrt{3}}, -\frac{1}{12}}$ we have $\delta  P_3=\delta P_4= 0$, which implies with \cite[Prop.\,3.12]{L2} the homogeneity of $\mathcal{H}_{\frac{2\sqrt{2}}{3\sqrt{3}}, -\frac{1}{12}}$ as claimed. It remains to determine the connected components of $\left\{h_{\frac{2\sqrt{2}}{3\sqrt{3}}, -\frac{1}{12}}=1\right\}\cap\left\{\text{hyperbolic points of }h_{\frac{2\sqrt{2}}{3\sqrt{3}}, -\frac{1}{12}}\right\}$ and show that they are, as quartic CCGPSR curves, equivalent. We find
			\begin{equation*}
				h_{\frac{2\sqrt{2}}{3\sqrt{3}},-\frac{1}{12}}\left(\left(\begin{smallmatrix}
				1\\ y
				\end{smallmatrix}\right)\right)=\frac{1}{12}\left(y-\sqrt{6}\right)^3\left(y+\frac{\sqrt{2}}{\sqrt{3}}\right)
			\end{equation*}
		and
			\begin{equation*}
				h_{\frac{2\sqrt{2}}{3\sqrt{3}},-\frac{1}{12}}\left(\left(\begin{smallmatrix}
				x\\ 1
				\end{smallmatrix}\right)\right)=\left(x-\frac{1}{\sqrt{6}}\right)^3\left(x+\frac{\sqrt{3}}{\sqrt{2}}\right).
			\end{equation*}
		Hence,
			\begin{align*}
				\left\{h_{\frac{2\sqrt{2}}{3\sqrt{3}},-\frac{1}{12}}>0\right\}&=\mathbb{R}_{>0}\cdot\left\{\left(\begin{smallmatrix}
				1\\ y
				\end{smallmatrix}
				\right)\in\mathbb{R}^2\ \left|\  y\in\left(-\tfrac{\sqrt{2}}{\sqrt{3}},\sqrt{6}\right)\right\}\right.\notag\\
				&\ \dot\cup\ \mathbb{R}_{>0}\cdot \left\{-\left(\begin{smallmatrix}
				1\\ y
				\end{smallmatrix}
				\right)\in\mathbb{R}^2\ \left|\  y\in\left(-\tfrac{\sqrt{2}}{\sqrt{3}},\sqrt{6}\right)\right\}\right.,
			\end{align*}
		that is $\left\{h_{\frac{2\sqrt{2}}{3\sqrt{3}},-\frac{1}{12}}>0\right\}$ has precisely two connected components, both of which only contain hyperbolic points and each a unique quartic CCGPSR curve. These two curves are the connected components of $\left\{h_{\frac{2\sqrt{2}}{3\sqrt{3}}, -\frac{1}{12}}=1\right\}$, and they are equivalent via $\left(\begin{smallmatrix}
		x\\ y
		\end{smallmatrix}
		\right)\to-\left(\begin{smallmatrix}
		x\\ y
		\end{smallmatrix}
		\right)$. Note that $\mathcal{H}_{\frac{2\sqrt{2}}{3\sqrt{3}},-\frac{1}{12}}\subset\mathbb{R}_{>0}\cdot\left\{\left(\begin{smallmatrix}
		1\\ y
		\end{smallmatrix}
		\right)\in\mathbb{R}^2\ \left|\ y\in\left(-\frac{\sqrt{2}}{\sqrt{3}},\sqrt{6}\right)\right\}\right.$. 
		
		In order to find the automorphism group $G^h$ of $h_{\frac{2\sqrt{2}}{3\sqrt{3}}, -\frac{1}{12}}$, we now only need to determine $G^{h_{\frac{2\sqrt{2}}{3\sqrt{3}}, -\frac{1}{12}}}_0$ and check that there are no additional discrete symmetries of $h_{\frac{2\sqrt{2}}{3\sqrt{3}}, -\frac{1}{12}}$ mapping $\mathcal{H}_{\frac{2\sqrt{2}}{3\sqrt{3}}, -\frac{1}{12}}\subset\mathbb{R}^2$ to itself. For $G^{h_{\frac{2\sqrt{2}}{3\sqrt{3}}, -\frac{1}{12}}}_0$ consider, similar to $h_{0,\frac{1}{4}}$ and $G^{h_{0,\frac{1}{4}}}_0$, equation \eqref{eqn_A_pulled_back_to_dom} and calculate the derivative of the map $\mathcal{A}:\mathrm{dom}\left(\mathcal{H}_{\frac{2\sqrt{2}}{3\sqrt{3}}, -\frac{1}{12}}\right)\to\mathrm{GL}(2)$ at $z=0\in\mathrm{dom}\left(\mathcal{H}_{\frac{2\sqrt{2}}{3\sqrt{3}}, -\frac{1}{12}}\right)$, cf. \eqref{eqn_differential_A_formula}. The linear map $b\in\mathrm{Lin}(\mathbb{R},\mathfrak{so}(1))$ in \eqref{eqn_differential_A_formula} vanishes identically since $\dim(\mathfrak{so}(1))=0$. We obtain
			\begin{equation*}
				\D\mathcal{A}_0(\partial_y)=\left(\begin{matrix}
				0 & \frac{1}{2}\\
				1 & \frac{\sqrt{2}}{\sqrt{3}}
				\end{matrix}
				\right).
			\end{equation*}
		With $\widetilde{a}:=\sqrt{2}\left(\begin{matrix}
		0 & \frac{1}{2}\\
		1 & \frac{\sqrt{2}}{\sqrt{3}}
		\end{matrix}
		\right)$, one can check that
			\begin{equation*}
				\left.\D h_{\frac{2\sqrt{2}}{3\sqrt{3}}, -\frac{1}{12}}\right|_{\left(\begin{smallmatrix}
				x\\ y
				\end{smallmatrix}
				\right)}\left(\widetilde{a}\cdot\left(\begin{smallmatrix}
				x\\ y
				\end{smallmatrix}
				\right)\right)\equiv 0
			\end{equation*}
		as expected. Explicitly, the induced action is given by
			\begin{equation}\label{eqn_nontriv_R_action_b}
				\exp(t\widetilde{a})=\left(\begin{matrix}
					1 & \tfrac{\sqrt{3}}{\sqrt{2}}\\
					1 & -\tfrac{1}{\sqrt{6}}
				\end{matrix}\right)\left(\begin{matrix}
					e^{\sqrt{3}t} & 0\\
					0 & e^{-\tfrac{t}{\sqrt{3}}}
				\end{matrix}\right)\left(\begin{matrix}
					\tfrac{1}{4} & \tfrac{3}{4}\\
					\tfrac{\sqrt{3}}{2\sqrt{2}} & -\tfrac{\sqrt{3}}{2\sqrt{2}}
				\end{matrix}\right).
			\end{equation}
		The quartic CCGPSR curve thus fulfils $\mathcal{H}_{\frac{2\sqrt{2}}{3\sqrt{3}}, -\frac{1}{12}}\cong\mathbb{R}$ as a Riemannian homogeneous space via the corresponding action of $G^{h_{\frac{2\sqrt{2}}{3\sqrt{3}}, -\frac{1}{12}}}_0$. Now, again similar to the cases $\mathcal{H}_{0,K}$ with $K\leq \frac{1}{4}$, we need to find all $A\in G^{h_{\frac{2\sqrt{2}}{3\sqrt{3}}, -\frac{1}{12}}}\subset\mathrm{GL}(2)$ which are not contained in $G^{h_{\frac{2\sqrt{2}}{3\sqrt{3}}, -\frac{1}{12}}}_0$, such that $A\left(\mathcal{H}_{\frac{2\sqrt{2}}{3\sqrt{3}}, -\frac{1}{12}}\right)=\mathcal{H}_{\frac{2\sqrt{2}}{3\sqrt{3}}, -\frac{1}{12}}$. With
			\begin{equation*}
				A=\left(\begin{matrix}
				a_{11} & a_{12} \\
				a_{21} & a_{22}
				\end{matrix}
				\right)
			\end{equation*}
		it find that $a_{11}\ne 0$, since otherwise
			\begin{equation*}
				h_{\frac{2\sqrt{2}}{3\sqrt{3}}, -\frac{1}{12}}\left(\left(\begin{smallmatrix}
				a_{11} \\ a_{21}
				\end{smallmatrix}
				\right)x\right)=h_{\frac{2\sqrt{2}}{3\sqrt{3}},  -\frac{1}{12}}\left(\left(\begin{smallmatrix}
				0 \\ a_{21}
				\end{smallmatrix}
				\right)x\right)=-\frac{1}{12}a_{21}^4 x^4,
			\end{equation*}
		but $h_{\frac{2\sqrt{2}}{3\sqrt{3}}, -\frac{1}{12}}\left(\left(\begin{smallmatrix}
		a_{11} \\ a_{21}
		\end{smallmatrix}
		\right)x\right)=x^4$ is a necessary requirement for $A$ to be an automorphism of $h_{\frac{2\sqrt{2}}{3\sqrt{3}}, -\frac{1}{12}}$. Furthermore, $\left(\begin{smallmatrix}
		1\\ 0
		\end{smallmatrix}
		\right)\in\mathcal{H}_{\frac{2\sqrt{2}}{3\sqrt{3}}, -\frac{1}{12}}$ is mapped to $\left(\begin{smallmatrix}
		a_{11} \\ a_{21}
		\end{smallmatrix}
		\right)$, which is required to be an element of $\mathcal{H}_{\frac{2\sqrt{2}}{3\sqrt{3}}, -\frac{1}{12}}$. The hyperbolicity of $\mathcal{H}_{\frac{2\sqrt{2}}{3\sqrt{3}}, -\frac{1}{12}}\subset\mathbb{R}^2$ then implies that $a_{11}\geq 1$. Since $\left(\begin{smallmatrix}
		a_{11} \\ a_{21}
		\end{smallmatrix}
		\right)\in\mathcal{H}_{\frac{2\sqrt{2}}{3\sqrt{3}}, -\frac{1}{12}}$, $A$ must thus be of the form \eqref{eqn_A_explicit}. With the fact that $G^{h_{\frac{2\sqrt{2}}{3\sqrt{3}}, -\frac{1}{12}}}_0$ acts transitively on $\mathcal{H}_{\frac{2\sqrt{2}}{3\sqrt{3}}, -\frac{1}{12}}$, this shows that $A$ can be written as
			\begin{equation*}
				A=A_0\cdot \widetilde{A},
			\end{equation*}
		where $A_0\in G^{h_{\frac{2\sqrt{2}}{3\sqrt{3}}, -\frac{1}{12}}}_0$ and $\widetilde{A}\in G^{h_{\frac{2\sqrt{2}}{3\sqrt{3}}, -\frac{1}{12}}}$ is contained in the stabilizer of the point $\left(\begin{smallmatrix}
		1\\ 0
		\end{smallmatrix}
		\right)\in\mathcal{H}_{\frac{2\sqrt{2}}{3\sqrt{3}}, -\frac{1}{12}}$. Hence, we need to determine all $\widetilde{A}\in G^{h_{\frac{2\sqrt{2}}{3\sqrt{3}}, -\frac{1}{12}}}$, such that $\widetilde{A}\cdot\left(\begin{smallmatrix}
		1\\ 0
		\end{smallmatrix}
		\right)=\left(\begin{smallmatrix}
		1\\ 0
		\end{smallmatrix}
		\right)$. $\widetilde{A}$ must be of the form
			\begin{equation*}
				\widetilde{A}=\left(\begin{matrix}
				1 & \left.-\frac{\partial_y h_{\frac{2\sqrt{2}}{3\sqrt{3}},  -\frac{1}{12}}}{\partial_x h_{\frac{2\sqrt{2}}{3\sqrt{3}}, -\frac{1}{12}}}\right|_{\left(\begin{smallmatrix}
				1\\ 0
				\end{smallmatrix}
				\right)} \cdot r\\
				0 & r
				\end{matrix}
				\right)=\left(\begin{matrix}
				1 & 0\\
				0 & r
				\end{matrix}
				\right)
			\end{equation*}
		for some $r\ne 0$. Then
			\begin{equation*}
				h_{\frac{2\sqrt{2}}{3\sqrt{3}}, -\frac{1}{12}}\left(\widetilde{A}\cdot\left(\begin{smallmatrix}
				x \\ y
				\end{smallmatrix}
				\right)\right)=x^4-r^2x^2y^2+r^3\frac{2\sqrt{2}}{3\sqrt{3}}xy^3-\frac{r^4}{12}y^4,
			\end{equation*}
		which shows that $r=1$. Summarising, we have shown that
			\begin{equation*}
				G^{h_{\frac{2\sqrt{2}}{3\sqrt{3}}, -\frac{1}{12}}}\cong \mathbb{R}\times\mathbb{Z}_2,
			\end{equation*}
		where $\mathbb{R}$ acts as described in \eqref{eqn_nontriv_R_action_b} and $\mathbb{Z}_2$ acts via $\left(\begin{smallmatrix}
		x\\ y
		\end{smallmatrix}
		\right)\mapsto -\left(\begin{smallmatrix}
		x\\ y
		\end{smallmatrix}
		\right)$.
		
		In order to complete the proof of Theorem \ref{thm_quartic_CCGPSR_curves_classification}, we still need to prove that the quartic CCGPSR curves \hyperref[eqn_qCCPSRcurves_class_a]{$a)$}, \hyperref[eqn_qCCPSRcurves_class_b]{$b)$}, \hyperref[eqn_qCCPSRcurves_class_c]{$c)$}, and elements in the family of curves \hyperref[eqn_qCCPSRcurves_class_d]{$d)$} are pairwise inequivalent. We already have seen that this is true if one considers two elements in the one-parameter family \hyperref[eqn_qCCPSRcurves_class_d]{$d)$}. Since the quartic CCGPSR curve \hyperref[eqn_qCCPSRcurves_class_a]{$a)$}, that is $\mathcal{H}_{0,\frac{1}{4}}$, has a transitive action of the corresponding Lie group $G^{h_{0,\frac{1}{4}}}_0$ it might only be equivalent to the quartic CCGPSR curve $\mathcal{H}_{\frac{2\sqrt{2}}{3\sqrt{3}},-\frac{1}{12}}$, that is \hyperref[eqn_qCCPSRcurves_class_b]{$b)$}, which also has a transitive $G^{h_{\frac{2\sqrt{2}}{3\sqrt{3}},-\frac{1}{12}}}_0$-action. But with
			\begin{equation*}
				\mathrm{dom}\left(\mathcal{H}_{0,\frac{1}{4}}\right)=\left(-\sqrt{2},\sqrt{2}\right)
			\end{equation*}
		and
			\begin{equation*}
				\mathrm{dom}\left(\mathcal{H}_{\frac{2\sqrt{2}}{3\sqrt{3}},-\frac{1}{12}}\right)=\left(-\frac{\sqrt{2}}{\sqrt{3}},\sqrt{6}\right)
			\end{equation*}
		we find that
			\begin{equation*}
				\left.\D h_{0,\frac{1}{4}}\right|_{\left(\begin{smallmatrix}
				1\\ -\sqrt{2}
				\end{smallmatrix}
				\right)}=\left.\D h_{0,\frac{1}{4}}\right|_{\left(\begin{smallmatrix}
				1\\ \sqrt{2}
				\end{smallmatrix}
				\right)}=0
			\end{equation*}
		and
			\begin{equation*}
				\left.\D  h_{\frac{2\sqrt{2}}{3\sqrt{3}},-\frac{1}{12}}\right|_{\left(\begin{smallmatrix}
				1\\ -\frac{\sqrt{2}}{\sqrt{3}}
				\end{smallmatrix}
				\right)}= \frac{64}{27}\D x+\frac{32\sqrt{2}}{9\sqrt{3}}\D y\ne0,\quad \left. \D h_{\frac{2\sqrt{2}}{3\sqrt{3}},-\frac{1}{12}}\right|_{\left(\begin{smallmatrix}
				1\\ \sqrt{6}
				\end{smallmatrix}
				\right)}=0.
			\end{equation*}
		This means that $\D h_{0,\frac{1}{4}}$ vanishes on $\partial\left(\mathbb{R}_{>0}\cdot\mathcal{H}_{0,\frac{1}{4}}\right)$, but $\D h_{\frac{2\sqrt{2}}{3\sqrt{3}},-\frac{1}{12}}$ vanishes only on one of the two connected components of $\partial\left(\mathbb{R}_{>0}\cdot\mathcal{H}_{\frac{2\sqrt{2}}{3\sqrt{3}},-\frac{1}{12}}\right)\setminus\{0\}$. Hence, the quartic CCGPSR curves \hyperref[eqn_qCCPSRcurves_class_a]{a)} and \hyperref[eqn_qCCPSRcurves_class_b]{b)} cannot be equivalent. Note that with the above calculations we have also shown that both curves in Thm. \ref{thm_quartic_CCGPSR_curves_classification} \hyperref[eqn_qCCPSRcurves_class_b]{a)} and \hyperref[eqn_qCCPSRcurves_class_b]{b)} are singular at infinity.  If we cared only for inequivalence, we could have used that
		$G^{h_{\frac{2\sqrt{2}}{3\sqrt{3}}, -\frac{1}{12}}}$ has $2$ connected components, while $G^{h_{0,\frac{1}{4}}}$ has $8$ connected components. Now, in order to prove that there exist no quartic CCGPSR curve $\mathcal{H}_{0,K}$ in the one-parameter family \hyperref[eqn_qCCPSRcurves_class_d]{d)} that is equivalent to the quartic CCGPSR curve $\mathcal{H}_{\frac{2}{3\sqrt{3}},0}$, that is \hyperref[eqn_qCCPSRcurves_class_c]{c)}, we will use a similar argument. We find that for $K<0$, the equation
			\begin{equation*}
				\D h_{0,K}\left(\partial_y\right)=y\left(-2+4Ky^2\right)
			\end{equation*}
		has no solutions other than $y=0$. For $0<K<\frac{1}{4}$,
			\begin{equation*}
				\D h_{0,K}\left(\partial_y\right)=0,\ y\ne 0\quad\Leftrightarrow\quad y=\pm\frac{1}{\sqrt{2K}}.
			\end{equation*}
		But then
			\begin{equation*}
				h_{0,K}\left(\left(\begin{smallmatrix}
				1\\ \pm\frac{1}{\sqrt{2K}}
				\end{smallmatrix}\right)\right)=1-\frac{1}{4K}<0.
			\end{equation*}
		Hence, for all $K\in\left(0,\frac{1}{4}\right)$ the points $\pm\frac{1}{\sqrt{2K}}$ are not contained in $\partial\left(\mathrm{dom}\left(\mathcal{H}_{0,K}\right)\right)$. Summarising, we have shown that for all $K<\frac{1}{4}$, $\left.\D h_{0,K}\right|_{\partial\left(\mathbb{R}_{>0}\cdot\mathcal{H}_{0,K}\right)}$ is nowhere zero, meaning that $\mathcal{H}_{0,K}$ is not singular at infinity for any $K<\frac{1}{4}$. But $\mathcal{H}_{\frac{2}{3\sqrt{3}},0}$ is singular at infinity. Explicitly, we find $\mathrm{dom}\left(\mathcal{H}_{\frac{2}{3\sqrt{3}},0}\right)=\left(-\frac{\sqrt{3}}{2},\sqrt{3}\right)$ and
			\begin{equation*}
				\left.\D h_{\frac{2}{3\sqrt{3}},0}\right|_{\left(\begin{smallmatrix}
				1\\ \sqrt{3}
				\end{smallmatrix}
				\right)}\equiv 0.
			\end{equation*}
		Hence, for all $K<\frac{1}{4}$ the quartic CCGPSR curves $\mathcal{H}_{0,K}$ and $\mathcal{H}_{\frac{2}{3\sqrt{3}},0}$ cannot be equivalent.
		
		This finishes the proof of Theorem \ref{thm_quartic_CCGPSR_curves_classification}.
	\end{proof}
	\end{Th}
	
	Next, we will classify all maximal non-closed quartic GPSR curves up to equivalence. The maximal non-closed quartic GPSR curves are precisely the maximal incomplete quartic GPSR curves since in dimension one being closed is equivalent to geodesic completeness for all (G)PSR curves by \cite[Thm.\,2.9]{CNS}.
	
	\begin{Th}\label{thm_incomplete_quartics}
		Let $\mathcal{H}\subset\{h=1\}$ be a maximal connected incomplete quartic GPSR curve. Let further
			\begin{equation*}
				U:=\frac{1}{24}\left(4+\left(31+3\sqrt{57}\right)\left(46+6\sqrt{57}\right)^{-1/3}+\left(100+12\sqrt{57}\right)\left(46+6\sqrt{57}\right)^{-2/3}\right).
			\end{equation*}
		Then $\mathcal{H}$ is equivalent to precisely one of the following maximal incomplete quartic GPSR curves in standard form:
			\begin{enumerate}[a)]
				\item $h=x^4-x^2y^2+Ky^4$, $K>\frac{1}{4}$, $\{h=1\}\cap\{\text{hyperbolic points of }h\}$ has $4$ equivalent connected components,\label{eqn_incomplete_qGPSRcurves_class_a} and
					\begin{equation*}
						G^h\cong\mathbb{Z}_2\ltimes\mathbb{Z}_2,
					\end{equation*}
				\item $h=x^4-x^2y^2+xy^3+Ky^4=1$, $K\in\left(-\frac{25}{72},U\right)$, $\{h=1\}\cap\{\text{hyperbolic points of }h\}$ has $4$ equivalent connected components,\label{eqn_incomplete_qGPSRcurves_class_b} and
					\begin{equation*}
						G^h\cong\mathbb{Z}_2\times\mathbb{Z}_2,
					\end{equation*}
				\item $h=x^4-x^2y^2+xy^3-\frac{25}{72}y^4=1$, $\{h=1\}\cap\{\text{hyperbolic points of }h\}$ has $4$ equivalent connected components,\label{eqn_incomplete_qGPSRcurves_class_c} and
					\begin{equation*}
						G^h\cong\mathbb{Z}_2\times\mathbb{Z}_2,
					\end{equation*}
				\item $h=x^4-x^2y^2+xy^3+Uy^4=1$, $\{h=1\}\cap\{\text{hyperbolic points of }h\}$ has $4$ equivalent connected components,\label{eqn_incomplete_qGPSRcurves_class_d} and
					\begin{equation*}
						G^h\cong\mathbb{Z}_2\times\mathbb{Z}_2.
					\end{equation*}
			\end{enumerate}
		The curve \hyperref[eqn_incomplete_qGPSRcurves_class_d]{d)} is singular at infinity, and the curves \hyperref[eqn_incomplete_qGPSRcurves_class_a]{a)}, \hyperref[eqn_incomplete_qGPSRcurves_class_b]{b)}, and \hyperref[eqn_incomplete_qGPSRcurves_class_c]{c)} are not singular at infinity.
		\begin{proof}
			Firstly note that we could have combined the polynomials in \hyperref[eqn_incomplete_qGPSRcurves_class_b]{b)}, \hyperref[eqn_incomplete_qGPSRcurves_class_c]{c)}, and \hyperref[eqn_incomplete_qGPSRcurves_class_d]{d)} into one family by allowing $K\in\left[-\frac{25}{72},U\right]$ in \hyperref[eqn_incomplete_qGPSRcurves_class_b]{b)}. The boundary cases have however unique properties regarding their geometry and their asymptotic behaviour which we will see in Proposition \ref{prop_limit_geos} and this proof.

			Recall the definition of $\mathbf{u}:\left(-\frac{1}{12},\frac{1}{4}\right)\to\mathbb{R}$ in \eqref{eqn_upper_bdr}, and that $\mathbf{u}$ can be continuously extended to the closure of its domain. We have seen that the graph of $\mathbf{u}$ coincides with the image of a maximal integral curve of the vector field $\mathcal{V}=\left(\frac{9}{2}L^2+4K-1\right)\partial_L + L\left(6K+\frac{1}{2}\right)\partial_K$ restricted to $\{\mathcal{V}\ne0\}$, cf. equation \eqref{eqn_partial_u} and the subsequent discussion. Observe that $\mathbf{u}$ is also well defined on $\left(\frac{1}{4},\infty\right)$. This follows from the fact that
				\begin{equation*}
					1-36K+\sqrt{(1+12K)^3}>0
				\end{equation*}
			for all $K>\frac{1}{4}$, which can quickly be checked with the help of any computer algebra system. The vector field $\mathcal{V}$ vanishes at no point in the graph of $\mathbf{u}:\left(\frac{1}{4},\infty\right)\to\mathbb{R}$. Hence, by equation \eqref{eqn_partial_u}, which is symbolically the same for $K>\frac{1}{4}$, it follows that the graph of  $\mathbf{u}:\left(\frac{1}{4},\infty\right)\to\mathbb{R}$ is the image of a maximal integral curve of $\mathcal{V}|_{\{\mathcal{V}\ne0\}}$. Recall at this point that the graph of $\mathbf{w}:\left(-\infty,-\frac{1}{12}\right)\to\mathbb{R}$ \eqref{eqn_sharp_bdr} is also the image of a maximal integral curve of $\mathcal{V}|_{\{\mathcal{V}\ne0\}}$, cf. \eqref{eqn_partial_w}, and keep the fact that $\mathbf{w}\left(-\frac{25}{72}\right)=1$ in mind. Lastly, observe that
				\begin{equation*}
					\mathcal{V}|_{\left\{L>\frac{2\sqrt{2}}{3\sqrt{3}},\ K=-\frac{1}{12}\right\}}=\left(\frac{9}{2}L^2+\frac{2}{3}\right)\partial_L
				\end{equation*}
			together with $\mathcal{V}|_{\left(\begin{smallmatrix}\frac{2\sqrt{2}}{3\sqrt{3}}\\ -\frac{1}{12}\end{smallmatrix}\right)}=0$ implies with
				\begin{equation*}
					\mathbf{m}:\left(\frac{2\sqrt{2}}{3\sqrt{3}},\infty\right)\to\mathbb{R},\quad L\mapsto -\frac{1}{12},
				\end{equation*}
			that the image of $L\mapsto\left(\mathbf{m}(L),L\right)^T$, $L\in\left(\frac{2\sqrt{2}}{3\sqrt{3}},\infty\right)$, coincides with the image of a maximal integral curve of $\mathcal{V}|_{\{\mathcal{V}\ne0\}}$.
			
			In the following we will denote the restriction of $\mathbf{u}$ to $\left(-\frac{1}{12},\frac{1}{4}\right)$ with $\mathbf{u}_{\text{lower}}$, and the restriction of $\mathbf{u}$ to $\left(\frac{1}{4},\infty\right)$ with $\mathbf{u}_{\text{upper}}$. Note that $\mathbf{w}^{-1}(1)=-\frac{25}{72}$ and that the constant $U\approx 1,054784062$ defined in the statement of this theorem fulfils $U=\mathbf{u}_{\mathrm{upper}}^{-1}(1)$. See Figure \ref{fig_u_lower_u_upper_w_m_gen_set} for an overview of the functions we will be working with in this proof together with the set
				\begin{equation*}
					I:=\left\{\left(\begin{smallmatrix}1\\ K\end{smallmatrix}\right)\ \left|\ K\in\left(-\tfrac{25}{72},U\right)\right\}\right.,
				\end{equation*}
			whose closure corresponds to the polynomials in Thm. \ref{thm_incomplete_quartics} \hyperref[eqn_incomplete_qGPSRcurves_class_b]{b)}--\hyperref[eqn_incomplete_qGPSRcurves_class_d]{d)}.
				\begin{figure}[H]
					\centering
					\includegraphics[scale=0.2]{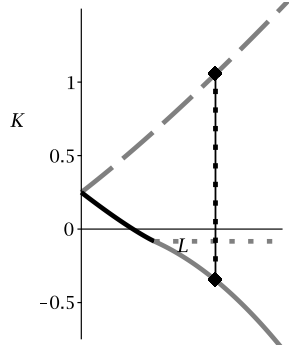}
					\caption{Graphs of $\mathbf{u}_{\mathrm{lower}}$ (black line), $\mathbf{u}_{\mathrm{upper}}$ (dashed grey line), $\mathbf{w}$ (grey line), and $\mathbf{m}$ (dotted grey line), and the set $I$ (thin black line with dots and diamonds on its boundary).}\label{fig_u_lower_u_upper_w_m_gen_set}
				\end{figure}
				
			Our approach for this proof is as follows. By Theorem \ref{thm_quartic_CCGPSR_curves_classification} we know that up to equivalence, any maximal connected incomplete quartic GPSR curve $\mathcal{H}$ in standard form corresponds to a defining polynomial of the form $h=x^4-x^2y^2+Lxy^3+Ky^4$ with
				\begin{align}
					\left(\begin{smallmatrix}L\\K\end{smallmatrix}\right)\in
						&\left\{L>\mathbf{w}(K),\ K<-\tfrac{1}{12}\right\}
						\cup\left\{L>\tfrac{2\sqrt{2}}{3\sqrt{3}},\ K=-\tfrac{1}{12}\right\}\notag\\
						\cup&\left\{L>\mathbf{u}_{\text{lower}}(K),\ -\tfrac{1}{12}<K\leq\tfrac{1}{4}\right\}
						\cup\left\{L>\mathbf{u}_{\text{upper}}(K),\ \tfrac{1}{4}<K\right\}\notag\\
						\cup&\left\{0\leq L<\mathbf{u}_{\text{upper}}(K),\ \tfrac{1}{4}<K\right\},\label{eqn_set_LK_incomplete}
				\end{align}
			see Figure \ref{fig_incomplete_area}.
				\begin{figure}[H]
					\centering
					\includegraphics[scale=0.2]{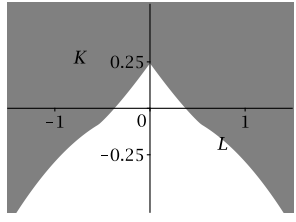}
					\caption{The set in equation \eqref{eqn_set_LK_incomplete} marked in grey.}\label{fig_incomplete_area}
				\end{figure}
			\noindent
			Let $\mathcal{H}_{L,K}$ denote the maximal connected incomplete quartic GPSR curve corresponding to $(L,K)^T$ in the above set \eqref{eqn_set_LK_incomplete}. First we will show that every $\mathcal{H}_{L,K}$ with
				\begin{equation*}
					\left(\begin{smallmatrix}L\\K\end{smallmatrix}\right)\in\left\{0\leq L<\mathbf{u}_{\text{upper}}(K),\ \tfrac{1}{4}<K\right\}
				\end{equation*}
			is equivalent to $\mathcal{H}_{1,\widetilde{K}}$ for some $\widetilde{K}>\frac{1}{4}$. Then we will show that every $\mathcal{H}_{L,K}$ with
				\begin{equation*}
					\left(\begin{smallmatrix}L\\K\end{smallmatrix}\right)\in\left\{L>\mathbf{w}(K),\ K<-\tfrac{1}{12}\right\}
				\end{equation*}
			is equivalent to $\mathcal{H}_{1,\widetilde{K}}$ for some $\mathbf{w}^{-1}(1)=-\frac{25}{72}<\widetilde{K}<-\frac{1}{12}$. Lastly we will show that every $\mathcal{H}_{L,K}$ with
				\begin{equation*}
					\left(\begin{smallmatrix}L\\K\end{smallmatrix}\right)\in\left\{L>\mathbf{u}_{\text{lower}}(K),\ -\tfrac{1}{12}<K\leq\tfrac{1}{4}\right\}\cup\left\{L>\mathbf{u}_{\text{upper}}(K),\ \tfrac{1}{4}<K\right\}
				\end{equation*}
			is equivalent to $\mathcal{H}_{1,\widetilde{K}}$ for some $\widetilde{K}\in\left(-\frac{1}{12},U\right)$. That, together with the information that the graphs of $\mathbf{u}_{\text{upper}}$, $\mathbf{u}_{\text{lower}}$, $\mathbf{w}$, and $\boldsymbol{\lambda}$ are maximal integral curves of $\mathcal{V}$, shows that every maximal connected incomplete quartic GPSR curve is equivalent to at least one of Thm. \ref{thm_quartic_CCGPSR_curves_classification} \hyperref[eqn_incomplete_qGPSRcurves_class_a]{a)}--\hyperref[eqn_incomplete_qGPSRcurves_class_d]{d)}. We then need to show that they are pairwise inequivalent, describe their respective connected components of $\{h=1\}\cap\{\text{hyperbolic points of }h\}$, and determine the respective automorphism groups $G^h$.
			
			Let $(L,K)^T\in\left\{0\leq L<\mathbf{u}_{\text{upper}}(K),\ \tfrac{1}{4}<K\right\}$. In order to show the existence of $\widetilde{K}>\frac{1}{4}$, such that $\mathcal{H}_{L,K}$ is equivalent to $\mathcal{H}_{0,\widetilde{K}}$, we will proceed as we did in the proof of Theorem \ref{thm_quartic_CCGPSR_curves_classification} for the set in \eqref{eqn_set_upper_bell}. If $L=0$, there is nothing to show. Consider the diffeomorphism
				\begin{align}
					\widetilde{F}&:\left\{0< L<\mathbf{u}_{\text{upper}}(K),\ \tfrac{1}{4}<K\right\}\to \left\{0<\varphi<1,\ \frac{1}{4}<\psi\right\}\notag\\
					\widetilde{F}&:\left(\begin{matrix}L\\ K\end{matrix}\right)\mapsto \left(\begin{matrix}\frac{L}{\mathbf{u}_{\mathrm{upper}}(K)}\\ K\end{matrix}\right),\label{eqn_above_bell_trafo}
				\end{align}
			see Figure \ref{fig_above_bell_trafo}.
				\begin{figure}[H]
					\centering
					\begin{subfigure}[h]{0.3\linewidth}
						\includegraphics[scale=0.2]{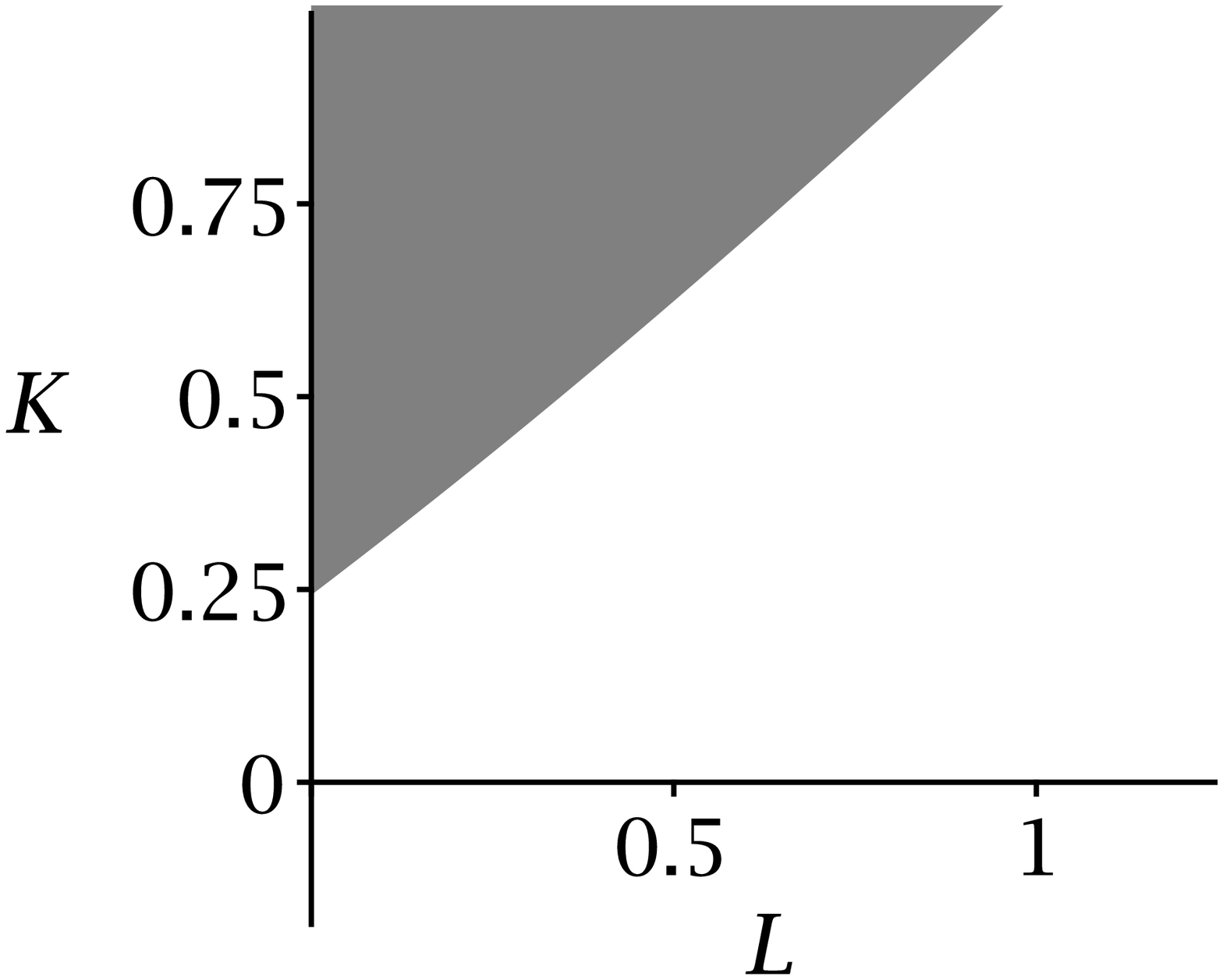}
					\end{subfigure}
					\begin{subfigure}[h]{0.3\linewidth}
						\includegraphics[scale=0.2]{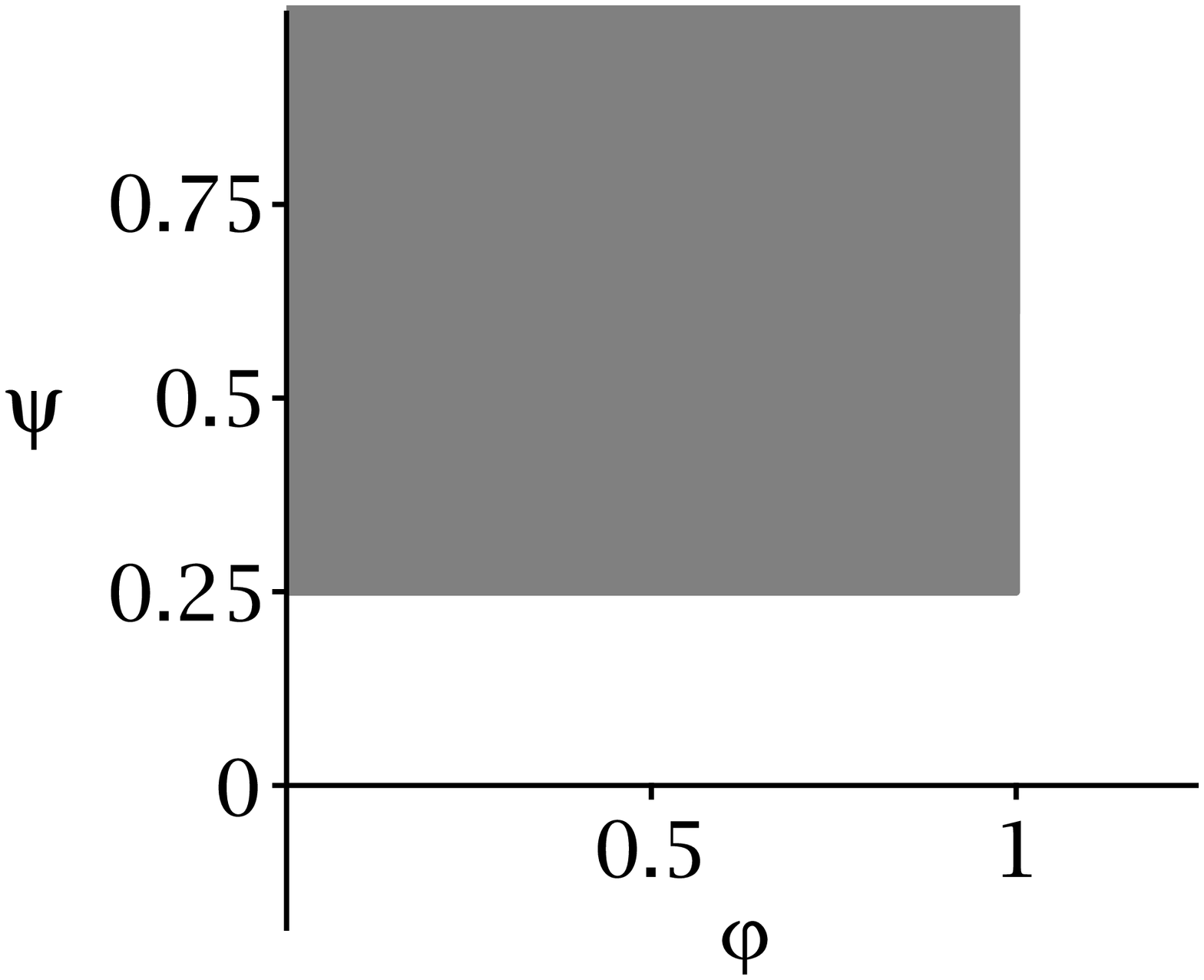}
					\end{subfigure}
					\caption{Domain and codomain of $\widetilde{F}$ \eqref{eqn_above_bell_trafo}.}\label{fig_above_bell_trafo}
				\end{figure}
			\noindent
			Symbolically, $\widetilde{F}$ fulfils \eqref{eqn_Ftilda_differential} and \eqref{eqn_Ftilda_inverse} with $\mathbf{u}$ replaced with $\mathbf{u}_{\mathrm{upper}}$. For $\widetilde{F}_*\mathcal{V}$ we obtain the same formula as in \eqref{eqn_Ftilda_V_pushforward}, now with domain $0<\varphi<1$ and $\frac{1}{4}<\psi$, on which we also define $\widetilde{\mathcal{R}}$ as in \eqref{eqn_Rtilda_def}. Note that on our current domain for $\varphi$ and $\psi$, $\widetilde{\mathcal{R}}$ is negative. We further define $\mathcal{Y}$ on $\left\{0<\varphi<1,\ \frac{1}{4}<\psi\right\}$ as in \eqref{eqn_Y_def} and $J(t)$ as in \eqref{eqn_J_def}. As in \eqref{eqn_gamma_t_in_Ftildastuff} we find that $\gamma(t)=\left(\begin{smallmatrix}\varphi(t)\\ t\end{smallmatrix}\right)$ is an integral curve of $\mathcal{Y}$ if and only if
				\begin{equation*}
					\varphi(t)=\sqrt{1-c \exp\left(-2\int\limits_{t_0}^tJ(s) \D s\right)}
				\end{equation*}
			for some $t_0>\frac{1}{4}$ with $c>0$ chosen so that the initial condition $(\varphi(t_0),t_0)^T$ for $\gamma$ is fulfilled. In order to show the existence of $\widetilde{K}>0$, it suffices to show that for all $c>0$ and all $t_0>\frac{1}{4}$ there exists $t\in\left(\frac{1}{4},t_0\right)$, such that $\varphi(t)=0$. In order to do so we proceed analogously to \eqref{eqn_J_int_infinity} by showing that
				\begin{equation*}
					\lim\limits_{t\searrow1/4}\int\limits_{t_0}^tJ(s)\D s=-\infty.
				\end{equation*}
			Without loss of generality we can choose $t_0=1$. Similar to \eqref{eqn_J_int_1stequiv_infinity}, the above statement is equivalent to
				\begin{equation}
					\lim\limits_{t\searrow1/4}\int\limits_{t}^1\frac{1-4s}{1-36s+\sqrt{(1+12s)^3}}\D s=-\infty,\label{eqn_s_limit_incomplete_case}
				\end{equation}
			with sign and integration range changed to fit our present case. We now substitute $s=r+\frac{1}{4}$ and obtain that \eqref{eqn_s_limit_incomplete_case} is equivalent to
				\begin{equation*}
					\lim\limits_{t\searrow 0}\int\limits_{t}^{\frac{3}{4}}\frac{r}{-9r-2+(6r+2)\sqrt{1+3r}}\D r=\infty.
				\end{equation*}
			Note that $-9r-2+(6r+2)\sqrt{1+3r}>0$ for all $r>0$. As in \eqref{eqn_J_int_3rdequiv_infinity} we show that there exists $\varepsilon>0$ and $A>0$, such that
				\begin{equation*}
					9r-2+(-6r+2)\sqrt{1-3r}\leq Ar^2
				\end{equation*}
			for all $r\in (0,\varepsilon)$. This follows from
				\begin{equation*}
					\lim\limits_{r\searrow 0} \frac{-9r-2+(6r+2)\sqrt{1+3r}}{r^2}=\frac{27}{4}.
				\end{equation*}
			We have thus shown, independently of $c>0$ and $t_0>\frac{1}{4}$, the existence of $t\in\left(\frac{1}{4},t_0\right)$, such that $\varphi(t)=0$, and thereby the existence of $\widetilde{K}$ having the required properties so that $\mathcal{H}_{L,K}$ is equivalent to $\mathcal{H}_{0,\widetilde{K}}$.
			
			Next, consider $\mathcal{H}_{L,K}$ with $(L,K)^T\in\left\{L>\mathbf{w}(K),\ K<-\tfrac{1}{12}\right\}$. In order to show that $\mathcal{H}_{L,K}$ is equivalent to $\mathcal{H}_{1,\widetilde{K}}$ for some $-\frac{25}{72}<\widetilde{K}<-\frac{1}{12}$, we will use the following approach. Let $X\in\mathfrak{X}((0,\infty)\times\mathbb{R})$ be a nowhere vanishing vector field, such that up to orientation preserving reparametrisation, $t\mapsto (t,0)^T$ and $t\mapsto (t,1)^T$ for $t>0$ are integral curves of $X$. Further assume that $\D x(X)\ne 0$ in $(0,\infty)\times(0,1)$, where $(x,y)$ denote the canonical coordinates on $(0,\infty)\times\mathbb{R}\subset\mathbb{R}^2$. Then every maximal integral curve $\gamma$ of $X$ with initial condition $\gamma(0)\in(0,\infty)\times(0,1)$ intersects the set $\{1\}\times(0,1)$ in precisely one point, which justifies calling it a generating set for the maximal integral curves of $X$ restricted to the strip $(0,\infty)\times(0,1)$. To see that this claim is true, observe that $\D x(X)\ne 0$ in $(0,\infty)\times(0,1)$ and $t\mapsto (t,0)^T$, $t\mapsto (t,1)^T$ being integral curves up to orientation preserving reparametrisation implies $\D x(X)>0$ in $(0,\infty)\times[0,1]$.
			The closure of the graph of any maximal integral curve $\gamma$ of $X$ with $\gamma(0)\in(0,\infty)\times(0,1)$ in $(0,\infty)\times\mathbb{R}$ cannot contain a point in $(0,\infty)\times\{0,1\}$, since that would by the smoothness of $X$ mean that $\gamma$ would, by the Picard–Lindel\"of theorem, coincide with either $t\mapsto (t,0)^T$ or $t\mapsto (t,1)^T$ up to reparametrisation, which contradicts the choice of its initial condition. Hence, every such integral curve meets $\{1\}\times[0,1]$ in either finite positive or finite negative time, which can be seen by minimising $\D x(X)$ over an appropriate compact subset of $(0,\infty)\times[0,1]$ containing the initial value $\gamma(0)$ and $\{1\}\times(0,1)$. Furthermore, any such integral curve cannot meet $\{1\}\times(0,1)$ more than once since this would imply the existence of some $t$ in the domain of $\gamma$, such that $\D x(\dot\gamma)=0$, which contradicts our assumptions on $X$.
			
			We can apply this method to our current case. Up to a fittingly chosen diffeomorphism of $\left\{L>\frac{2\sqrt{2}}{3\sqrt{3}}\right\}$,
				\begin{equation}
					\left(\begin{matrix}L\\K\end{matrix}\right)\mapsto
					\left(\begin{matrix}L-\frac{2\sqrt{2}}{3\sqrt{3}}\\
						\frac{K+\frac{1}{12}}{-\mathbf{w}^{-1}(L)+\frac{1}{12}}+1
					\end{matrix}\right),\label{eqn_strip_trafo}
				\end{equation}
			mapping the graph of $\mathbf{m}$ to $(0,\infty)\times\{1\}$ and the graph of $\mathbf{w}$ to $(0,\infty)\times\{0\}$, and the interior to $(0,\infty)\times(0,1)$, see Figure \ref{fig_strip_trafo}, we are precisely in the described setting.
				\begin{figure}[H]
					\centering
					\begin{subfigure}[h]{0.3\linewidth}
						\includegraphics[scale=0.2]{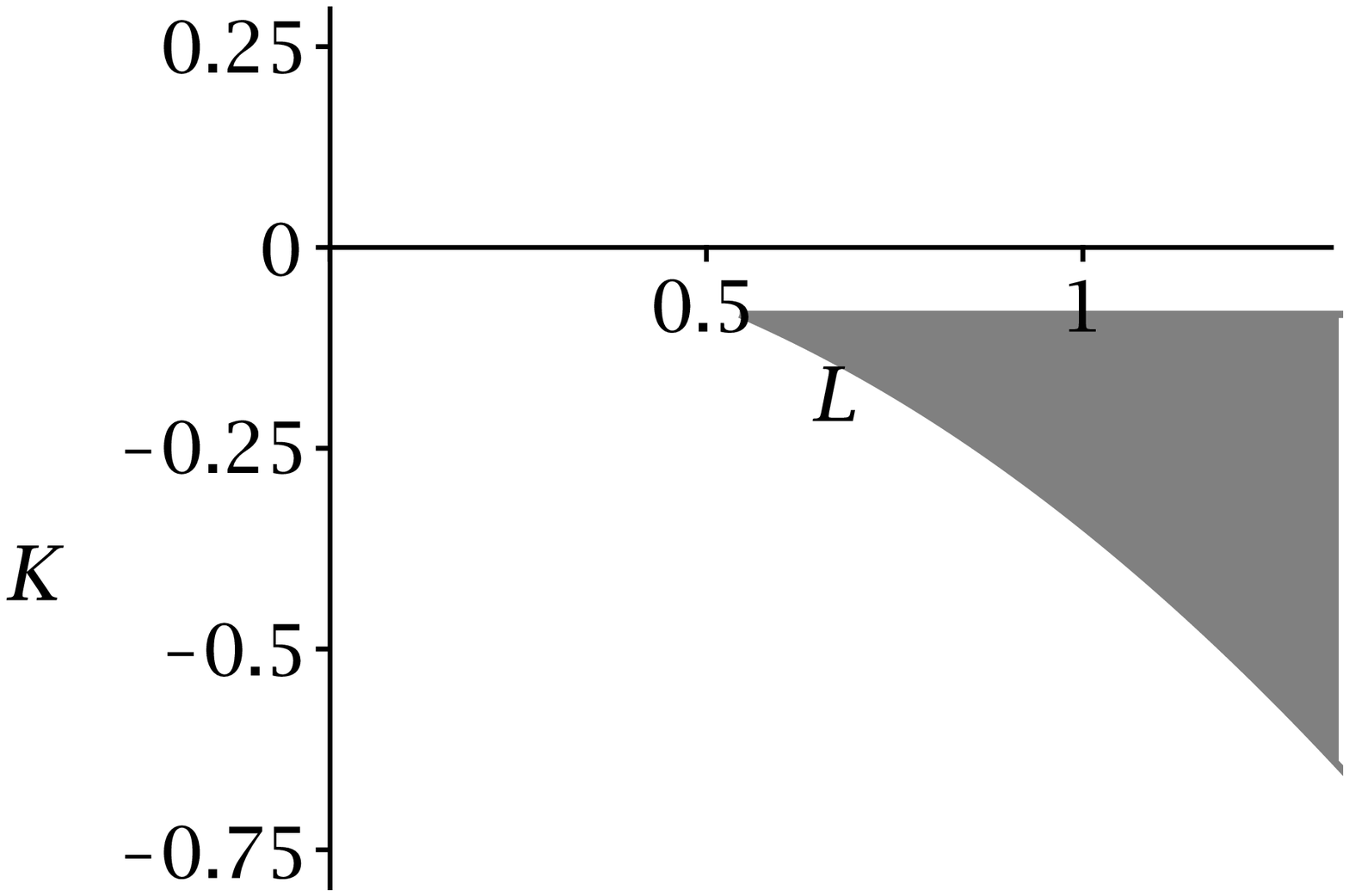}
					\end{subfigure}
					\begin{subfigure}[h]{0.3\linewidth}
						\includegraphics[scale=0.2]{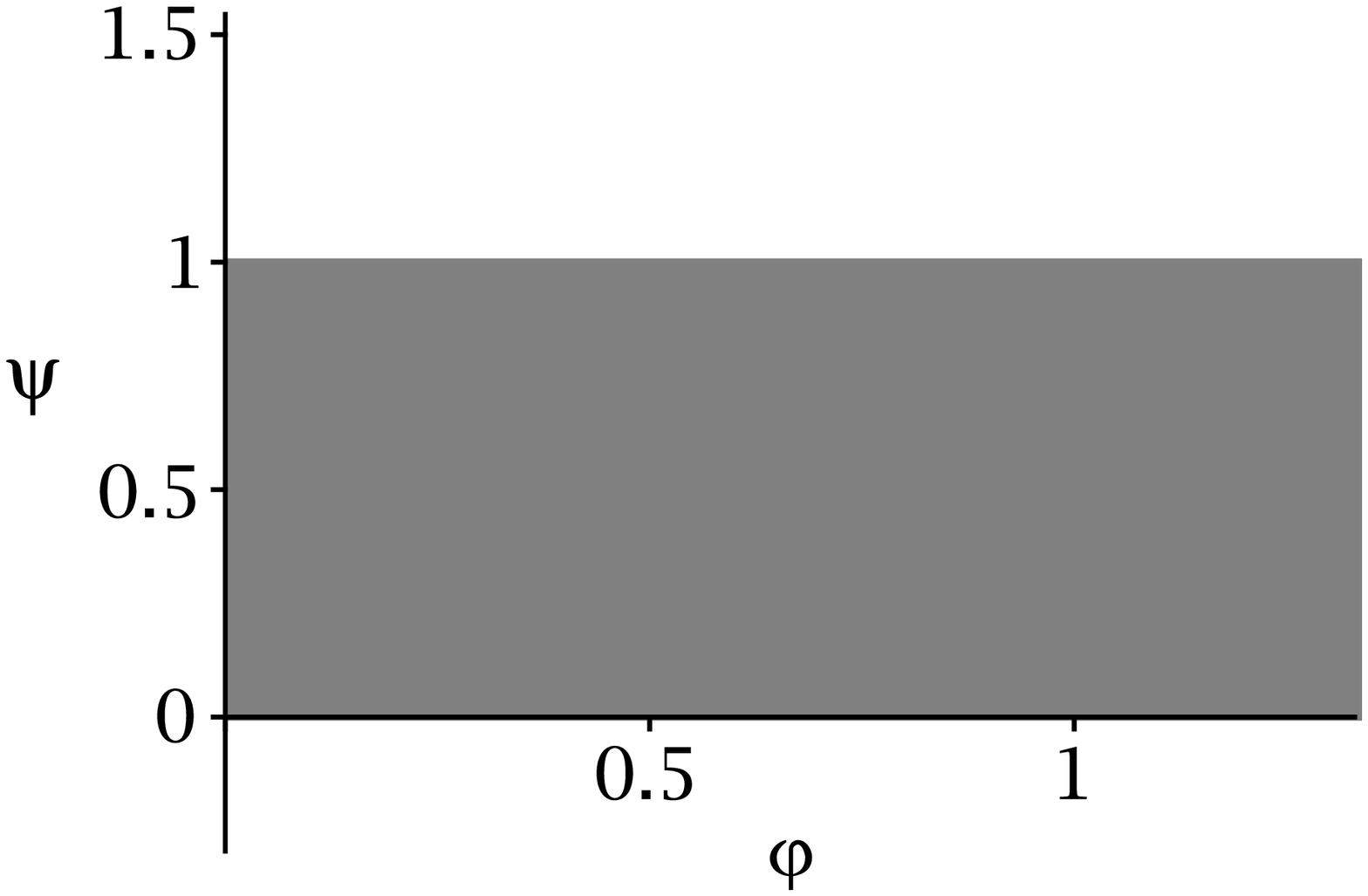}
					\end{subfigure}
					\caption{Domain and codomain of the map in \eqref{eqn_strip_trafo}.}\label{fig_strip_trafo}
				\end{figure}
			\noindent
			Since the above diffeomorphism is, up to a translation, a rescaling in the $K$-axis, it suffices to show that in $\left\{L>\frac{2\sqrt{2}}{3\sqrt{3}},\ -\mathbf{w}^{-1}(L)<K<-\frac{1}{12}\right\}$, $\D L(\mathcal{V})\ne 0$ in order to show that the set $\left\{L=1,\ -\frac{25}{72}<K<-\frac{1}{12}\right\}$ is a generating set for all maximal incomplete quartic GPSR curves $\mathcal{H}_{L,K}$ with $(L,K)^T\in \left\{L>\mathbf{w}(K),\ K<-\tfrac{1}{12}\right\}$. We find $\D L(\mathcal{V})=0$ if and only if $K=\frac{1}{4}-\frac{9}{8}L^2$, 
				\begin{equation*}
					\frac{1}{4}-\frac{9}{8}L^2=\mathbf{w}^{-1}(L)\quad \Leftrightarrow\quad L=\pm\frac{2\sqrt{2}}{3\sqrt{3}},
				\end{equation*}
			and $\left.\left(\frac{1}{4}-\frac{9}{8}L^2\right)\right|_{L=1}=-\frac{63}{72}<-\frac{25}{72}=\mathbf{w}^{-1}(1)$. Hence, $\D L(\mathcal{V})\ne 0$ in
				\begin{equation*}
					\left\{L>\mathbf{w}(K),\ K<-\tfrac{1}{12}\right\}=\left\{L>\frac{2\sqrt{2}}{3\sqrt{3}},\ -\mathbf{w}^{-1}(L)<K<-\frac{1}{12}\right\}
				\end{equation*}
			as needed. We conclude that the set $\left\{L=1,\ -\frac{25}{72}<K<-\frac{1}{12}\right\}$ is indeed a generating set as stated above.
			
			In the next step we need to deal with the set
				\begin{equation}
					\left\{L>\mathbf{u}_{\text{lower}}(K),\ -\tfrac{1}{12}<K\leq\tfrac{1}{4}\right\}\cup\left\{L>\mathbf{u}_{\text{upper}}(K),\ \tfrac{1}{4}<K\right\}.\label{eqn_set_homotopy}
				\end{equation}
			In order to show that $\left\{L=1,\ -\frac{1}{12}<K<U\right\}$ is a generating set for all maximal incomplete quartic GPSR curves $\mathcal{H}_{L,K}$ with $(L,K)^T$ contained in the above set we will construct a certain homotopy, filling the set \eqref{eqn_set_homotopy}, and having the property that each curve in that homotopy is transversal to $\mathcal{V}$. Consider the continuous homotopy of curves
				\begin{align}
					H:(0,\infty)\times(0,1)&\to\left\{L>\mathbf{u}_{\text{lower}}(K),\ -\tfrac{1}{12}<K\leq\tfrac{1}{4}\right\}\cup\left\{L>\mathbf{u}_{\text{upper}}(K),\ \tfrac{1}{4}<K\right\},\label{eqn_H_homotopy}\\
					(s,t)&\mapsto\left\{
						\begin{array}{ll}
							\left(\begin{matrix}
								\frac{2\sqrt{2}}{3\sqrt{3}}-\sqrt{6}(1-s)r(s)t\\
								r(s)t-\frac{1}{12}
							\end{matrix}\right), & s\in(0,1],\\
							\left(\begin{matrix}
								\frac{2\sqrt{2}}{3\sqrt{3}}+s-1\\
								r(s)t-\frac{1}{12}
							\end{matrix}\right), & s\in(1,\infty),
						\end{array}
					\right.\notag
				\end{align}
			where the continuous function $r:(0,\infty)\to\mathbb{R}_{>0}$ is defined by
				\begin{equation*}
					\left\{
						\begin{array}{ll}
							\mathbf{u}_{\mathrm{lower}}\left(r(s)-\frac{1}{12}\right) = \frac{2\sqrt{2}}{3\sqrt{3}}-\sqrt{6}(1-s)r(s), & s\in\left(0,\frac{1}{3}\right),\\
							r\left(\frac{1}{3}\right)=\frac{1}{3}, & \\
							\mathbf{u}_{\mathrm{upper}}\left(r(s)-\frac{1}{12}\right) = \frac{2\sqrt{2}}{3\sqrt{3}}-\sqrt{6}(1-s)r(s), & s\in\left(\frac{1}{3},1\right],\\
							\mathbf{u}_{\mathrm{upper}}\left(r(s)-\frac{1}{12}\right)=\frac{2\sqrt{2}}{3\sqrt{3}}+s-1, & s>1.
						\end{array}
					\right.
				\end{equation*}
			See Figure \ref{fig_homotopy} for plots of $r$ and the homotopy $H$.
				\begin{figure}[H]
					\centering
					\begin{subfigure}[h]{0.3\linewidth}
						\includegraphics[scale=0.2]{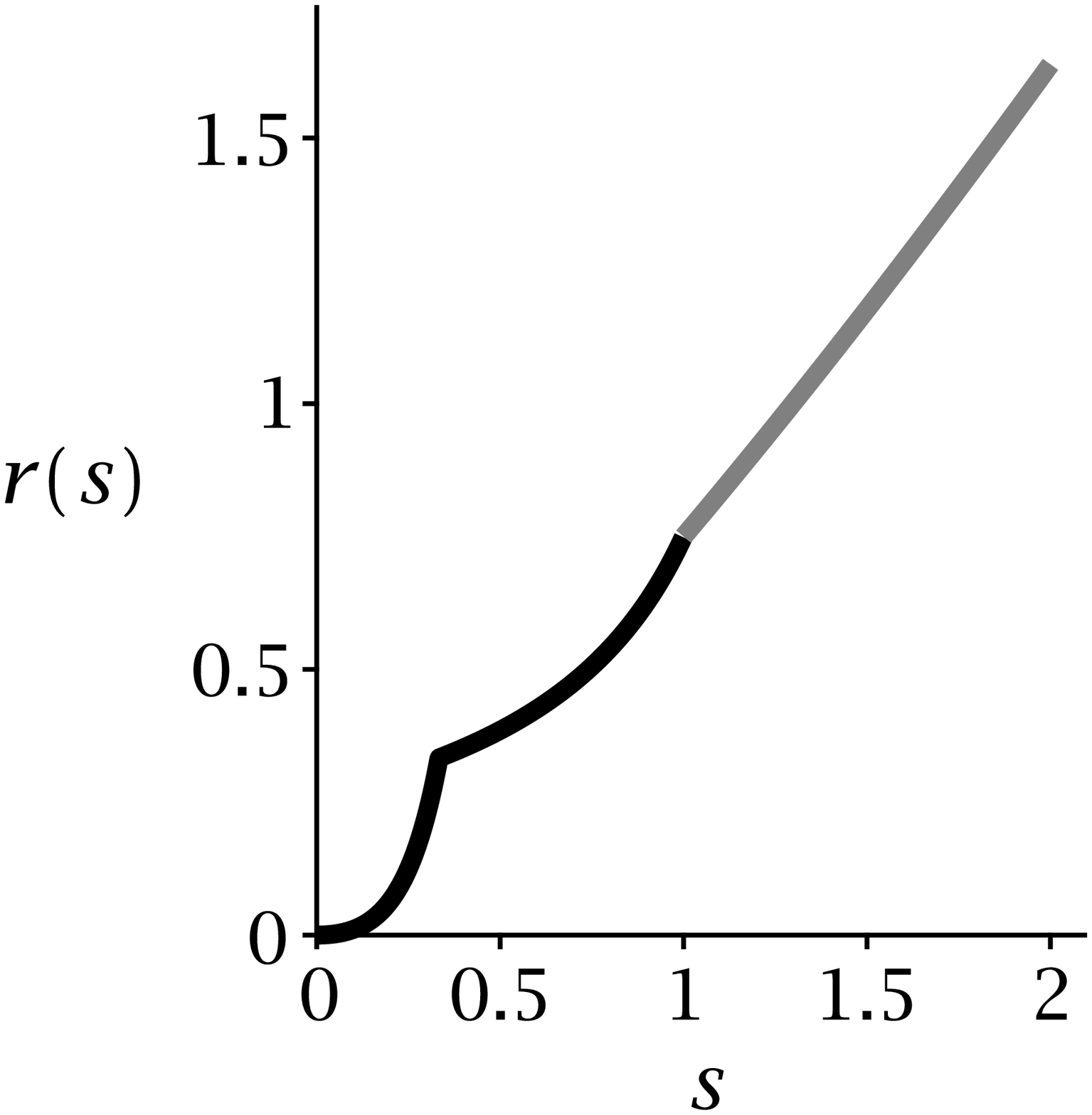}
					\end{subfigure}
					\begin{subfigure}[h]{0.3\linewidth}
						\includegraphics[scale=0.2]{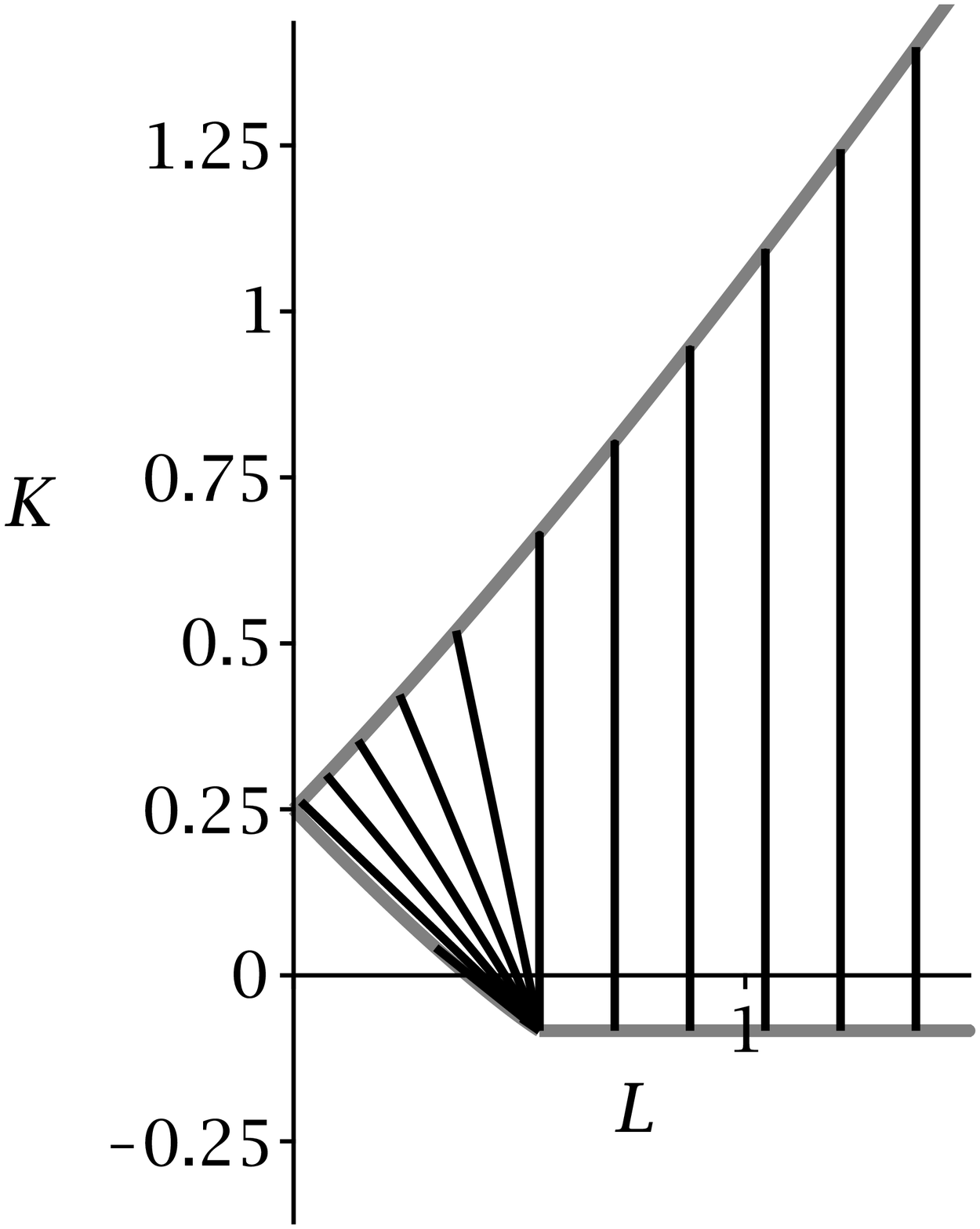}
					\end{subfigure}
					\caption{A plot of $r$ coloured black for $s\in(0,1]$ and grey for $s>1$, and graphs of $H(s,\cdot)$ for several values $s$.}\label{fig_homotopy}
				\end{figure}
			\noindent
			Note that $H$ is, in fact, a bijection with the property
				\begin{equation*}
					\mathrm{Im}\left(H\left(2-\frac{2\sqrt{2}}{3\sqrt{3}},\cdot\right)\right)=\left\{L=1,\ -\frac{1}{12}<K<U\right\}.
				\end{equation*}
			In order to see that, recall that $\mathbf{u}_{\mathrm{lower}}$ is a convex function. For every fixed $s\in(0,\infty)$, the curve $H(s,\cdot)$ is transversal to $\mathcal{V}$. In order to see this we will first assume $s\geq 1$. In that case we need to show $\D L(\mathcal{V})\ne 0$ along $H(s,\cdot)$. We have already seen that $\D L(\mathcal{V})\ne 0$ if and only if $K=\frac{1}{4}-\frac{9}{8}L^2$, and by $\frac{1}{4}-\frac{9}{8}L^2\leq -\frac{1}{12}$ and \eqref{eqn_set_homotopy} we obtain that $H(s,\cdot)$ is indeed transversal to $\mathcal{V}$ for $s\geq 1$ as claimed. For $0<s<1$ fixed, $H(s,\cdot)$ being transversal to $\mathcal{V}$ is equivalent to
				\begin{equation}
					\left\langle \mathcal{V}_{H(s,t)},\left(\begin{smallmatrix}
						r(s)\\ \sqrt{6}(1-s)r(s)
					\end{smallmatrix}\right)\right\rangle=0\label{eqn_H_transversal_s01}
				\end{equation}
			having no solution in $t$ with $t\in(0,1)$. In order to study \eqref{eqn_H_transversal_s01} we drop the overall positive scale $r(s)$, substitute $t$ with $\frac{1}{r(s)}\theta$, and obtain the equation
				\begin{equation*}
					-9\theta\left((1-s)^2\theta - \frac{4}{9}s\right)=0.
				\end{equation*}
			Solving for $\theta$, we obtain the set of solutions $\theta=0$ and $\theta=\frac{4s}{9(1-s)^2}$. The first one is to be expected by the fact that $\mathcal{V}$ vanishes at $(L,K)^T=\left(\frac{2\sqrt{2}}{3\sqrt{3}},-\frac{1}{12}\right)^T$. In order to proceed we first observe that $\frac{2\sqrt{2}}{3\sqrt{3}}-\sqrt{6}(1-s)\theta$ evaluated at $\theta=\frac{4s}{9(1-s)^2}$ is given by $\frac{2\sqrt{2}}{3\sqrt{3}}\cdot\frac{1-3s}{1-s}$, meaning that it is non-positive for all $s\in\left[\frac{1}{3},1\right)$. This already shows that $H(s,\cdot)$ is transversal to $\mathcal{V}$ for all $s\in\left[\frac{1}{3},1\right)$. For $s\in\left(0,\frac{1}{3}\right)$ we need to compare the graph of the curve
				\begin{equation}
					s\mapsto H\left(\frac{1}{r(s)}\cdot\frac{4s}{9(1-s)^2},s\right)\label{eqn_H_graph_check}
				\end{equation}
			with the graph of $\mathbf{u}_{\mathrm{lower}}$, cf. Figure \ref{fig_transversality_u}.
				\begin{figure}[H]
					\centering
					\includegraphics[scale=0.2]{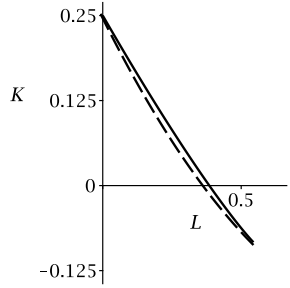}
					\caption{Graphs of $\mathbf{u}_{\mathrm{lower}}$ (solid line) and the map in \eqref{eqn_H_graph_check} (dashed line).}\label{fig_transversality_u}
				\end{figure}
			\noindent
			To do so we first verify that we can write the graph of \eqref{eqn_H_graph_check} over $s\in\left(0,\frac{1}{3}\right)$ as the graph of the function
				\begin{equation*}
					K\mapsto F(K):=\frac{\sqrt{6}\left(-5-12K+2\sqrt{7+36K}\right)}{-6+3\sqrt{7+36K}},\quad K\in\left(-\frac{1}{12},\frac{1}{4}\right).
				\end{equation*}
			In order to show the transversality of the curves $H(s,\cdot)$ for $s\in\left(0,\frac{1}{3}\right)$ it now suffices to show that $F(K)<\mathbf{u}_{\mathrm{lower}}(K)$ for all $K\in\left(-\frac{1}{12},\frac{1}{4}\right)$. We find that
				\begin{equation*}
					\lim\limits_{K\nearrow 1/4}(F(K)-\mathbf{u}_{\mathrm{lower}}(K))=\lim\limits_{K\searrow -1/12}(F(K)-\mathbf{u}_{\mathrm{lower}}(K))=0.
				\end{equation*}
			Now we verify that $F(0)-\mathbf{u}_{\mathrm{lower}}(0))<0$. Preferably with the help of computer algebra software, one can now quickly check that $F(K)=\mathbf{u}_{\mathrm{lower}}(K)$ has indeed no further solutions in $\left(-\frac{1}{12},\frac{1}{4}\right)$, proving that each curve $H(s,\cdot)$ for $s\in\left(0,\frac{1}{3}\right)$ is in fact transversal to $\mathcal{V}$ as claimed. For the next step we use the transversality property of $H$ with respect to $\mathcal{V}$ to show that every maximal integral curve $\gamma$ of $\mathcal{V}$ with initial condition $\gamma(0)$ contained in the set \eqref{eqn_set_homotopy} meets $\mathrm{Im}\left(H\left(2-\frac{2\sqrt{2}}{3\sqrt{3}},\cdot\right)\right)$ precisely once.
			
			Firstly, suppose that $\gamma(0)\notin\mathrm{Im}\left(H\left(2-\frac{2\sqrt{2}}{3\sqrt{3}},\cdot\right)\right)$ and that $\mathrm{Im}(\gamma)\cap\mathrm{Im}\left(H\left(2-\frac{2\sqrt{2}}{3\sqrt{3}},\cdot\right)\right)=\emptyset$. For $L(\gamma(0))>\frac{2\sqrt{2}}{3\sqrt{3}}$ the argument is analogous to the one used in our study of the set $\left\{L>\mathbf{w}(K),\ K<-\tfrac{1}{12}\right\}$ and we obtain a contradiction. Next suppose that $L(\gamma(0))<\frac{2\sqrt{2}}{3\sqrt{3}}$. Let $s_0\in\left(0,2-\frac{2\sqrt{2}}{3\sqrt{3}}\right)$ be the unique real number, such that $\mathrm{Im}(\gamma)\cap\mathrm{Im}(H(s_0,\cdot))\ne\emptyset$, but for all $s\in\left(s_0,2-\frac{2\sqrt{2}}{3\sqrt{3}}\right)$, $\mathrm{Im}(\gamma)\cap\mathrm{Im}(H(s,\cdot))=\emptyset$. The existence of such $s_0$ can be seen as follows. $H$ being a bijection implies that there exists some $\widetilde{s}\in\left(0,2-\frac{2\sqrt{2}}{3\sqrt{3}}\right)$, such that $\mathrm{Im}(\gamma)\cap\mathrm{Im}(H(\widetilde{s},\cdot))\ne\emptyset$. Suppose that this holds true for all $\widetilde{s}\in\left(0,2-\frac{2\sqrt{2}}{3\sqrt{3}}\right)$. Then by the assumption that $\mathrm{Im}(\gamma)\cap\mathrm{Im}\left(H\left(2-\frac{2\sqrt{2}}{3\sqrt{3}},\cdot\right)\right)=\emptyset$, it follows with $\mathrm{Im}(\gamma)\cap\mathrm{Im}(H(\widetilde{s},\cdot))\ne\emptyset$ for all $\widetilde{s}\in\left(0,2-\frac{2\sqrt{2}}{3\sqrt{3}}\right)$ that we can find a convergent sequence $\{p_n\}_{n\in\mathbb{N}}\subset\mathrm{Im}(\gamma)$, such that $p:=\lim\limits_{n\to\infty}p_n\in\partial\,\mathrm{Im}\left(H\left(2-\frac{2\sqrt{2}}{3\sqrt{3}},\cdot\right)\right)$. By the continuity of $\mathcal{V}$, we can extend $\gamma$ continuously to $p$, but $p$ is an element of the graph of either $\mathbf{u}_{\mathrm{upper}}$ or $\mathbf{m}$, meaning that no maximal integral curve of $\mathcal{V}$ with initial condition in the set \eqref{eqn_set_homotopy} can have $p$ in its image. This is a contradiction to the assumption of $\gamma(0)$. Thus,
				\begin{equation}
					s_0:=\sup\limits_{\left\{0<\widetilde{s}<1,\ \mathrm{Im}(\gamma)\cap\mathrm{Im}(H(\widetilde{s},\cdot))\ne\emptyset\right\}} \widetilde{s}\label{eqn_s0_def}
				\end{equation}
			has the property $s_0<1$. It remains to show that $\mathrm{Im}(\gamma)\cap\mathrm{Im}(H(s_0,\cdot))\ne \emptyset$. To do so suppose that the latter is not true. Similarly to the sequence construction above, one proves that there exists no convergence sequence in $\mathrm{Im}(\gamma)$ with limit in the graph of either $\mathbf{u}_{\mathrm{upper}}$ or $\mathbf{m}$. From that it follows for any fixed $s\in\left(0,2-\frac{2\sqrt{2}}{3\sqrt{3}}\right)$ with $\mathrm{Im}(\gamma)\cap\mathrm{Im}(H(s,\cdot))=\emptyset$ that the Euclidean distance of the two compact sets $\overline{\mathrm{Im}(\gamma)}$ and $\overline{\mathrm{Im}(H(s,\cdot))}$ is positive, i.e.
				\begin{equation*}
					\D_{\,\mathrm{Eucl.}}\left(\overline{\mathrm{Im}(\gamma)},\overline{\mathrm{Im}(H(s,\cdot))}\right)=\varepsilon>0.
				\end{equation*}
			But since by construction $H$ can be continuously extended to the larger domain $(0,\infty)\times[0,1]$, we obtain that there exists $\delta>0$, such that for all $\widetilde{s}$ fulfilling $|s-\widetilde{s}|<\delta$
				\begin{equation*}
					\D_{\,\mathrm{Eucl.}}\left(\overline{\mathrm{Im}(\gamma)},\overline{\mathrm{Im}(H(\widetilde{s},\cdot))}\right)>\frac{\varepsilon}{2}.
				\end{equation*}
			This would under our current supposition in particular hold for $s=s_0$, which would be a contradiction to its definition in \eqref{eqn_s0_def}. Hence, $\mathrm{Im}(\gamma)\cap\mathrm{Im}(H(s_0,\cdot))\ne \emptyset$ as claimed. Now let $\widetilde{t}\in\gamma^{-1}\left(\mathrm{Im}(\gamma)\cap\mathrm{Im}(H(s_0,\cdot))\right)$. Note that $\widetilde{t}$ might in general not be unique, see Figure \ref{fig_gamma_nonunique}.
				\begin{figure}[H]
					\centering
					\includegraphics[scale=0.2]{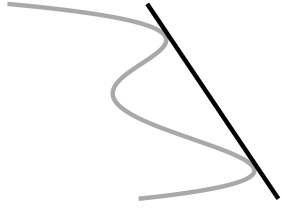}
					\caption{Example of $\gamma^{-1}\left(\mathrm{Im}(\gamma)\cap\mathrm{Im}(H(s_0,\cdot))\right)$ consisting of two elements.}\label{fig_gamma_nonunique}
				\end{figure}
			\noindent
			Then
				\begin{equation*}
					\dot\gamma(\widetilde{t})=R\left.\frac{\partial}{\partial t}\right|_{t=\widetilde{t}}H(s_0,t)
				\end{equation*}
			for some $R\in\mathbb{R}\setminus{0}$. Note that $R\ne 0$ is a consequence of $\mathcal{V}$ not vanishing at any point in the set \eqref{eqn_set_homotopy}, and the $\gamma$ being tangential to $H(s_0,\cdot)$ follows since otherwise $s_0$ would not be maximal in the sense of \eqref{eqn_s0_def}. But this is a contradiction to the transversality property of $H$, i.e. that for all $s>0$, $\mathcal{H}(s,\cdot)$ is transversal to $\mathcal{V}$, but $\gamma$ is an integral curve of $\mathcal{V}$.

			Together with the trivial case $L(\gamma(0))=\frac{2\sqrt{2}}{3\sqrt{3}}$ we have thus proven that for all possible initial conditions $\gamma(0)$ in the set \eqref{eqn_set_homotopy}, the image of the corresponding maximal integral curve $\gamma$ of $\mathcal{V}$ has at least one common element with the set $\mathrm{Im}\left(H\left(2-\frac{2\sqrt{2}}{3\sqrt{3}},\cdot\right)\right)=\left\{L=\frac{2\sqrt{2}}{3\sqrt{3}},\ -\frac{1}{12}<K<U\right\}$. Next we need to show that it is precisely one common element, and that that element has precisely one preimage under $\gamma$.
			
			Suppose that $\left|\mathrm{Im}(\gamma)\cap \mathrm{Im}\left(H\left(2-\frac{2\sqrt{2}}{3\sqrt{3}},\cdot\right)\right)\right|\geq 2$. Choose $a,b$ with $a<b$ in the domain of $\gamma$, such that $\gamma(a)$ and $\gamma(b)$ are contained in $\mathrm{Im}\left(H\left(2-\frac{2\sqrt{2}}{3\sqrt{3}},\cdot\right)\right)$ and so that for all $t\in(a,b)$, $\gamma(t)\notin \mathrm{Im}\left(H\left(2-\frac{2\sqrt{2}}{3\sqrt{3}},\cdot\right)\right)$. Note that such a choice of $a,b$ is possible by the transversality of $H(s,\cdot)$ and $\mathcal{V}$, which excludes the possibility of any interval being mapped to $\mathrm{Im}\left(H\left(2-\frac{2\sqrt{2}}{3\sqrt{3}},\cdot\right)\right)$ under $\gamma$. Firstly, suppose that for all $t\in(a,b)$, $L(\gamma(t))>\frac{2\sqrt{2}}{3\sqrt{3}}$. In this case the arguments leading to a contradiction is, again, analogous to our study of the set $\left\{L>\mathbf{w}(K),\ K<-\tfrac{1}{12}\right\}$ and we will not repeat them here. Secondly, suppose that for all $t\in(a,b)$, $L(\gamma(t))<\frac{2\sqrt{2}}{3\sqrt{3}}$. Let $s_1\in(0,1)$ be the unique value, such that for all $s\in(0,s_1)$, $\gamma([a,b])\cap\mathrm{Im}(H(s,\cdot))=\emptyset$ and $\gamma([a,b])\cap\mathrm{Im}(H(s_1,\cdot))\ne\emptyset$. The existence and uniqueness of $s_1$ follows similarly as the existence and uniqueness of $s_0$, cf. the discussion for \eqref{eqn_s0_def}. Then for all $\widetilde{t}\in\gamma^{-1}\left(\gamma([a,b])\cap\mathrm{Im}(H(s_1,\cdot)\right)$,
				\begin{equation*}
					\dot\gamma(\widetilde{t})=R\left.\frac{\partial}{\partial t}\right|_{t=\widetilde{t}} H(s_1,t)
				\end{equation*}
			for some $R\in\mathbb{R}\setminus\{0\}$. This, again, contradicts the transversality property of $H(s_1,\cdot)$ to $\mathcal{V}$.
			
			Summarizing, we have shown that every maximal integral curve $\gamma$ of $\mathcal{V}$ with initial condition $\gamma(0)$ contained in the set \eqref{eqn_set_homotopy} meets $\left\{L=1,\ -\frac{1}{12}<K<U\right\}$ precisely once, turning the latter set into a generating set of the maximal integral curves of the restriction of $\mathcal{V}$ to \eqref{eqn_set_homotopy} as desired. We have thus proven that indeed every maximal connected incomplete quartic GPSR curve has a standard form given by one of Thm. \ref{thm_incomplete_quartics} \hyperref[eqn_incomplete_qGPSRcurves_class_a]{a)}--\hyperref[eqn_incomplete_qGPSRcurves_class_d]{d)}. It remains to prove the statements about the respective equivalent connected components and their respective automorphism groups, and show that the curves in Thm. \ref{thm_incomplete_quartics} \hyperref[eqn_incomplete_qGPSRcurves_class_a]{a)}--\hyperref[eqn_incomplete_qGPSRcurves_class_d]{d)} are pairwise inequivalent.
			
			We start with Thm. \ref{thm_incomplete_quartics} \hyperref[eqn_incomplete_qGPSRcurves_class_a]{a)}, that is $h=x^4-x^2y^2+Ky^4$ with $K>\frac{1}{4}$. In that case, $h=0$ has no real solutions, meaning that $\{h=1\}$ is compact. In order to see that, we check that $h|_{x=1}$ has a local maximum at $y=0$ with value $1$, and two local minima at $\pm\frac{1}{\sqrt{2K}}$ with value $\frac{4K-1}{4K}>0$. Since $h\left(\left(\begin{smallmatrix}0\\ 1\end{smallmatrix}\right)\right)=K>0$, it follows that $h$ is positive on $\mathbb{R}^2\setminus\{0\}$ and, hence, that $\{h=1\}$ is compact. This also implies that these curves are not singular at infinity. We now need to count and describe the connected components of $\{\text{hyperbolic points of }h\}$. Define
				\begin{equation*}
					f:=\det\left(\partial^2h\right)=12(-2x^4-(1-12K)x^2y^2-2Ky^4).
				\end{equation*}
			A point $p\in\{h>0\}$ is a hyperbolic point of $h$ if and only if $f(p)<0$. We check that $f\left(\left(\begin{smallmatrix}0\\ 1\end{smallmatrix}\right)\right)=-24K<0$ for all $K>\frac{1}{4}$. Now we need to study $f|_{x=1}$. We find
				\begin{equation*}
					f|_{x=1}=0\quad\Leftrightarrow\quad y=\pm_1\frac{\sqrt{-1+12K\pm_2\sqrt{1-40K+144K^2}}}{2\sqrt{K}}.
				\end{equation*}
			In the above equation, $\pm_1$ and $\pm_2$ are independent, so we have $4$ solutions. Also note that $-1+12K>0$, $1-40K+144K^2>0$, and $-1+12K>\sqrt{1-40K+144K^2}$ for all $K>\frac{1}{4}$, meaning that all of the above solutions are real and pairwise distinct for fixed $K$. We further calculate
				\begin{align*}
					&\partial_y(f|_{x=1})\left(\frac{\sqrt{-1+12K-\sqrt{1-40K+144K^2}}}{2\sqrt{K}}\right)\\
					&=\frac{12\sqrt{1-40K+144K^2}\sqrt{-1+12K-\sqrt{1-40K+144K^2}}}{\sqrt{K}}>0,\\
					&\partial_y(f|_{x=1})\left(\frac{\sqrt{-1+12K+\sqrt{1-40K+144K^2}}}{2\sqrt{K}}\right)\\
					&=\frac{12\sqrt{1-40K+144K^2}\sqrt{-1+12K+\sqrt{1-40K+144K^2}}}{\sqrt{K}}<0,
				\end{align*}
			for all $K>\frac{1}{4}$. Together with the symmetry of $f$ in $y$, this implies that $\{\text{hyperbolic points of }h\}$ has $4$ connected components. See Figure \ref{fig_hyp_points_K_1_2} for a plot of $\{h=1\}$ and the boundary of the hyperbolic points of $h$ for $K=\frac{1}{2}$.
				\begin{figure}[H]
					\centering
					\includegraphics[scale=0.2]{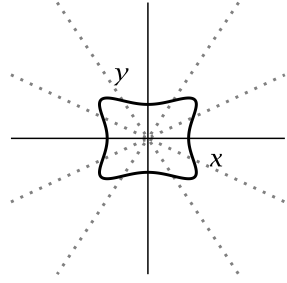}
					\caption{The set $\{h=1\}$ in black and the boundary of $\{\text{hyperbolic points of }h\}$ as dotted grey lines.}\label{fig_hyp_points_K_1_2}
				\end{figure}
			\noindent
			In order to see that all the connected components of $\{h=1\}\cap\{\text{hyperbolic points of }h\}$ are equivalent, it suffices by the homogeneity of degree $4$ of $h$ to show that the two connected components containing the points $\left(\begin{smallmatrix}1\\0\end{smallmatrix}\right)$ and $\left(\begin{smallmatrix}0\\\sqrt[4]{K}\end{smallmatrix}\right)$, respectively, are equivalent. We find that $A=\left(\begin{matrix} & -\sqrt[4]{K}\\ \frac{1}{\sqrt[4]{K}} & \end{matrix} \right)$, just as in the case $K\in\left(0,\frac{1}{4}\right)$ in Theorem \ref{thm_quartic_CCGPSR_curves_classification}, maps the first connected component to the latter. Thus, all connected components of $\{h=1\}\cap\{\text{hyperbolic points of }h\}$ are equivalent. Analogously to the case $K\in\left(0,\frac{1}{4}\right)$ in Theorem \ref{thm_quartic_CCGPSR_curves_classification} we obtain $G^h=\mathbb{Z}_2\ltimes\mathbb{Z}_2$.
			
			Next we consider Thm. \ref{thm_incomplete_quartics} \hyperref[eqn_incomplete_qGPSRcurves_class_b]{b)}, \hyperref[eqn_incomplete_qGPSRcurves_class_c]{c)}, and \hyperref[eqn_incomplete_qGPSRcurves_class_d]{d)}, that is $h=x^4-x^2y^2+xy^3+Ky^4=1$ with $K\in\left[-\frac{25}{72},U\right]$. In order to show that in each of these cases, $\{h=1\}\cap\{\text{hyperbolic points of }h\}$ has $4$ connected components and that they are pairwise equivalent, we will introduce new coordinates. Let $z_{\mathrm{max}}$ denote the biggest negative zero of $h|_{x=1,\ K=-25/72}$, and $z_{\mathrm{min}}$ the biggest negative zero of $h|_{x=1,\ K=U}$. Since $-\frac{25}{72}<U$, it follows that $z_{\mathrm{min}}<z_{\mathrm{max}}$. Explicitly, we have
				\begin{equation*}
					z_{\mathrm{min}}= \frac{2-\left(46+6\sqrt{57}\right)^{1/3}-4\left(46+6\sqrt{57}\right)^{-1/3}}{3} ,\quad z_{\mathrm{max}}=\frac{18 - 2\sqrt{6} - \sqrt{48 + 378\sqrt{6}}}{25}.
				\end{equation*}
			There is a one-to-one correspondents between $\left[-\frac{25}{72},U\right]$ and $\left[z_{\mathrm{min}},z_{\mathrm{max}}\right]$ determined by mapping $K$ to the biggest negative zero of $h|_{x=1}$.
			The inverse map has the simple form
				\begin{equation*}
					K=K(z)=\frac{-1+z^2-z^3}{z^4}.
				\end{equation*}
			We will first deal with the cases $K\in\left[-\frac{25}{72},U\right)$ or, equivalently, $z\in\left(z_{\mathrm{min}},z_{\mathrm{max}}\right]$. The case $K=U$ will be treated separately afterwards. For our desired coordinate change we want that the transformed polynomial vanishes along $\{y=0\}$, so we set for $z\in\left(z_{\mathrm{min}},z_{\mathrm{max}}\right]$, $\widetilde{y}=zx+y$. Now, for our transformed polynomial $\widetilde{h}:=h\left(\left(\begin{smallmatrix}x\\ \widetilde{y}\end{smallmatrix}\right)\right)$, we need to find a coordinate transformation leaving $\{y=0\}$ invariant, such that our yet again transformed $\widetilde{h}$ vanishes along $\{x=0\}$. To do so, we need to solve the equation $\widetilde{h}|_{y=1}=0$ for $x$, and we obtain as one of the solutions
				\begin{center}
					\scalebox{0.65}{
						$F:=\frac{\begin{array}{l}\left(12\sqrt{3}(z^3 - 2z^2 + 4)\sqrt{23z^6 - 36z^5 + 20z^4 + 40z^3 - 32z^2 + 16} + 100z^6 - 288z^5 + 288z^4 - 288z^2\right)^{1/3} -6z^3 + 10z^2 - 12\\
							+ (4z^4 - 48z^3 + 48z^2 - 48)\left(12\sqrt{3}(z^3 - 2z^2 + 4)\sqrt{23z^6 - 36z^5 + 20z^4 + 40z^3 - 32z^2 + 16} + 100z^6 - 288z^5 + 288z^4 - 288z^2\right)^{-1/3}\end{array}}{\begin{array}{l}6(z^3 - 2z^2 + 4)z\end{array}}.$}
				\end{center}
			Note that $F$ is real for all $z\in\left(z_{\mathrm{min}},z_{\mathrm{max}}\right]$. The main difficulty in proving the latter is showing that the term $s:=23z^6 - 36z^5 + 20z^4 + 40z^3 - 32z^2 + 16$ is positive for all such $z$. The easiest way to show this is to check that all real solutions of $\partial^2_z s$, which can be found with the help of any modern computer algebra software like Maple, are bigger than $z_{\mathrm{max}}$, then showing that $\partial^2_z s$ is positive in $\left[z_{\mathrm{min}},z_{\mathrm{max}}\right]$, checking that $\partial_z s(z_{\mathrm{max}})<0$, and finally confirming that $s(z_{\mathrm{max}})>0$. We set $\widetilde{x}=x+Fy$ and $\overline{h}:=h\left(\left(\begin{smallmatrix}\widetilde{x}\\ \widetilde{y}\end{smallmatrix}\right)\right)$. By construction, $\overline{h}$ vanishes along $\{x=0\}$ and $\{y=0\}$, and, hence, $xy$ divides $\overline{h}$. To construct our final coordinate transformation we need to show that the two terms
				\begin{equation*}
					f_1 :=\frac{1}{2}\partial^2_x\left(\frac{\overline{h}}{xy}\right)=\frac{-z^3 + 2z^2 - 4}{z},\quad
					f_2 :=\frac{1}{2}\partial^2_y\left(\frac{\overline{h}}{xy}\right),
				\end{equation*}
			are positive for all $z\in\left(z_{\mathrm{min}},z_{\mathrm{max}}\right]$. Showing this for $f_1$ is not difficult. But the term $f_2$ is more complicated, so before even writing it out, we introduce the abbreviations
				\begin{align*}
					p_1&:=z^4 - 12z^3 + 12z^2 - 12,\\
					p_2&:=z^3 - 2z^2 + 4,\\
					p_3&:=23z^6 - 36z^5 + 20z^4 + 40z^3 - 32z^2 + 16,\\
					p_4&:=z^4 - 72/25z^3 + 72/25z^2 - 72/25,\\
					p_5&:=100z^6 - 288z^5 + 288z^4 - 288z^2,\\
					p_6&:=z^{12} - \frac{504}{89}z^{11} + \frac{7848}{623}z^{10} - \frac{5184}{623}z^9 - \frac{1368}{89}z^8 + \frac{19008}{623}z^7\\
					&\ - \frac{1728}{623}z^6 - \frac{20736}{623}z^5 + \frac{14688}{623}z^4 + \frac{10368}{623}z^3 - \frac{10368}{623}z^2 + \frac{3456}{623},
				\end{align*}
			and
				\begin{equation*}
					\begin{array}{lll}
						q_1:=p_1,&
						q_2:=zp_2,&
						q_3:=\frac{3\sqrt{3}}{25}p_2\sqrt{p_3}+z^2p_4,\\
						q_4:=z^2p_4,&
						q_5:=12\sqrt{3}p_2\sqrt{p_3}+p_5,&
						q_6:=\frac{75\sqrt{3}}{623}z^2p_2p_4\sqrt{p_3}+p_6.
					\end{array}
				\end{equation*}
			Note that $q_1$ has non-constant sign and a zero in $\left(z_{\mathrm{min}},z_\mathrm{max}\right)$. The term $q_2$ is negative on $\left(z_{\mathrm{min}},z_\mathrm{max}\right]$ and vanishes at $z_{\mathrm{min}}$, and $q_4$ is positive on $\left(z_{\mathrm{min}},z_\mathrm{max}\right)$ and vanishes at $z_\mathrm{max}$. The remaining terms $q_3,q_5,q_6$ are positive on $\left(z_{\mathrm{min}},z_\mathrm{max}\right)$. We leave it as a small exercise to the reader to prove these claims. We can write $f_2$ as
				\begin{equation*}
					f_2=\frac{q_5^{1/3}}{-1200z^2q_2q_3}\left(16q_1^2+4q_1q_5^{2/3}+q_5^{4/3}\right),
				\end{equation*}
			which allows us to see that $f_2>0$ on $\left(z_{\mathrm{min}},z_\mathrm{max}\right]$. Having shown the positivity of $f_1$ and $f_2$, we can construct our final coordinate transformation. Using
				\begin{align*}
					\overline{x}:=f_1^{-3/8}f_2^{1/8}x,\quad
					\overline{y}:=f_1^{1/8}f_2^{-3/8}y,
				\end{align*}
			we define $\widehat{h}:=\overline{h}\left(\left(\begin{smallmatrix}\overline{x}\\ \overline{y}\end{smallmatrix}\right)\right)$, which is of the form
				\begin{equation*}
					\widehat{h}=xy(x^2+cxy+y^2),
				\end{equation*}
			with $c$ being given by
				\begin{equation*}
					c=\frac{-4\sqrt{3}\left(25q_1^2q_3q_5^{2/3}+1246q_1q_5^{1/3}q_6+31150q_4q_6\right)\sqrt{q_5}}{\left(50q_1^2q_3q_5^{1/3}+2492q_1q_6+623q_5^{2/3}q_6\right)\sqrt{400q_1q_3+16q_1^2q_5^{1/3}+q_5^{5/3}}}.
				\end{equation*}
			Our next aim is to show that $c:\left(z_{\mathrm{min}},z_\mathrm{max}\right]\to\mathbb{R}$ is strictly monotonously increasing and that its image is given by $(-2,0]$. We first show, with the help of computer algebra software, that $\{\partial_z c=0\}\subset\{7z^4 - 8z^3 + 8z^2 - 8=0\}$, and we verify that $\{7z^4 - 8z^3 + 8z^2 - 8=0\}\cap\left(z_{\mathrm{min}},z_\mathrm{max}\right]=\emptyset$. After checking the sign of $\partial_z c$ at $\frac{z_{\mathrm{min}}+z_\mathrm{max}}{2}$ we conclude that $c$ is, in fact, strictly monotonously increasing. We now simply evaluate $c$ at $z_\mathrm{min}$ and $z_\mathrm{max}$ and obtain $-2$ and $0$, respectively. We conclude $\mathrm{Im}(c)=(-2,0]$ as claimed, and furthermore that $c$ is a bijection. Hence, counting the connected components of $\{h=1\}\cap\{\text{hyperbolic points of }h\}$ and showing their pairwise equivalence for each $K\in\left[-\frac{25}{72},U\right)$ is equivalent to doing so for $\{\widehat{h}=1\}\cap\{\text{hyperbolic points of }\widehat{h}\}$ for each $c\in(-2,0]$. For all $c\in(-2,0]$, $\widehat{h}$ is invariant under switching $x$ and $y$, and under the negative identity transformation. Furthermore, ${\widehat{h}>0}=\{x>0,\ y>0\}\cup\{x<0,\ y<0\}$. Let
				\begin{equation*}
					\varphi:=\det\left(\left.\partial^2\widehat{h}\right|_{x=1/2+s,\ y=1/2-s}\right).
				\end{equation*}
			The function $\varphi$ is symmetric in $s$ and we find that $\varphi\left(\pm\frac{1}{2}\right)=-9<0$. We further calculate $\varphi(0)=-\frac{3}{4}c^2-\frac{3}{2}c$, which is $0$ for $c=0$ and positive for $c\in(-2,0)$. For $c\in(-2,0)$, the real zeros of $\varphi$ in $s$ are given by
				\begin{equation*}
					\pm\frac{\sqrt{c^2 + 2\sqrt{-2c^2 + 9} - 6}}{\sqrt{c(c - 2)}}.
				\end{equation*}
			For $c=0$, the only real solution of $\varphi=0$ is $s=0$. We conclude that for each $c\in(-2,0]$, $\{\widehat{h}=1\}\cap\{\text{hyperbolic points of }\widehat{h}\}$ has $4$ equivalent connected components as claimed. Note that in the boundary case $c=0$ or, equivalently, $K=-\frac{25}{72}$, the cone $\{h>0\}\cap\{\text{non-hyperbolic points of }h\}$ is a degenerate cone, meaning that it consists only of two rays in $\mathbb{R}^2$. Hence, each connected component of the curve $\{h=1\}$ has a single flat point.
			
			Lastly, we need to study the case $K=U$. In that case, $h$ is not equivalent to $\widehat{h}$ with $c=-2$. This follows from the fact that for $K=U$, respectively $c=-2$, $\{h=0\}$ is just a line in $\mathbb{R}^2$, but $\{\widehat{h}=0\}$ consists of three distinct lines in $\mathbb{R}^2$. Instead, $\widehat{h}$ for $c=-2$ is equivalent to Thm. \ref{thm_quartic_CCGPSR_curves_classification} \hyperref[eqn_qCCPSRcurves_class_c]{c)} via
				\begin{equation}
					A=\left(\begin{matrix}
						2^{1/2}3^{-3/4} & 0\\
						2^{-1/2}3^{3/4} & -2^{-1/2}3^{-1/4}
					\end{matrix}\right),\label{eqn_u_graph_A_trafo}
				\end{equation}
			so that $A^*\left(x^4-x^2y^2+\frac{2}{3\sqrt{3}}xy^3\right)=\widehat{h}$ for $c=-2$. As for $z\in\left(z_\mathrm{min},z_\mathrm{max}\right]$ we start with transforming $h$ via $\widetilde{y}=z_\mathrm{min}x+y$ into $\widetilde{h}:=h\left(\left(\begin{smallmatrix}x\\ \widetilde{y}\end{smallmatrix}\right)\right)$, so that $\widetilde{h}$ vanishes on $\{y=0\}$. In contrast to the previous cases, $\widetilde{h}$ does not vanish and is in fact positive on $\{y\ne0\}$, see Figure \ref{fig_h_tilda_U_case} for a plot of $\{\widetilde{h}=1\}$.
				\begin{figure}[H]
					\centering
					\includegraphics[scale=0.2]{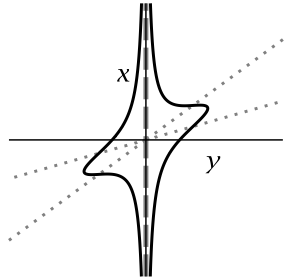}
					\caption{The set $\{\widetilde{h}=1\}$ in black, the boundary of $\{\text{hyperbolic points of }\widetilde{h}\}$ as dotted grey lines, and the set $\{\widetilde{h}=0\}$ as a dashed grey line.}\label{fig_h_tilda_U_case}
				\end{figure}
			\noindent
			We want to find a linear transformation so that $\widetilde{h}$ is symmetric in both $x$ and $y$. To find such a transformation, we make the following ansatz. For $\overline{h}:=\widetilde{h}\left(\left(\begin{smallmatrix}x+Cy\\y\end{smallmatrix}\right)\right)$, $C\in\mathbb{R}$, we solve the equation
				\begin{equation*}
					\left.\partial_x\overline{h}\right|_{x=0,\ y=1}=0
				\end{equation*}
			for $C$ and find the unique solution
				\begin{equation*}
					C=\frac{24 +(2\sqrt{57} - 18)(46 + 6\sqrt{57})^{1/3} +  (\sqrt{57} - 9)(46 + 6\sqrt{57})^{2/3}}{ - 80 + (2\sqrt{57} - 2)(46 + 6\sqrt{57})^{1/3} + (-5\sqrt{57} + 37)(46 + 6\sqrt{57})^{2/3}}.
				\end{equation*}
			Using the above value for $C$, we check that $\overline{h}$ is indeed symmetric in $x$ and $y$. It is now straightforward to check that $\left.\det\left(\partial^2\overline{h}\right)\right|_{y=1}$ has precisely two real zeros, both of which are simple zeros, and that it is positive at $x=0$. Together with the symmetry of $\overline{h}$ we conclude that $\{\overline{h}=1\}\cap\{\text{hyperbolic points of }\overline{h}\}$ has $4$ equivalent connected components and, hence, the same holds true for $h$.
			
			Now we will show that for all $K\in\left[-\frac{25}{72},U\right]$, $G^h_0$ consists only of the identity transformation. Every element of $G^h_0$ is necessarily of the form $A(p)$ as in equation \eqref{eqn_A_explicit}. Suppose that there exists $p\in\mathcal{H}$, such that $A(p)^*h=h$ and $p\ne\left(\begin{smallmatrix}1\\0\end{smallmatrix}\right)$. Since $\mathcal{V}$ vanishes at no point of the form $(1,K)^T$ for $K\in\left[-\frac{25}{72},U\right]$ and by the correspondence of maximal integral curves of $\mathcal{V}|_{\{\mathcal{V}\ne0\}}$ and equivalence classes of inhomogeneous maximal quartic GPSR curves this would imply the existence of a periodic orbit of $\mathcal{V}$ restricted to the set consisting of points $(L,K)^T\in\mathbb{R}^2$ corresponding to incomplete maximal quartic GPSR curves. Finally, by the connectedness and simply connectedness of that set, this would imply the existence of a fix point of $\mathcal{V}$ therein. This is a contradiction by the fact that points $(L,K)^T\in\mathbb{R}^2$ corresponding to incomplete maximal quartic GPSR curves are contained in $\{\mathcal{V}\ne0\}$. Hence, $G^h_0$ consist only of the identity, and our arguments would have also worked for the family of curves Thm. \ref{thm_incomplete_quartics} \hyperref[eqn_incomplete_qGPSRcurves_class_a]{a)}.
			
			Since for each case of Thm. \ref{thm_incomplete_quartics} \hyperref[eqn_incomplete_qGPSRcurves_class_b]{b)}, \hyperref[eqn_incomplete_qGPSRcurves_class_c]{c)}, and \hyperref[eqn_incomplete_qGPSRcurves_class_d]{d)} the maximal incomplete GPSR curves are inhomogeneous, we obtain by our investigations that in these cases $G^h=\mathbb{Z}_2\times\mathbb{Z}_2$ by the fact that, in the right coordinates, $\{h=1\}$ is symmetric under two distinct mirror transformations and that $G^h_0$ consists only of the identity.
			
			Next, we need to show that the curves Thm. \ref{thm_incomplete_quartics} \hyperref[eqn_incomplete_qGPSRcurves_class_b]{b)} and \hyperref[eqn_incomplete_qGPSRcurves_class_c]{c)} are not singular at infinity and that the curve Thm. \ref{thm_incomplete_quartics} \hyperref[eqn_incomplete_qGPSRcurves_class_d]{d)} is singular at infinity. Recall the definitions of $t_m$ \eqref{eqn_t_m_def} and $t_M$ \eqref{eqn_t_M_def} as the non-zero solutions of $\dot{f}_{L,K}(t)=0$ for $K\ne 0$, where $f_{L,K}(t)$ was defined in \eqref{eqn_fLK}. A maximal connected quartic GPSR curve $\mathcal{H}_{L,K}$ with $K\ne 0$ might thus only be singular at infinity if $t_m\in \partial\,\mathrm{dom}(\mathcal{H}_{L,K})$ or $t_M\in \partial\,\mathrm{dom}(\mathcal{H}_{L,K})$. For $K=0$ and $L>0$, the unique non-zero solution of $\dot{f}_{L,0}(t)=0$ is $t=\frac{2}{3L}$. Since $f_{1,0}\left(\frac{2}{3}\right)=\frac{23}{27}>0$, it thus follows that $\mathcal{H}_{1,0}$ is not singular at infinity. We further obtain that for all $K\in\left[-\frac{25}{72},U\right]\setminus\{0\}$
				\begin{align*}
					f_{1,K}(t_m)&=\frac{-27-144K-128K^2+512K^3+(9+32K)^{\frac{3}{2}}}{512K^3},\\ f_{1,K}(t_M)&=\frac{-27-144K-128K^2+512K^3-(9+32K)^{\frac{3}{2}}}{512K^3}.
				\end{align*}
			Solving $f_{1,K}(t_M)=0$, we obtain precisely one real solution, namely $K=U$, and we find that $f_{1,K}(t_m)=0$ has no real solutions. This shows that Thm. \ref{thm_incomplete_quartics} \hyperref[eqn_incomplete_qGPSRcurves_class_b]{b)} and \hyperref[eqn_incomplete_qGPSRcurves_class_c]{c)} are not singular at infinity. For all $K>0$ the term $t_M$ is negative. Hence, in order to obtain that Thm. \ref{thm_incomplete_quartics} \hyperref[eqn_incomplete_qGPSRcurves_class_d]{d)} is singular at infinity, it suffices to show that for all $y\in\left(t_M|_{K=U},0\right]$, $(1,y)^T$ is a hyperbolic point of $h_{1,U}$. We check this by solving $\det(\partial^2 h_{1,U})|_{x=1}=0$ and obtaining that the unique non-positive solution is given precisely by the unique real solution of $f_{1,U}=0$.

			We need to show that none of the curves Thm. \ref{thm_incomplete_quartics} \hyperref[eqn_incomplete_qGPSRcurves_class_b]{b)}, \hyperref[eqn_incomplete_qGPSRcurves_class_c]{c)}, or \hyperref[eqn_incomplete_qGPSRcurves_class_d]{d)} are equivalent to Thm. \ref{thm_incomplete_quartics} \hyperref[eqn_incomplete_qGPSRcurves_class_a]{a)}. This follows by comparing their respective automorphism groups.
			
			The curve Thm. \ref{thm_incomplete_quartics} \hyperref[eqn_incomplete_qGPSRcurves_class_d]{d)} is neither equivalent to the curve Thm. \ref{thm_incomplete_quartics} \hyperref[eqn_incomplete_qGPSRcurves_class_b]{b)} or any curve in Thm. \ref{thm_incomplete_quartics} \hyperref[eqn_incomplete_qGPSRcurves_class_c]{c)} since Thm. \ref{thm_incomplete_quartics} \hyperref[eqn_incomplete_qGPSRcurves_class_d]{d)} is the only maximal incomplete quartic CCGPSR curve up to equivalence that is singular at infinity.
			
			Any curve in Thm. \ref{thm_incomplete_quartics} \hyperref[eqn_incomplete_qGPSRcurves_class_b]{b)} is not equivalent to Thm. \ref{thm_incomplete_quartics} \hyperref[eqn_incomplete_qGPSRcurves_class_c]{c)} since the latter curve has an isolated flat point, while the former curves do not admit such a point. This can be checked by bringing the defining polynomials of curves Thm. \ref{thm_incomplete_quartics} \hyperref[eqn_incomplete_qGPSRcurves_class_b]{b)} to the form $h=xy(x^2+cxy+y^2)$ for fitting $c$ so that $\mathcal{H}$ is contained in the upper right quadrant, and checking that the negative Hessian of $h$ at $(x,y)^T=(1,1)^T$ is negative.
			
			Lastly, no two distinct curves in Thm. \ref{thm_incomplete_quartics} \hyperref[eqn_incomplete_qGPSRcurves_class_b]{b)} are equivalent. This follows from the fact that up to a sign change in the $y$-coordinate, every transformation preserving the standard form of $h$ must be of the form \eqref{eqn_A_explicit}. In the proof of this theorem we have however seen that if one such transformation, and thereby a curve of such transformations connecting two distinct elements of Thm. \ref{thm_incomplete_quartics} \hyperref[eqn_incomplete_qGPSRcurves_class_b]{b)}, would exist, the constructed homotopy $H$ \eqref{eqn_H_homotopy} would not have the described transversality property.

			This finishes the proof of Theorem \ref{thm_incomplete_quartics}.
		\end{proof}
	\end{Th}
	
\section{Applications}\label{sect_applications}
	The defining polynomial of a PSR manifold in standard form is of the form $h=x^3-x\langle y,y\rangle+P_3(y)$. In \cite[Prop.\,4.6]{L2} it is shown that homogeneous PSR manifolds are singular at infinity, which was shown to imply \cite[Rem.\,4.7]{L2} that their respective standard form has the property that every critical point $p$ of $P_3$ restricted to $S^{n-1}=\{\langle y,y\rangle=1\}$ fulfils $P_3(p)=\frac{2}{3\sqrt{3}}$. For PSR curves the other direction also holds true, meaning that any PSR curve in standard form with the property that $P_3(1)=\pm\frac{2}{3\sqrt{3}}$ is automatically homogeneous. For higher dimensional PSR manifolds it is not known whether this is true as well and this is the topic of current research of the author. In this work, we will show that homogeneous quartic PSR manifolds are, just like their PSR counterparts, singular at infinity. We will also characterise the corresponding terms $P_3$ and $P_4$ of their defining polynomials in standard form by studying their possible critical values on $S^{n-1}$.
		
	\begin{Prop}\label{prop_hom_quartics_sing_at_inf}
		Homogeneous quartic GPSR manifolds are singular at infinity.
		\begin{proof}
			Let $\mathcal{H}$ be a homogeneous quartic GPSR manifold. Since $\mathcal{H}$, as a Riemannian homogeneous space, is geodesically complete, we can assume without loss of generality that $\mathcal{H}$ is a quartic CCGPSR manifold by restriction to one of its connected components. Furthermore we can assume that $\mathcal{H}$ is in standard form and that $\dim(\mathcal{H})\geq 2$, since in the $1$-dimensional case we have already proven in Theorem \ref{thm_quartic_CCGPSR_curves_classification} \hyperref[eqn_qCCPSRcurves_class_a]{a)} and \hyperref[eqn_qCCPSRcurves_class_b]{b)} that being homogeneous implies singular at infinity. Then by Proposition \ref{prop_hom_iff_deltaPs_vanish}, the corresponding terms $\delta P_3(y)$ \eqref{eqn_deltaP3_quartics} and $\delta P_4(y)$ \eqref{eqn_deltaP4_quartics} vanish identically for a choice of the linear map $b:\mathbb{R}^n\to\mathfrak{so}(n)$, which we from now on assume to be fixed to fulfil that requisite. Let $\overline{y}$ be a critical point of the restriction of $h\left(\left(\begin{smallmatrix}1\\y\end{smallmatrix}\right)\right)=1-\langle y,y\rangle + P_3(y) + P_4(y)$ to $S^{n-1}=\{\langle y,y\rangle=1\}$. Then $\overline{y}$ fulfils
				\begin{equation}\label{eqn_crit_point_hom_quartics}
					-2\langle\overline{y},\D y\rangle + \D P_3|_{\overline{y}} + \D P_4|_{\overline{y}}=r\langle\overline{y},\D y\rangle
				\end{equation}
			for some $r\in\mathbb{R}$, and we thus have
				\begin{equation}\label{eqn_crit_point_on_kernel}
					\D P_3|_{\overline{y}}(u) + \D P_4|_{\overline{y}}(u)=0
				\end{equation}
			for all $u\in{\overline{y}}^\bot=\ker\langle \overline{y},\cdot\rangle$. After a possible rotation in the $y$-coordinates, we can assume without loss of generality that $\overline{y}=(0,\ldots,0,1)^T$ and, with the notation $y=(v_1,\ldots,v_{n-1},w)^T=\left(\begin{smallmatrix}v\\w\end{smallmatrix}\right)$, have that $P_3$ and $P_4$ are of the form
				\begin{align*}
					P_3&= aw^3 + w^2\langle \eta,v\rangle + w q(v) + c(v),\\
					P_4&= Aw^4 + w^3\langle E,v\rangle + w^2 Q(v) + w C(v) + F(v),
				\end{align*}
			with fittingly chosen $a,A\in\mathbb{R}$, $\eta,E\in\mathbb{R}^{n-1}$, $q,Q$ quadratic forms, $c,C$ cubic forms, and $F$ a quartic form. In these coordinates, every $u\in{\overline{y}}^\bot$ is of the form $u=\left(\begin{smallmatrix}U\\0\end{smallmatrix}\right)$. By writing out equation \eqref{eqn_crit_point_on_kernel} applied to $\left(\begin{smallmatrix}U\\0\end{smallmatrix}\right)$ for $U\in\mathbb{R}^{n-1}$ arbitrary, we obtain that $\eta=-E$ must hold. We will use that and the homogeneity of $\mathcal{H}$ that $\eta$ and $E$ must, if $r\ne 0$, vanish. Applying both sides of equation \eqref{eqn_crit_point_hom_quartics} to $b \overline{y} + \frac{1}{4}\partial^2P_3|_{\overline{y}}\D y$ yields
				\begin{equation*}
					\left(\D P_3|_{\overline{y}}+\D P_4|_{\overline{y}}\right)\left(b \overline{y} + \frac{1}{4}\partial^2P_3|_{\overline{y}}\D y\right) = \frac{r+2}{2}\D P_3|_{\overline{y}}.
				\end{equation*}
			Hence, $\delta P_3(y)=0$, $\delta P_4(y)=0$, and consequently $\delta P_3(y)+\delta P_4(y)=0$ implies that
				\begin{equation*}
					\left(\frac{1}{2}P_3(\overline{y})-1\right)\langle \overline{y},\D y\rangle + \frac{r+2}{2}\D P_3|_{\overline{y}} + \D P_4|_{\overline{y}}=0.
				\end{equation*}
			The above equation together with with \eqref{eqn_crit_point_on_kernel} now implies $\frac{r}{2}\D P_3|_{\overline{y}}(u)=0$ for all $u\in {\overline{y}}^\bot$, which in our chosen coordinates is equivalent to $r\eta=0$. Hence, $\eta=-E=0$ for $r\ne 0$. If $r=0$, we obtain by applying \eqref{eqn_crit_point_hom_quartics} to $\overline{y}$ that $-2+3a+4A=0$. Since the intersection of $\mathcal{H}$ with the plane $\mathrm{span}\left\{\left(\begin{smallmatrix}1\\0\\0\end{smallmatrix}\right),\left(\begin{smallmatrix}1\\0\\1\end{smallmatrix}\right)\right\}$ is a quartic CCGPSR curve, we have by the proof of Theorem \ref{thm_quartic_CCGPSR_curves_classification} that $A\leq \frac{1}{4}$ must hold and that $|a|\leq \mathbf{u}(A)$ if $A\in\left[-\frac{1}{12},\frac{1}{4}\right]$ or $|a|<\mathbf{w}(A)$ if $A\in\left(-\infty,-\frac{1}{12}\right)$, cf. equations \eqref{eqn_upper_bdr} and \eqref{eqn_sharp_bdr}. But $r=0$ implies $a=\frac{2-4A}{3}$. We leave it as an easy exercise for the reader to check that
				\begin{equation*}
					\frac{2-4A}{3}>\mathbf{u}(A)
				\end{equation*}
			for all $A\in\left[-\frac{1}{12},\frac{1}{4}\right]$, and
				\begin{equation*}
					\frac{2-4A}{3}>\mathbf{w}(A)
				\end{equation*}
			for all $A\in\left(-\infty,-\frac{1}{12}\right)$, showing that $r=0$ is in contradiction to $\mathcal{H}$ being closed in its ambient space. Hence, we can exclude the case $r=0$ and have thereby shown that $\eta$ and $E$ must vanish. Next, we calculate
				\begin{align*}
					\delta P_3(y)\left(\left(\begin{smallmatrix}0\\1\end{smallmatrix}\right)\right)&= \left(\frac{9}{2}a^2+4A-1\right)w^3 + \text{ terms of lower order in }w,\\
					\delta P_4(y)\left(\left(\begin{smallmatrix}0\\1\end{smallmatrix}\right)\right)&= a\left(6A+\frac{1}{2}\right)w^4 + \text{ terms of lower order in }w,
				\end{align*}
			which, up to a relabelling of the variables, coincides with \eqref{eqn_deltaP3_quarticcurves} and \eqref{eqn_deltaP4_quarticcurves} in the highest order $w$ term. Since $\delta P_3(y)$ and $\delta P_4(y)$ vanish identically by the homogeneity of $\mathcal{H}$, we obtain that $a$ and $A$ must fulfil either $a=0$ and $A=\frac{1}{4}$, or $a=\pm\frac{2\sqrt{2}}{3\sqrt{3}}$ and $A=-\frac{1}{12}$, meaning that the quartic CCGPSR curve obtained from intersecting $\mathcal{H}$ with $\mathrm{span}\left\{\left(\begin{smallmatrix}1\\0\\0\end{smallmatrix}\right),\left(\begin{smallmatrix}1\\0\\1\end{smallmatrix}\right)\right\}$ is equivalent to either Theorem \ref{thm_quartic_CCGPSR_curves_classification} \hyperref[eqn_qCCPSRcurves_class_a]{a)} or \hyperref[eqn_qCCPSRcurves_class_b]{b)}. Since both of these curves are singular at infinity, we obtain with Lemma \ref{lem_sing_at_inf_iff_restr_to_curve} that $\mathcal{H}$ itself singular at infinity.
		\end{proof}
	\end{Prop}
	
	In the proof of the above proposition we have in particular shown the following property of homogeneous quartic GPSR manifolds that was announced in Remark \ref{rem_motivation_quartic_homs_sing_at_inf}.
	
	\begin{Cor}\label{cor_crit_values_quartic_homs}
		For each critical value $\overline{y}$ of the terms $P_3(y)$ and $P_4(y)$ of a homogeneous quartic GPSR manifold in standard form restricted to $S^{n-1}$, it holds that either
			\begin{equation*}
				P_3(\overline{y})=0\quad\text{and}\quad P_4(\overline{y})=\frac{1}{4},
			\end{equation*}
		or
			\begin{equation*}
				|P_3(\overline{y})|=\frac{2\sqrt{2}}{3\sqrt{3}}\quad\text{and}\quad P_4(\overline{y})=-\frac{1}{12}.
			\end{equation*}
	\end{Cor}
	
	One possible way to tell whether two given maximal connected GPSR manifolds defined by polynomials of the same degree are equivalent is to study their limit geometries. We use an analogue of \cite[Def.\,1.3]{L3} allowing for degrees higher than three and without the restriction to closed connected PSR manifolds to give that term a precise meaning.
	
	\begin{Def}\label{def_limit_geo}
		Let $\mathcal{H}\subset\{h=1\}$ be an $n\geq 1$-dimensional maximal connected (G)PSR manifold. Then a maximal connected (G)PSR manifold $\overline{\mathcal{H}}\subset\{\overline{h}=1\}$, $\deg(h)=\tau\geq 3$, is called \textbf{limit geometry} of $\mathcal{H}$ if the following it can be constructed in the following way:\\
		Assume that exists a point $p\in\mathcal{H}$ and a standard form $\widetilde{h}:=A(p)^*h$ of $h$ in the sense of Definition \ref{def_standard_form}, such that there exists an affine-linear curve $\gamma:[0,1)\to\mathbb{R}^{n}$ of the form $\gamma(t)=\left(\begin{smallmatrix}1\\tv\end{smallmatrix}\right)$ for some $v\in\mathbb{R}^n\setminus0$, fulfilling $\widetilde{h}(\gamma(1))=0$ and $\gamma([0,1))\subset\mathrm{dom}(A(p)^{-1}(\mathcal{H}))$, and a curve of homogeneous hyperbolic polynomials in standard form $(A\circ\gamma)^*h$ coinciding with $\widetilde{h}$ at $t=0$, such that the limit
			\begin{equation}
				\widehat{h}:=\lim\limits_{t\to 1} A(\gamma(t))^*h\label{eqn_limit_geo_defining_eqn}
			\end{equation}
		exists and such that $\left(\begin{smallmatrix}1\\0\end{smallmatrix}\right)$ is a hyperbolic point of $\widehat{h}$. Then the maximal connected (G)PSR manifold in standard form $\widehat{\mathcal{H}}\subset\{\widehat{h}=1\}$ is called a limit geometry of $\mathcal{H}$, and $\overline{\mathcal{H}}$ is called a limit geometry of $\mathcal{H}$ if it is equivalent to a maximal connected (G)PSR manifold $\widehat{\mathcal{H}}$ constructed in this way.
	\end{Def}
	
	The above definition is at first look admittedly slightly unintuitive. But it turns out that it leads to interesting results. In \cite{L3} we have shown that when restricting to complete CCPSR manifolds, there indeed always exists a limit geometry, independently of the choice of the plane $P$ or the considered direction. In dimension two every possible limit geometry of CCPSR surfaces is a homogeneous space, cf. \cite[Thm.\,1.11\,(ii)]{L3}. We do not expect this statement to be true in the case of quartic surfaces. But, as we will show next, if a maximal connected quartic GPSR curve has a limit geometry, then it is automatically homogeneous, meaning either equivalent to Thm. \ref{thm_quartic_CCGPSR_curves_classification} \hyperref[eqn_qCCPSRcurves_class_a]{a)} or \hyperref[eqn_qCCPSRcurves_class_b]{b)}. Geometrically, a limit geometry gives insight into the asymptotic behaviour of the centro-affine fundamental form $g_\mathcal{H}$ along a centrally projected affine line to a maximal connected (G)PSR manifold, cf. \cite[Prop.\,1.5]{L3}. If a limit geometry $\overline{\mathcal{H}}$ with centro-affine fundamental form $g_{\overline{\mathcal{H}}}$ of $\mathcal{H}$ corresponding to a certain curve is defined, then for all $\varepsilon>0$ and all compactly embedded neighbourhoods $U$ of $\left(\begin{smallmatrix}1\\0\end{smallmatrix}\right)$ in $\overline{\mathcal{H}}$ there exists a diffeomorphism $F:\overline{U}\to F(\overline{U})\subset\mathcal{H}$, such that
		\begin{equation*}
			\left\|g_{\overline{\mathcal{H}}}-F^*g_{\mathcal{H}}\right\|_{g_{\overline{\mathcal{H}}}}<\varepsilon
		\end{equation*}
	in $\overline{U}$.
	
	\begin{Prop}\label{prop_limit_geos}
		If a maximal connected quartic GPSR curve has a limit geometry, it is equivalent to one of the homogeneous curves Thm. \ref{thm_quartic_CCGPSR_curves_classification} \hyperref[eqn_qCCPSRcurves_class_a]{a)} or \hyperref[eqn_qCCPSRcurves_class_b]{b)}. The limit geometries are characterised as follows:
			\begin{center}
				\scalebox{0.75}{
					$\begin{array}{|c|c|c|c|c|c|}
						\hline \mathcal{H} & \text{closed} & \text{homogeneous} & \text{limit geometries} & \sharp\text{ directions with limit geometry}\\
						\hline
						\hline
						\text{Thm. \ref{thm_quartic_CCGPSR_curves_classification} \hyperref[eqn_qCCPSRcurves_class_a]{a)}}\vphantom{\bigg|} & \text{yes} & \text{yes} & \text{Thm. \ref{thm_quartic_CCGPSR_curves_classification} \hyperref[eqn_qCCPSRcurves_class_a]{a)}} & 2 \\
						\hline
						\text{Thm. \ref{thm_quartic_CCGPSR_curves_classification} \hyperref[eqn_qCCPSRcurves_class_b]{b)}}\vphantom{\bigg|} & \text{yes} & \text{yes} & \text{Thm. \ref{thm_quartic_CCGPSR_curves_classification} \hyperref[eqn_qCCPSRcurves_class_b]{b)}} & 2 \\
						\hline
						\text{Thm. \ref{thm_quartic_CCGPSR_curves_classification} \hyperref[eqn_qCCPSRcurves_class_c]{c)}}\vphantom{\bigg|} & \text{yes} & \text{no} & \text{Thm. \ref{thm_quartic_CCGPSR_curves_classification} \hyperref[eqn_qCCPSRcurves_class_a]{a)} \& \hyperref[eqn_qCCPSRcurves_class_b]{b)}} & 2 \\
						\hline
						\text{Thm. \ref{thm_quartic_CCGPSR_curves_classification} \hyperref[eqn_qCCPSRcurves_class_d]{d)}}\vphantom{\bigg|} & \text{yes} & \text{no} & \text{Thm. \ref{thm_quartic_CCGPSR_curves_classification} \hyperref[eqn_qCCPSRcurves_class_b]{b)}} & 2 \\
						\hline
						\text{Thm. \ref{thm_incomplete_quartics}\vphantom{\bigg|} \hyperref[eqn_incomplete_qGPSRcurves_class_a]{a)}} & \text{no} & \text{no} & \text{none} & 0 \\
						\hline
						\text{Thm. \ref{thm_incomplete_quartics}\vphantom{\bigg|} \hyperref[eqn_incomplete_qGPSRcurves_class_b]{b)}} & \text{no} & \text{no} & \text{Thm. \ref{thm_quartic_CCGPSR_curves_classification} \hyperref[eqn_qCCPSRcurves_class_b]{b)}} & 1 \\
						\hline
						\text{Thm. \ref{thm_incomplete_quartics}\vphantom{\bigg|} \hyperref[eqn_incomplete_qGPSRcurves_class_c]{c)}} & \text{no} & \text{no} & \text{Thm. \ref{thm_quartic_CCGPSR_curves_classification} \hyperref[eqn_qCCPSRcurves_class_b]{b)}} & 1 \\
						\hline
						\text{Thm. \ref{thm_incomplete_quartics}\vphantom{\bigg|} \hyperref[eqn_incomplete_qGPSRcurves_class_d]{d)}} & \text{no} & \text{no} & \text{Thm. \ref{thm_quartic_CCGPSR_curves_classification} \hyperref[eqn_qCCPSRcurves_class_a]{a)}} & 1 \\
						\hline
					\end{array}$
				}
			\end{center}
		\begin{proof}
			Since we are dealing with curves, the term $E$ in the definition of $A$ \eqref{eqn_A_explicit} is uniquely determined up to sign. Hence, in \eqref{eqn_limit_geo_defining_eqn} we do not need to concern ourselves with the question whether the existence of the limit might depend on the choice of $E$ which in higher dimensions will depend on an $\mathrm{O}(n)$-section over the considered curve.
			
			We start with the closed cases. The curves Thm. \ref{thm_quartic_CCGPSR_curves_classification} \hyperref[eqn_qCCPSRcurves_class_a]{a)} and \hyperref[eqn_qCCPSRcurves_class_b]{b)} are homogeneous spaces. Any standard form of their defining polynomial $h_{L,K}$ corresponds to points $(L,K)^T\in\{\mathcal{V}=0\}$. Thus, the term on the right hand side of \eqref{eqn_limit_geo_defining_eqn} is for any allowed choice of $A$ and $\gamma$ constant in $t$. This means that Thm. \ref{thm_quartic_CCGPSR_curves_classification} \hyperref[eqn_qCCPSRcurves_class_a]{a)} and \hyperref[eqn_qCCPSRcurves_class_b]{b)} are their own limit geometries, respectively.
			
			For Thm. \ref{thm_quartic_CCGPSR_curves_classification} \hyperref[eqn_qCCPSRcurves_class_c]{c)}, recall that all the their standard forms $h_{L,K}$, assuming without loss of generality $L\geq 0$, of the defining polynomial fulfil $(L,K)^T\in \mathrm{Im}(\mathbf{u}(K),K)^T$, $K\in\left(-\frac{1}{12},\frac{1}{4}\right)$, cf. \eqref{eqn_partial_u} and the subsequent discussion. Thus the limit geometries of Thm. \ref{thm_quartic_CCGPSR_curves_classification} \hyperref[eqn_qCCPSRcurves_class_c]{c)} correspond to points in the boundary of the graph of $\mathbf{u}$, which are the two homogeneous curves Thm. \ref{thm_quartic_CCGPSR_curves_classification} \hyperref[eqn_qCCPSRcurves_class_a]{a)} and \hyperref[eqn_qCCPSRcurves_class_b]{b)}.
			
			The curves in Thm. \ref{thm_quartic_CCGPSR_curves_classification} \hyperref[eqn_qCCPSRcurves_class_d]{d)} need a little more work. On $U_1:=\left\{K\in\left[-\frac{1}{12},\frac{1}{4}\right),\ |L|<\mathbf{u}(K)\right\}$, $\partial_K$ is transversal to $\mathcal{V}$, which follows from $\D L(\mathcal{V})=0$ having no solution in $U_1$ and $\{\mathcal{V}=0\}\cap U_1=\emptyset$. The boundary of $U_1$ consists of images of maximal integral curves of $\mathcal{V}$. Hence, every maximal integral curve of $\mathcal{V}$ with initial condition in $U_1$ has image in $U_1$. By the transversality of $\partial_K$ and $\mathcal{V}$, and $\mathcal{V}$ nowhere vanishing on $U_1$ it thus follows that the boundary of the image of every such integral curve is contained in $\{\mathcal{V}=0\}$. Since we have shown in the proof of Theorem \ref{thm_quartic_CCGPSR_curves_classification} that every such integral curve meets $\left\{L=0,\ K<\frac{1}{4}\right\}$ precisely once and by the symmetry of $\D L(\mathcal{V})$ in $L$ we conclude that the boundary of the image of every maximal integral curve of $\mathcal{V}$ with initial condition in $U_1$ is given by $\left\{\left(\pm\frac{2\sqrt{2}}{3\sqrt{3}},-\frac{1}{12}\right)^T\right\}$. Hence, a limit geometry of $\mathcal{H}_{0,K}$ for $K\in\left[-\frac{1}{12},\frac{1}{4}\right)$ exists independent of the direction of $\gamma$ in \eqref{eqn_limit_geo_defining_eqn} and is in any case given by Thm. \ref{thm_quartic_CCGPSR_curves_classification} \hyperref[eqn_qCCPSRcurves_class_b]{b)}.
			
			Now consider $\mathcal{H}_{0,K}$ with $K<-\frac{1}{12}$ and recall the coordinate change described in \eqref{eqn_first_F_trafo}. Since $\mathcal{V}$ is nowhere vanishing on $U_2:=\left\{K<\frac{1}{4},\ |L|<\mathbf{w}(K)\right\}$ and, similar to $U_1$, the boundary of $U_2$ consists of images of maximal integral curves of $\mathcal{V}$, we obtain that images of maximal integral curves of $\mathcal{V}$ with initial condition in $U_2$ are contained in $U_2$. By $\mathcal{V}$ being nowhere vanishing on $U_2$, $\D \varphi(F_*\mathcal{V})<0$ and $\mathrm{sgn}(\D \psi(F_*\mathcal{V}))=\mathrm{sgn}(\varphi)$ on $F(U_2)=\{|\varphi|<1\}\times\left(-\infty,-\frac{1}{12}\right)$, see \eqref{eqn_first_F_V_pushforward}, and the symmetry of $\D \varphi(F_*\mathcal{V})$ we conclude that the boundary of the image of every maximal integral curve of $\mathcal{V}$ starting in $U_2$ is given by $\left\{\left(\pm\frac{2\sqrt{2}}{3\sqrt{3}},-\frac{1}{12}\right)^T\right\}$. As for the cases with $K\in\left[0,\frac{1}{4}\right)$, this implies that every limit geometry of $\mathcal{H}_{0,K}$ for any $K<-\frac{1}{12}$ and any direction of $\gamma$ in \eqref{eqn_limit_geo_defining_eqn} is given by Thm. \ref{thm_quartic_CCGPSR_curves_classification} \hyperref[eqn_qCCPSRcurves_class_b]{b)}.
			
			Next we consider the incomplete cases. The curves in Thm. \ref{thm_incomplete_quartics} \hyperref[eqn_incomplete_qGPSRcurves_class_a]{a)} do not have a well-defined limit geometry. This follows simply by definition, as we have seen that the set $\{h_{0,K}=0\}$ consists only of $(x,y)^T=(0,0)^T$ for all $K>\frac{1}{4}$, meaning that the intersection of $\{1\}\times\overline{\mathrm{dom}(\mathcal{H}_{0,K})}$ and $\{h_{0,K}=0\}$ is empty in these cases.
			
			For Thm. \ref{thm_incomplete_quartics} \hyperref[eqn_incomplete_qGPSRcurves_class_b]{b)}. Using the construction of the homotopy $H$ in \eqref{eqn_H_homotopy} in the proof of Theorem \ref{thm_incomplete_quartics} and the transversality argument of curves of the form $H(s,\cdot)$ and $\mathcal{V}$ we immediately obtain that each curve in Thm. \ref{thm_incomplete_quartics} \hyperref[eqn_incomplete_qGPSRcurves_class_a]{a)} has precisely one limit geometry given by Thm. \ref{thm_quartic_CCGPSR_curves_classification} \hyperref[eqn_qCCPSRcurves_class_b]{b)} which is associated to the choice $\D y(\gamma)<0$ in \eqref{eqn_limit_geo_defining_eqn} for the corresponding polynomial $h_{1,K}$, $K\in\left(-\frac{25}{72},U\right)$, in standard form. For $\D y(\gamma)>0$, no limit geometry exists as every integral curve of $\mathcal{V}$ with initial condition in $\{1\}\times \left(-\frac{25}{72},U\right)$ is unbounded in positive time direction. This also follows from the construction of $H$.
			
			The curve Thm. \ref{thm_incomplete_quartics} \hyperref[eqn_incomplete_qGPSRcurves_class_c]{c)} corresponds to all maximal connected quartic GPSR curves in standard form $\mathcal{H}_{L,K}$ with $(\pm L,K)^T$ being precisely the points in the image of $(\mathbf{w}(K),K)^T$, $K<-\frac{1}{12}$, which is the image of a maximal integral curve of $\mathcal{V}$, cf. \eqref{eqn_partial_w}. Hence, Thm. \ref{thm_incomplete_quartics} \hyperref[eqn_incomplete_qGPSRcurves_class_c]{c)} has precisely one limit geometry, corresponding to $\D y(\gamma)<0$ in \eqref{eqn_limit_geo_defining_eqn} for $h_{1,-\frac{25}{72}}$, given by Thm. \ref{thm_quartic_CCGPSR_curves_classification} \hyperref[eqn_qCCPSRcurves_class_b]{b)}. For $\D y(\gamma)>0$, the corresponding points in the image of $(\mathbf{w}(K),K)^T$ have diverging norm, meaning simply that the corresponding maximal integral curve of $\mathcal{V}$ is unbounded in one, and bounded in the other time direction.
			
			Lastly, the curve Thm. \ref{thm_incomplete_quartics} \hyperref[eqn_incomplete_qGPSRcurves_class_d]{d)} corresponds the maximal integral curve of $\mathcal{V}$ that has coinciding image with $(\mathbf{u}_{\mathrm{upper}}(K),K)^T$, $K>\frac{1}{4}$, cf. the initial discussion in the proof of Theorem \ref{thm_incomplete_quartics}. With a similar argument as for the above case Thm. \ref{thm_incomplete_quartics} \hyperref[eqn_incomplete_qGPSRcurves_class_c]{c)} we obtain that $\mathcal{H}_{1,U}$ has precisely one well-defined limit geometry corresponding given by Thm. \ref{thm_quartic_CCGPSR_curves_classification} \hyperref[eqn_qCCPSRcurves_class_a]{a)} corresponding to $\D y(\gamma)<0$ in \eqref{eqn_limit_geo_defining_eqn} for $h_{1,U}$. We see that $\D y(\gamma)>0$ is associated to studying the limit behaviour in positive time direction of the maximal integral curve of $\mathcal{V}$ with initial condition $(1,U)^T$, which is divergent. Hence, $\mathcal{H}_{1,U}$ does not have a limit geometry for $\D y(\gamma)>0$.
		\end{proof}
	\end{Prop}
	
	\begin{rem}
			Another approach for the classification of hyperbolic homogeneous quartic polynomials $h=h\left(\left(\begin{smallmatrix}x\\y\end{smallmatrix}\right)\right)$ in two real variables is as follows. When studying the complete cases, that is when $\{h=1\}\cap\{\text{hyperbolic points of }h\}$ has at least one closed connected component $\mathcal{H}$, we know by the precompactness of intersections of the cone $\mathbb{R}_{>0}\cdot\mathcal{H}$ with affinely embedded tangent spaces of $\mathcal{H}$ that, when written in standard form, $h|_{x=1}$ mus have at least two real zeros in $y$. Hence, $\{h=0\}$ contains at least two distinct lines, so we can assume without loss of generality after a suitable coordinate transformation that $h$ is of the form
				\begin{equation*}
					h=xy(\varepsilon_x x^2 + cxy+ \varepsilon_y y^2),
				\end{equation*}
			where $\varepsilon_x,\varepsilon_y\in\{0,1\}$ and $c\in\mathbb{R}$. Note that we obtain the condition on $\varepsilon_x$ and $\varepsilon_y$ by assuming, again without the loss of generality, that $h$ is positive at $\left(\begin{smallmatrix}1\\ \delta\end{smallmatrix}\right)$ and at $\left(\begin{smallmatrix}\delta\\1 \end{smallmatrix}\right)$ for $\delta>0$ small enough.
					
			The upside of this approach is that we can comparatively quickly obtain a partial classification which in particular contains all quartic CCGPSR curves, but we are stuck at the point were we need to classify irreducible homogeneous quartic polynomials on $\mathbb{R}^2$ that are positive everywhere except at the origin, and the one equivalent to Thm. \ref{thm_incomplete_quartics} \hyperref[eqn_incomplete_qGPSRcurves_class_c]{c)} which has the property that $\{h=0\}$ contains precisely one line. After this we would need to study them for hyperbolicity. To our knowledge, such a classification did not exists when writing this work. Now, however, such a classification is known as we have shown that every hyperbolic homogeneous quartic polynomial on $\mathbb{R}^2$ that vanishes only at the origin is equivalent to exactly one of Thm. \ref{thm_incomplete_quartics} \hyperref[eqn_incomplete_qGPSRcurves_class_a]{a)} and that the unique hyperbolic homogeneous quartic polynomial on $\mathbb{R}^2$ that vanishes along precisely one line is equivalent to Thm. \ref{thm_incomplete_quartics} \hyperref[eqn_incomplete_qGPSRcurves_class_c]{c)}.
			
			The downsides of this approach are twofold. Firstly, this does not allow us to determine the limit geometries of the so-obtained curves. We would need to bring our curves in standard form and check, by hand, which polynomials $P_3$ and $P_4$ we can obtain. This would lead to equations \eqref{eqn_L_T_pmoving} and \eqref{eqn_K_T_pmoving}, both of which are extremely complicated and we cannot expect to realistically use them to obtain our desired result about limit geometries. Hence, we would after all need to describe the equivalence classes in terms of maximal integral curves of $\mathcal{V}$. Furthermore as explained above, our method allows us to classify all maximal incomplete connected quartic GPSR curves in addition to the closed ones. Additionally, interpreting equivalence classes of quartics as maximal integral curves of $\mathcal{V}$ allows us to better understand the moduli space of hyperbolic homogeneous quartic polynomials, at least in two real variables, that is $\mathrm{Sym}^4(\mathbb{R}^2)^*_{\mathrm{hyp.}}/\mathrm{GL}(2)$. This is mainly motivated to obtain a comparison with the cubic case, cf. \cite{L2}. The moduli space $\mathrm{Sym}^3(\mathbb{R}^{n+1})^*_{\mathrm{hyp.}}/\mathrm{GL}(n+1)$ admits a compact convex generating set for any $n\in\mathbb{N}_0$. But in the quartic case we have now seen that this is not true, even for $n=1$. This is one of the big difficulties in proving, or disproving, that closed quartic GPSR curves are complete with respect to the centro-affine fundamental form, while in the cubic case there are three different known proofs \cite[Thm.\,2.5]{CNS} \cite[Prop.\,4.17]{L1} \cite[Prop.\,5.3]{L2}. However, from Theorem \ref{thm_quartic_CCGPSR_curves_classification} it follows that the curves that are singular at infinity still are compactly generated, in fact by three points. It might still hold true for higher dimensions that the moduli space of closed connected quartic GPSR manifolds with non-regular boundary behaviour is also compactly generated, which would imply the completeness result by using \cite[Thm.\,1.18]{CNS}.
	\end{rem}
%\newpage


\begin{thebibliography}{ABCD}
	\bibitem[ACD]{ACD}  D. V. Alekseevsky, V. Cort\'es, and C. Devchand, \textit{Special complex manifolds}, J. Geom. Phys. {\bf 42} (2002), No. 1--2, 85--105.
	
	\bibitem[C]{C} V. Cort\'es, {\it Alekseevskian spaces}, Differential Geom. Appl. {\bf 6} (1996), No. 2, 129--168.
	
	\bibitem[CDL]{CDL} V. Cort\'es, M. Dyckmanns, and D. Lindemann, \textit{Classification of complete projective special real surfaces}, Proc. London Math. Soc. \textbf{109} (2014), No. 2, 423--445. 
	
	\bibitem[CDJL]{CDJL} V. Cort\'es, M. Dyckmanns, M. J\"ungling, and D. Lindemann, \textit{A class of cubic hypersurfaces and quaternionic K\"ahler manifolds of co-homogeneity one} (2017), to appear in Asian J. Math.,  \href{https://arxiv.org/abs/1701.07882}{arxiv:1701.07882}.
	
	\bibitem[CHM]{CHM} V. Cort\'es, X. Han, and T. Mohaupt, {\it Completeness in supergravity constructions}, Commun. Math. Phys. {\bf 311} (2012), No. 1, 191--213. 
	
	\bibitem[CNS]{CNS} V. Cort\'es, M. Nardmann, and S. Suhr, {\it Completeness of hyperbolic centroaffine hypersurfaces}, Comm. Anal. Geom., Vol. \textbf{24}, No. 1 (2016), 59--92.
	
	\bibitem[CST]{CST} V. Cort\'es, D. Thung, A. Saha, \textit{Curvature of quaternionic Kähler manifolds with $S^1$-symmetry}, \href{https://arxiv.org/abs/2001.10032}{arxiv:2001.10032}.
	
	\bibitem[DV]{DV} B. de Wit, A. Van Proeyen, {\it Special geometry, cubic polynomials and homogeneous quaternionic spaces}, Comm. Math. Phys. {\bf 149} (1992), No. 2, 307--333.
	
	\bibitem[DP]{DP} J.-P. Demailly and M. Paun, \textit{Numerical characterization of the K\"ahler cone of a compact K\"ahler manifold}, Annals of Mathematics, Vol. \textbf{159} (2004), No. 3.
	
	\bibitem[F]{F} D. S. Freed, \textit{Special K\"ahler manifolds}, Comm. Math. Phys. \textbf{203} (1999), No. 1, 31--52.
	
	\bibitem[GST]{GST} M. G{\"u}naydin, G. Sierra, and P. K. Townsend, \textit{The geometry of $N=2$ Maxwell--Einstein supergravity and Jordan algebras}, Nucl. Phys. {\bf B242} (1984), 244--268.
	
	\bibitem[KW]{KW} A.B. Korchagin and D.A. Weinberg, \textit{The isotopy classification of affine quartic curves}, Rocky Mt. J. Math., Vol. \textbf{32} (2002), No. 1, 255–-347. 
	
	\bibitem[L1]{L1} D. Lindemann, \textit{Structure of the class of projective special real manifolds and their generalisations}, PhD-thesis (2018).
	
	\bibitem[L2]{L2} D. Lindemann, \textit{Properties of the moduli set of complete connected projective special real manifolds}, \href{https://arxiv.org/abs/1907.06791}{arXiv:1907.06791}.
	
	\bibitem[L3]{L3} D. Lindemann, \textit{Limit geometry of complete projective special real manifolds},\hfill\textcolor{white}{.} \href{https://arxiv.org/abs/2009.12956}{arXiv:2009.12956}.
	
	\bibitem[M]{M} G.Th. Magn\'usson, \textit{Cohomological expression of the curvature of K\"ahler moduli}, \href{https://arxiv.org/abs/2004.06881}{arXiv:2004.06881}.
	
	\bibitem[M1]{M1} G. Th. Magn\'usson, \textit{The geometry of K\"ahler cones} (2012), \href{https://arxiv.org/abs/1211.6934}{arxiv:1211.6934}.
	
	\bibitem[M2]{M2} G. Th. Magn\'usson, \textit{Cohomological expression of the curvature of K\"ahler moduli
	} (2020), \href{https://arxiv.org/abs/2004.06881}{arxiv:2004.06881}.
	
	\bibitem[N]{N} I. Newton, \textit{Curves}, related entry in Lexicon Technicum Vol. \textbf{II} by John Harris (1710).
	
	\bibitem[P]{P} E.M. Prodanov, \textit{Classification of the Real Roots of the Quartic Equation and their Pythagorean Tunes}, Int. J. Appl. Comput. Math \textbf{7} (2021).
	
	\bibitem[T]{T} B. Totaro, \textit{The curvature of a Hessian metric}, Int. J. Math., \textbf{15}, 369 (2004).
	
	\bibitem[W]{W} P.M.H. Wilson, \textit{Sectional curvatures of K\"ahler moduli}, Math. Ann. \textbf{330} (2004) 631--664.
\end{thebibliography}
\end{document}